\documentclass[12pt,a4paper]{amsart}
\usepackage{amsfonts, amssymb, amsmath, amsthm}
\usepackage{graphicx} 
\usepackage{amscd}
\usepackage{color}
\usepackage{mathtools}
\usepackage[update,prepend]{epstopdf}
\usepackage{hyperref} 
\usepackage{epsfig} 
\usepackage{float}
\usepackage[margin=5mm]{caption}
\newtheorem{theorem}{Theorem}
\theoremstyle{plain}

\newtheorem{corollary}[theorem]{Corollary}
\newtheorem{lemma}[theorem]{Lemma}\newtheorem{remark}{Remark}

\newtheorem{proposition}[theorem]{Proposition}
\theoremstyle{definition}\newtheorem{example}{Example}

\numberwithin{equation}{section}

\def\eps{\varepsilon}

\def\ri{{\rm i}}

\def\util{\tilde{u}}
\def\vtil{\tilde{v}}
\def\uhat{\hat{u}}

\def\Atil{\tilde{A}}
\def\Ahat{\widehat{A}}
\def\Uvec{\vec{U}}
\def\Uvecext{\vec{U}^\text{ext}}
\def\Uvecapp{\vec{U}^\text{app}}
\def\UvecOext{\vec{U}^{0,\text{ext}}}
\def\UvecOapp{\vec{U}^{0,\text{app}}}

\def\Evec{\vec{E}}

\newcommand{\R}{\mathbb{R}}
\newcommand{\C}{\mathbb{C}}
\newcommand{\Q}{\mathbb{Q}}

\newcommand{\B}{\mathbb{B}}

\newcommand{\N}{\mathbb{N}}

\newcommand{\Z}{\mathbb{Z}}

\newcommand{\cD}{{\mathcal D}}

\newcommand{\cK}{{\mathcal K}}
\newcommand{\cL}{{\mathcal L}}

\newcommand{\cS}{{\mathcal S}}
\newcommand{\cT}{{\mathcal T}}

\newcommand{\cX}{{\mathcal X}}

\renewcommand{\phi}{\varphi}

\newcommand{\pa}{\partial}

\newcommand{\Res}{\text{\rm Res}}

\newcommand{\supp}{\text{\rm supp}}

\newcommand{\bspm}{\left(\begin{smallmatrix}}\newcommand{\espm}{\end{smallmatrix}\right)}
\newcommand{\bpm}{\begin{pmatrix}}\newcommand{\epm}{\end{pmatrix}}

\def\blem{\begin{lemma}}\def\elem{\end{lemma}}
\def\bthm{\begin{theorem}}\def\ethm{\end{theorem}}
\def\bprop{\begin{proposition}}\def\eprop{\end{proposition}}
\def\bcor{\begin{corollary}}\def\ecor{\end{corollary}}
\def\beq{\begin{equation}}\def\eeq{\end{equation}}
\def\bpf{\begin{proof}}\def\epf{\end{proof}}
\def\bex{\begin{example}}\def\eex{\end{example}}

\DeclareFontEncoding{FMS}{}{}
\DeclareFontSubstitution{FMS}{futm}{m}{n}
\DeclareFontEncoding{FMX}{}{}
\DeclareFontSubstitution{FMX}{futm}{m}{n}
\DeclareSymbolFont{fouriersymbols}{FMS}{futm}{m}{n}
\DeclareSymbolFont{fourierlargesymbols}{FMX}{futm}{m}{n}
\DeclareMathDelimiter{\VERT}{\mathord}{fouriersymbols}{152}{fourierlargesymbols}{147}

\def\bi{\begin{itemize}}\def\ei{\end{itemize}}

\def\ri{{\rm i}}

\def\brem{\begin{remark}}\def\erem{\end{remark}}

\def\uapp{u_{\text{app}}}

\begin{document}
\title[Coupled Mode Asymptotics for Localized Wavepackets]{Justification of the Coupled Mode Asymptotics for Localized Wavepackets in the Periodic Nonlinear Schr\"odinger Equation} 


\author{Tom\'{a}\v{s} Dohnal}
\address{T. Dohnal \hfill\break 
Department of Mathematics, Technical University Dortmund\hfill\break
D-44221 Dortmund, Germany}
\email{tomas.dohnal@math.tu-dortmund.de}

\author{Lisa Helfmeier}
\address{L. Helfmeier \hfill\break 
Department of Mathematics, Technical University Dortmund\hfill\break
D-44221 Dortmund, Germany}
\email{lisa.helfmeier@math.tu-dortmund.de}

\date{\today}


\begin{abstract} We consider wavepackets composed of two modulated carrier Bloch waves with opposite group velocities in the one dimensional periodic Nonlinear Schr\"odinger/Gross-Pitaevskii equation. These can be approximated by first order coupled mode equations (CMEs) for the two slowly varying envelopes. Under a suitably selected periodic perturbation of the periodic structure the CMEs possess a spectral gap of the corresponding spatial operator and allow families of exponentially localized solitary waves parametrized by velocity. This leads to a family of approximate solitary waves in the periodic nonlinear Schr\"odinger equation. Besides a formal derivation of the CMEs a rigorous justification of the approximation and an error estimate in the supremum norm are provided. Several numerical tests corroborate the analysis.
\end{abstract}

\maketitle
\section{Introduction}

Propagation of localized wavepackets in nonlinear media with a periodic structure is a classical problem in the field of nonlinear dispersive equations. This paper contributes to the mathematical analysis of asymptotic approximations of small wavepackets in periodic structures of finite (rather than infinitesimal) contrast. Typical examples of coherent wavepackets in nonlinear periodic media are optical pulses in nonlinear photonic crystals \cite{SE03} or matter waves in Bose-Einstein condensates with superimposed optical lattices \cite{LOKS03}. In optics such pulses can be applied as bit carriers in future devices for optical logic and computation, which are likely to heavily exploit photonic crystals \cite{Brod_98,Pelusi_etal_08}.

Small spatially broad wavepackets in nonlinear problems with periodic coefficients can be effectively studied using a slowly varying envelope approximation. The scaling of the envelope variables is, however, not unique and can lead to different effective amplitude equations with qualitatively different solutions. On the one hand there is the nonlinear Schr\"odinger scaling, in which the envelope multiplies a single selected Bloch wave and depends on the moving frame variable $x-c_g t$ and slowly on time to describe a slow temporal modulation, see e.g. \cite{BSTU06}. For semilinear equations of second order in space with a nonlinearity that is quadratic or cubic near zero (like, e.g., a nonlinear wave equation with the nonlinearity $u^2$ or $u^3$ or the periodic nonlinear Schr\"odinger equation (PNLS) with the nonlinearity $|u|^2u$) the ansatz
\beq\label{E:NLS_ans}
u(x,t)\sim \eps A(\eps(x-c_g t),\eps^2 t)p(x,k_0)^{\ri (k_0x-\omega_0 t)}, \quad 0<\eps \ll 1
\eeq
leads to the cubic nonlinear Schr\"odinger equation (NLS) with constant coefficients for the envelope $A(X,T)$. The function $p(x,k_0)^{\ri (k_0x-\omega_0 t)}$, with $(k_0,\omega_0)$ being a point in the graph of the dispersion relation, is the selected Bloch wave and $c_g$ is its group velocity. The localized bound state solutions of the NLS predict via \eqref{E:NLS_ans} approximate solitary wave solutions $u$ of the original problem. At a given frequency $\omega_0$ the solitary wave has a velocity that is asymptotically close to the group velocity $c_g$ of the selected Bloch wave. In one spatial dimension, where the dispersion relation is monotonous and symmetric there is only one group velocity (up to the sign) for a given frequency.

On the other hand one can consider a wavepacket composed of several carrier Bloch waves. Such coupling was studied, for instance, in \cite{GMS08} for the cubic PNLS in $d\in \N$ dimensions, where the Bloch waves are allowed to have different frequencies but need to form a so called closed mode system. When restricted to $d=1$ and a single frequency, the ansatz has the form
\beq\label{E:CME_ans}
u(x,t)\sim \eps^{1/2} e^{-\ri \omega_0 t}\left( A_+(\eps x,\eps t)p_+(x)e^{\ri k_0x}+A_-(\eps x,\eps t)p_-(x)e^{-\ri k_0x}\right), \quad 0<\eps \ll 1,
\eeq
 where $p_\pm(x)e^{\ri (\pm k_0x-\omega_0 t)}$ are the two Bloch waves. Clearly, the two envelopes $A_\pm$ are not prescribed to be functions of any moving frame variable. The scaling in \eqref{E:CME_ans} leads to a system of first order Dirac type equations - so called \textit{coupled mode equations} (CMEs). In order for the ansatz to predict approximate solitary waves of the original system, e.g. the PNLS, the envelope pair $(A_+,A_-)$ has to be a solitary wave solution of the CMEs with $A_+$ and $A_-$ propagating at the same velocity. The CMEs in \cite{GMS08} do not possess a spectral gap of the corresponding linear spatial operator such that (exponentially) localized solitary waves of the CMEs are not expected. Indeed, the linear part of the CMEs in  \cite{GMS08} is 
$$
\begin{aligned}
&\ri(\pa_T A_++c_g \pa_XA_+)=0\\
&\ri(\pa_T A_--c_g \pa_XA_-)=0
\end{aligned}
$$
and the spectrum of $\ri c_g \bspm \pa_X & 0 \\ 0 & -\pa_X\espm$ is the whole $\R$. Note that $X:=\eps x$ and $T:=\eps t$.

Our aim is to find a setting which leads to CMEs with solitary waves and to rigorously justify the asymptotic approximation. If a family of CME-solitary waves parametrized by velocity exists, then ansatz \eqref{E:CME_ans} predicts a family of approximate solitary waves with frequency close to $\omega_0$ but with an $O(1)$- range of velocities, i.e. not only velocities asymptotically close to the group velocity of the Bloch wave $p_+$ or $p_-$.

We restrict to the one dimensional periodic PNLS and show that a large class of $\eps$-small periodic perturbations of the underlying periodic structure leads to CMEs with the linear part
$$
\begin{aligned}
&\ri(\pa_T A_++c_g \pa_XA_+)+\kappa A_-=0\\
&\ri(\pa_T A_--c_g \pa_XA_-)+\kappa A_+=0
\end{aligned}
$$
with $\kappa>0$, where, clearly, the operator $\bspm \ri c_g \pa_X & \kappa \\ \kappa & -\ri c_g \pa_X\espm$ has the spectral gap $(-\kappa,\kappa)$. In certain cases the CMEs are identical to those for small contrast periodic structures \cite{AW89,GWH01,SU01,Pelinov_2011}, which possess explicit solitary wave families, see \cite{AW89} and Section \ref{S:num} here.

\bigskip 

Hence, we consider 
\beq\label{E:PNLS}
\ri \pa_t u + \pa_x^2 u - (V(x) + \eps W(x))u -\sigma(x) |u|^2 u =0, \qquad x,t \in \R,
\eeq
where the real functions $V,W$, and $\sigma$ satisfy $V(x+2\pi)=V(x), \sigma(x+2\pi)=\sigma(x)$ and $W(x+2\pi/k_W)=W(x)$ for all $x\in \R$ with $k_W>0$ and where $W$ is of mean zero, i.e. 
\beq\label{E:W}
W(x)=\sum_{n\in \Z\setminus \{0\}}a_n e^{\ri n k_W x}, \quad a_{-n}=\overline{a_n} \text{ for all } n\in \N.
\eeq
We assume (see assumptions (H1)-(H4)) that $V$ and $\sigma$ are continuous, $W\in C^2(\R)$ and that if $k_W\notin \Q$, then $W$ has only finitely many nonzero Fourier components, i.e. there is an $M\in \N$ such that $a_n=0$ for all $|n|>M$. Note that there is no loss of generality in assuming $a_0=0$ since solutions for a nonzero $a_0\in \R$ can be obtained by a simple phase rotation factor $e^{-\ri \eps a_0 t}$. We point out that since $V$ is independent of $\eps$, the periodic structure has finite contrast, unlike in \cite{AW89,GWH01,SU01,Pelinov_2011}, where infinitesimal contrast is considered.

Our approximate wavepacket ansatz is
\beq\label{E:uapp}
\uapp(x,t) := \eps^{1/2}e^{-\ri \omega_0 t}\left(A_+(\eps x,\eps t)p_+(x)e^{\ri k_+ x} + A_-(\eps x,\eps t)p_-(x)e^{\ri k_- x}\right),
\eeq
where $p_\pm$ are $2\pi$-periodic Bloch eigenfunctions at the ``band structure coordinates'' $(k_\pm,\omega_0)$ respectively, see Sec. \ref{S:Bloch} for details.

This problem was previously studied also in \cite{D14}, where a specific form of $V$ was considered, namely a finite band potential. This choice guarantees the presence of points $(k_0,\omega_0)$
in the band structure, where at $\eps=0$ eigenvalue curves cross transversally and the perturbation $\eps W$ generates a small spectral gap of the operator $-\pa_x^2+V+\eps W$. Here we show that this special choice of $V$ is not necessary for CMEs with a spectral gap. In addition we provide a rigorous justification of the CME-approximation by proving an estimate of the supremum norm of the error $u-\uapp$ on time intervals of length $O(\eps^{-1})$.
In \cite{GMS08}, where the $d-$dimensional problem with $W\equiv 0$ is considered, such a justification is performed in $H^s(\R^d)$ using mainly the semigroup theory, the unitary nature of the Gross-Pitaevskii group and the Gronwall inequality. The proof is performed in the $(x,t)$-variables. On the other hand the approach in \cite{BSTU06} (for a nonlinear wave equation and the ansatz \eqref{E:NLS_ans}) is based on reformulating the problem into one for the Bloch coefficients. First the Bloch transform is applied to the equation as well as the ansatz. The transform is an isomorphism between $H^s(\R)$ and the $L^2(\B,H^s(0,2\pi))$-space over the Brillouin zone $\B$. Subsequently one expands the solution and the ansatz in the Bloch eigenfunctions. This leads to an infinite dimensional ODE system for the expansion coefficients parametrized by the wavenumber $k$. Due to the concentration of the ansatz near $k=k_0$ the problem can be approximately decomposed into one on a neighborhood of $k_0$ and one on neighborhoods of the $k-$points generated by the nonlinearity applied to the ansatz. The problem near $k_0$ leads to the amplitude equation, i.e. the NLS, and the problem near the generated $k-$points can be solved explicitly. The error estimate is provided also here using a Gronwall argument. We present a relatively detailed application of the Bloch transform approach to \eqref{E:PNLS} and \eqref{E:uapp}, where besides the different scaling compared to \eqref{E:NLS_ans} in \cite{BSTU06} a major change is that the action of the potential $\eps W$ on the ansatz generates new $k-$concentration points, which need to be accounted for. Unlike \cite{BSTU06} we are forced to work in $L^1(\B,H^s(0,2\pi))$ as too many powers of $\eps$ are lost in $L^2$. This was observed also in \cite{SU01}. No isomorphism is available for the $L^1$ space and we cannot carry over estimates in the original $u(x,t)$-variable to estimates for the Bloch coefficients. Nevertheless, we take advantage of the fact that the inverse Bloch transform applied to $L^1(\B,H^s(0,2\pi))$ with $s>1/2$ produces a continuous function, see Sec. \ref{S:fn_spaces}. Hence an $L^1$ estimate of the error in the Bloch variables translates into a supremum norm estimate in the $u(x,t)$-variable.

\subsection{The Structure of the Paper}\label{S:struct}

In Section \ref{S:Bloch} the concept of Bloch waves and of the band structure is reviewed and the choice of carrier Bloch waves $p_\pm(x)e^{\ri(k_\pm x-\omega_0t)}$ for the ansatz \eqref{E:uapp} is explained. Two cases of the choice of the wavenumbers $k_+$ and $k_-$ are distinguished, namely case (a) with simple Bloch eigenvalues at $(k_+,\omega_0)$ and $(k_-,\omega_0)$ with $k_-=-k_+$, and case (b) with a double eigenvalue at $(k_+,\omega_0)=(k_-,\omega_0)$. Section \ref{S:formal} presents the formal derivation of the effective coupled mode equations and explains why it makes sense to distinguish the cases of $k_\pm$ rational and $k_\pm$ irrational. Namely, the rational case can always be reduced to case (b) with a double eigenvalue at $(0,\omega_0)$. In Section \ref{S:justif} we formulate and prove the main approximation result. After defining the function spaces and the Bloch transformation in Section \ref{S:fn_spaces}, we first present the proof for the case of rational $k_\pm$ in Section \ref{S:rational_pf}. In Section \ref{S:irrational_pf} we treat the case of irrational $k_\pm$, which works in an analogous way but several mainly notational changes are necessary. In Section \ref{S:discuss} we discuss some extensions of the analytical results. Finally, Section \ref{S:num} provides numerical examples and convergence tests confirming the analysis.

\section{Linear Bloch waves; choice of the carrier waves}\label{S:Bloch}

For a $P$-periodic ($P>0$) potential $V$, i.e. $V(x+P)=V(x)$ for all $x\in \R$ and the corresponding Brillouin zone $\B_{P}:=(-\pi/P,\pi/P]$, we consider first the Bloch eigenvalue problem
\beq\label{E:Bloch_ev_prob}
\cL(x,k)p:=-(\pa_x + \ri k)^2p +V(x)p = \omega p, \quad x\in (0,P)
\eeq
with $p(x+P)=p(x)$. There is a countable set of eigenvalues $\omega_n(k), n\in \N=\{1,2,\dots\}$ for each $k\in \B_{P}$, which we order in the natural way $\omega_n(k)\leq \omega_{n+1}(k)$. The graph $(k,\omega_n(k))_{n\in \N}$ over $k\in \B_P$ is called the band structure. As functions of $k$ the eigenvalues $\omega_n(k)$ are $2\pi/P-$periodic and analytic away from points of higher multiplicity, which can occur only at $k\in \{0,\pi/P\}$. Moreover, on $(0,\pi/P)$ the eigenvalues are strictly monotone, i.e. $\omega_n'(k)\neq 0$ for all $n\in \N, k\in(0,\pi/P)$ and the band structure is symmetric: $\omega_n(-k)=\omega_n(k)$ for all $k\in \R$. Due to the $2\pi/P$-periodicity we have also $\omega_n(\pi/P-k)=\omega_n(\pi/P+k)$ for all $k\in \R$. In addition, the multiplicity is at most two as the eigenvalue problem is an ordinary differential equation of second order. 
The $L^2$-spectrum of the operator $-\pa_x^2+V$ is $\text{spec}(-\pa_x^2+V)=\cup_{k\in \B_P} \text{spec} (\cL(\cdot,k))=\cup_{n\in \N}\omega_n(\B_{P})$. For a review of problems with periodic coefficients see \cite{Eastham} or \cite{RS4}.

The $P$-periodic eigenfunction corresponding to $\omega_n(k)$ is denoted by $p_n(x,k)$ and called a Bloch eigenfunction. After normalization the eigenfunctions satisfy 
$$\langle p_n(\cdot,k),p_m(\cdot,k)\rangle_{P}:=\langle p_n(\cdot,k),p_m(\cdot,k)\rangle_{L^2(0,P)}=\delta_{n,m} \ \text{for all} \ k\in \B_{P}.$$ 
For each fixed $k\in \B_{P}$ the set $(p_n(\cdot,k))_{n\in \N}$ is complete in $L^2(0,P)$. As functions of $k$ the eigenfunctions satisfy the periodicity 
\beq\label{E:Bloch_k_per}
p_n(x,k+2\pi/P)=p_n(x,k)e^{-\ri (2\pi/P) x}.
\eeq 
Due to the equivalence of complex conjugation of \eqref{E:Bloch_ev_prob} and replacing $k$ by $-k$, we have $p_n(\cdot,-k)=\overline{p_n(\cdot,k)}$ for simple eigenvalues $\omega_n(k)$. Hence, certainly,
$$p_n(\cdot,-k)=\overline{p_n(\cdot,k)} \qquad \text{for all } k\in(0,\pi/P).$$

For our $2\pi$-periodic $V$ in \eqref{E:PNLS} we assume that if the selected eigenvalue $\omega_0$ of \eqref{E:Bloch_ev_prob} at $k=k_0\in \{0,1/2\}$ is double, i.e. $\omega_{n_0}(k_0)=\omega_{n_0+1}(k_0)$ for some $n_0\in  \N$, then 
one can define $C^2-$smooth eigenvalue curves $\tilde{\omega}_{\pm}(k)$, see assumption (H1), as 
\begin{equation}\label{E:om_til}
\tilde{\omega}_{+}(k)=\begin{cases}\omega_{n_0}(k), &k< k_0 \\  \omega_{n_0+1}(k), &k> k_0 \end{cases} \quad \text{ and } \quad  \tilde{\omega}_{-}(k)=\begin{cases}\omega_{n_0+1}(k), &k< k_0 \\  \omega_{n_0}(k), &k> k_0 \end{cases},
\end{equation}
where $1-$periodicity of the eigenvalues is used if $k_0=1/2$. Also, we assume that the corresponding eigenfunction families
\begin{equation}\label{E:p_til}
\tilde{p}_{+}(\cdot, k)=\begin{cases}p_{n_0}(\cdot,k), &k< k_0 \\  p_{n_0+1}(\cdot,k), &k> k_0 \end{cases} \quad \text{ and } \quad  \tilde{p}_{-}(\cdot,k)=\begin{cases}p_{n_0+1}(\cdot,k), &k< k_0 \\  p_{n_0}(\cdot,k), &k> k_0 \end{cases}
\end{equation}
are Lipschitz continuous in $k$ in the $H^2(0,2\pi)$-norm, i.e. the maps $\phi_\pm:\overline{\B}\to H^2(0,2\pi), k\mapsto \tilde{p}_\pm(\cdot,k)$ are Lipschitz continuous on $\overline{\B}$, see assumption (H1). Due to the above even symmetries of the eigenvalues each such point $k_0\in \{0,1/2\}$ of double multiplicity then satisfies $\tilde{\omega}'_+(k_0)=-\tilde{\omega}'_{-}(k_0)>0.$ The conjugation and the periodicity symmetries in $k$ then imply
$$\tilde{p}_{+}(x,0)=\overline{\tilde{p}_{-}(x,0)} \quad \text{and} \quad \tilde{p}_{+}(x,1/2)=\overline{\tilde{p}_{-}(x,1/2)}e^{-\ri x}.$$

\bigskip
In \eqref{E:uapp} we first choose $\omega_0\in \text{spec}(-\pa_x^2+V)$ such that there are two linearly independent eigenfunctions of \eqref{E:Bloch_ev_prob} at $\omega=\omega_0$ and we denote by $k_+,k_-$ the corresponding values in the level set (within $\B_{2\pi}$) of $\omega_0$. Due to the band structure symmetry $\omega_n(-k)=\omega_n(k)$ and monotonicity $\omega_n'(k)\neq 0$ for all $k\in(0,1/2)$, we get that only the following two cases are possible
\begin{enumerate}
\item[(a)] \textit{simple eigenvalues at $k_+,k_-=-k_+$:} \quad  $k_+ \in (0,1/2)$, $\omega_0=\omega_{n_0}(k_+)=\omega_{n_0}(k_-)$ for some $n_0\in \N$,
\item[(b)] \textit{double eigenvalue at $k_+=k_-$:} \quad $k_+=k_- \in \{0,1/2\}$, $\omega_0=\omega_{n_0}(k_+)=\omega_{n_0+1}(k_+)$ for some $n_0\in \N$.
\end{enumerate}
In both cases we denote by $\tilde{\omega}_{+}(k), \tilde{\omega}_{-}(k)$ the eigenvalue curves with $C^2$-smoothness at $k=k_+$ and $k=k_-$ respectively with $\tilde{\omega}_{+}(k_+)=\tilde{\omega}_{-}(k_-)=\omega_0$. The corresponding eigenfunction families are denoted (in both cases) by $\tilde{p}_\pm(\cdot,k)$. The group velocity at $k=k_\pm$ is given by
\beq\label{E:cg}
c_g:=\tilde{\omega}'_{+}(k_+)=-\tilde{\omega}'_{-}(k_-)=\lim_{k\to k_+}\omega'_{n_0}(k)=-\lim_{k\to k_-}\omega'_{n_0}(k).
\eeq
To simplify the notation, we also define the Bloch eigenfunctions at $k=k_\pm$ of the families $\tilde{p}_\pm(\cdot,k)$ by $p_\pm$, i.e.
$$p_+:=\tilde{p}_+(\cdot,k_+), \quad p_-\coloneqq\tilde{p}_-(\cdot,k_-).$$

In summary, we have for the two above cases
\begin{enumerate}
\item[(a)] $\tilde{\omega}_+\equiv \tilde{\omega}_-\equiv \omega_{n_0}$; $\tilde{p}_+\equiv \tilde{p}_- \equiv p_{n_0}$, $p_-=\overline{p_+}$,
\item[(b)] $\tilde{\omega}_\pm$ given by \eqref{E:om_til}; $\tilde{p}_\pm$ given by \eqref{E:p_til}, $p_-=\overline{p_+}e^{-2\ri k_+ \cdot}$, $\langle p_+,p_-\rangle_{2\pi}=0, c_g>0$. 
\end{enumerate}
For later use we note that differentiating the eigenvalue problem \eqref{E:Bloch_ev_prob} in $k$, we obtain also the formula
\beq\label{E:cg_int}
c_g=2\langle k_+p_+ -\ri \pa_x p_+, p_+\rangle_{2\pi} = -2\langle k_-p_- -\ri \pa_x p_-, p_-\rangle_{2\pi},
\eeq
where the second equality can be obtained using the above symmetries between $p_+$ and $p_-$.

\section{Formal Asymptotics}\label{S:formal}

Substituting the formal ansatz $\uapp + \eps^{3/2}u_1(x,\eps x, \eps t)e^{-\ri \omega_0 t}$ in \eqref{E:PNLS} and collecting terms on the left hand side with the same power of $\eps$, we get at $O(\eps^{1/2})$ 
$$e^{-\ri \omega_0 t}\sum_\pm A_\pm(X,T)(\omega_0 p_\pm +(\pa_x+\ri k_\pm)^2 p_\pm -V(x)p_\pm)e^{\ri k_\pm x}, \quad X:=\eps x, T:= \eps t$$
where the expression in the parentheses vanishes for each $\pm$ due to the choice of $p_\pm$ and $\omega_0$. For $O(\eps^{3/2})$ let us first consider the term $W\uapp$. In order to identify terms of the form of a $2\pi$-periodic function times $e^{\ri k_\pm x}$, we note that we can write
\beq\label{E:W_split}
W(x)=W^{(1)}(x)+W_\pm^{(2)}(x)e^{-2\ri k_\pm x} +W^{(3)}_\pm(x),
\eeq
where 
\begin{align*}
& W^{(1)}(x)=\sum_{\stackrel{n\in \Z\setminus\{0\}}{nk_W\in \Z}}a_n e^{\ri nk_Wx}, \quad W^{(2)}_\pm(x)=\sum_{\stackrel{n\in \Z\setminus\{0\}}{nk_W\notin \Z, nk_W+2k_\pm\in \Z}}a_ne^{\ri (nk_W+2k_\pm)x},\\
 & W^{(3)}_\pm(x)=\sum_{\stackrel{n\in \Z\setminus\{0\}}{nk_W\notin \Z, nk_W+2k_\pm\notin \Z}}a_ne^{\ri nk_W x}.
\end{align*}
Clearly, $W^{(1)}$ and $W^{(2)}_\pm$ are $2\pi$-periodic while $W^{(3)}_\pm$ are not $2\pi$-periodic. Equation \eqref{E:W_split} defines two ways of splitting $W(x)$. The splitting of $W$ in \eqref{E:W_split} is motivated by the relation
$$
W(x)e^{\ri k_\pm x} = W^{(1)}(x)e^{\ri k_\pm x}+W^{(2)}_\pm(x)e^{\ri k_\mp x}+W^{(3)}_\pm(x)e^{\ri k_\pm x}.
$$
Hence, in the case $k_+=-k_-, k_+\in(0,1/2)$ the part $W_\pm^{(2)}(x)e^{-2\ri k_\pm x}$ guarantees coupling of the two carrier waves in \eqref{E:uapp} because its multiplication with $p_\pm(x)e^{\ri k_\pm x}$ produces a periodic function times $e^{\ri k_\mp x}$. Therefore, in \eqref{E:eps3/2} below we use the splitting with $W^{(2)}_+$ and $W^{(3)}_+$ for the multiplication of $W$ with $e^{\ri k_+ x}$ and the  splitting with $W^{(2)}_-$ and $W^{(3)}_-$ for the multiplication of $W$ with $e^{\ri k_- x}$. In the case $k_+=k_-\in \{0,1/2\}$, clearly, $W^{(2)}_\pm\equiv 0$ and the coupling is provided by $W^{(1)}$.
The coupling can be seen explicitly in the $\kappa$ coefficient in \eqref{E:CME}. For the above two cases we get
\begin{enumerate}
\item[(a)] $k_+=-k_-\in(0,1/2) \ \Rightarrow \ W^{(2)}_-\equiv \overline{W^{(2)}_+}$
\item[(b)] $k_+=k_-\in \{0,1/2\} \ \Rightarrow \ W^{(2)}_+\equiv W^{(2)}_-\equiv 0$.
\end{enumerate}
Because $a_{-n}=\overline{a_n},$ we also have that $W^{(1)}$ is real. 

The Bloch functions $p_\pm$ have a free complex phase, which we fix by requiring $p_-=\overline{p_+}$ and 
\begin{equation}\label{E:kap_real}
\begin{aligned}
&\text{(a)} \quad \text{Im}(\langle W^{(2)}_+p_+,p_-\rangle_{2\pi})\stackrel{!}{=}0 \text{ if } k_+=-k_-\in(0,1/2),\\
&\text{(b)} \quad \text{Im}(\langle W^{(1)}p_+,p_-\rangle_{2\pi})\stackrel{!}{=}0 \text{ if } k_+=k_-\in \{0,1/2\}.
\end{aligned}
\end{equation}
This choice makes the coefficient $\kappa$ in the effective equations \eqref{E:CME} for the envelopes $A_\pm$ real.

Hence at $O(\eps^{3/2})$ the left hand side of \eqref{E:PNLS} is (after multiplication by $e^{\ri \omega_0 t}$)
\beq\label{E:eps3/2}
\begin{aligned}
&(\omega_0+\pa_x^2 -V(x))u_1(x,X,T)\\
+&\left[\ri\left(p_+\pa_T A_+ + 2(k_+p_+-\ri p_+')\pa_XA_+\right)-W^{(1)}(x)p_+A_+-W^{(2)}_-(x)p_-A_-\right.\\
&\quad \left. -\sigma(x)\left(|p_+|^2|A_+|^2+2|p_-|^2|A_-|^2\right)p_+A_+\right]e^{\ri k_+ x }\\
+&\left[\ri\left(p_-\pa_T A_- + 2(k_-p_--\ri p_-')\pa_XA_-\right)-W^{(1)}(x)p_-A_--W^{(2)}_+(x)p_+A_+\right.\\
&\quad \left. -\sigma(x)\left(|p_-|^2|A_-|^2+2|p_+|^2|A_+|^2\right)p_-A_-\right]e^{\ri k_- x }\\
-& \sigma(x) \left(p_+^2\overline{p_-}A_+^2\overline{A_-}e^{\ri (2k_+-k_-) x }+p_-^2\overline{p_+}A_-^2\overline{A_+}e^{\ri(2k_--k_+) x}\right)\\
-& W^{(3)}_+(x)p_+A_+e^{\ri k_+ x} - W^{(3)}_-(x)p_-A_-e^{\ri k_- x},
\end{aligned}
\eeq
where $X=\eps x, T=\eps t$ and except for the last line all factors multiplying $e^{\ri k_\pm x}$ and $e^{\ri(2k_\mp-k_\pm) x}$ are $2\pi$-periodic in $x$.
Note that when $2k_+-k_- \in \{k_+,k_-\}+\Z$, then $e^{\ri(2k_+-k_-) x}$ can also be written as a $2\pi$-periodic function times $e^{\ri k_+ x}$ or $e^{\ri k_- x}$. Similarly for $e^{\ri(2k_--k_+) x}$. Within our allowed setting, i.e. within cases (a) and (b), this occurs if and only if $k_+=-k_-=1/4$ or $k_+=k_-\in\{0,1/2\}$:
\begin{align*}
k_+=-k_-=1/4 \ \Rightarrow \ &2k_+-k_-=k_-+1, \text{ i.e. }  e^{\ri(2k_+-k_-) x}= e^{\ri x}e^{\ri k_-x} \quad \text{and} \\ 
&2k_--k_+=k_+-1, \text{ i.e. }  e^{\ri(2k_--k_+) x}= e^{-\ri x}e^{\ri k_+x},\\
k_+=k_-\in\{0,1/2\} \ \Rightarrow \ &2k_+-k_-=2k_--k_+=k_+=k_-.
\end{align*}

In order to set the $O(\eps^{3/2})$ terms proportional to a $2\pi$-periodic function times $e^{\ri k_+ x}$ or $e^{\ri k_- x}$ to zero, we search for $u_1$ in form 
$$u_1(x,X,T)=U_{1,+}(X,T)s_+(x)e^{\ri k_+ x}+U_{1,-}(X,T)s_-(x)e^{\ri k_- x}$$
with $2\pi-$periodic functions $s_\pm$, such that these terms vanish. Formally, such $u_1$ exists if the Fredholm alternative holds, i.e. if the inhomogeneous terms (independent of $u_1$) in \eqref{E:eps3/2} having the form of a $2\pi$-periodic function times $e^{\ri k_\pm x}$ are $L^2(0,2\pi)$-orthogonal to $p_\pm(x)e^{\ri k_\pm x}$ respectively. In the case (b), when $k_+=k_-$, it means, of course, that all inhomogeneous terms having the form of a $2\pi$-periodic function times $e^{\ri k_+ x}=e^{\ri k_- x}$ must be orthogonal to both $p_+(x)e^{\ri k_+ x}$ and $p_-(x)e^{\ri k_- x}$. In this case we have $\langle p_+,p_-\rangle_{2\pi} =0$ and $\langle \pa_x p_+,p_-\rangle_{2\pi} =0$.  The latter identity follows for $k_+=0$ from $p_-=\overline{p_+}$ because $\langle \pa_x p_+,p_-\rangle_{2\pi} =\int_0^{2\pi}p_+(x)\pa_x p_+(x)dx=\tfrac{1}{2}\int_0^{2\pi}\pa_x (p_+(x))^2dx=0$ and for $k_+=1/2$ from $p_-=\overline{p_+}e^{-\ri x}$ because $\langle \pa_x p_+,p_-\rangle_{2\pi} =\tfrac{1}{2}\int_0^{2\pi}\pa_x (p_+(x))^2e^{\ri x}dx = -\tfrac{\ri}{2}\int_0^{2\pi}p_+^2(x)e^{\ri x}dx$ and  $\int_0^{2\pi}p_+^2(x)e^{\ri x}dx=\langle p_+,p_-\rangle_{2\pi} =0$.
Using these identities, the orthogonality conditions become the \textit{coupled mode equations} (CMEs)
\beq\label{E:CME}
 \begin{split}
  \ri \left(\pa_T +c_g \pa_X\right)A_+ + \kappa A_- + \kappa_s A_+ +\alpha (|A_+|^2+2|A_-|^2)A_+ & \\
	   + \beta (|A_-|^2+2|A_+|^2)A_- + \overline{\beta}A_+^2\overline{A_-} +\gamma A_-^2 \overline{A_+} &= 0,\\
  \ri \left(\pa_T -c_g \pa_X\right)A_- + \kappa A_+ + \kappa_s A_- +\alpha (|A_-|^2+2|A_+|^2)A_- & \\
	    + \overline{\beta} (|A_+|^2+2|A_-|^2)A_+ + \beta A_-^2\overline{A_+} +\overline{\gamma} A_+^2 \overline{A_-} &= 0,\\
 \end{split}
\eeq
where 
\begin{equation}\label{E:CME_coeffs}
\begin{aligned}
c_g & =2\langle k_+p_+ -\ri \pa_xp_+, p_+\rangle_{2\pi}= - 2\langle k_-p_- -\ri \pa_xp_-, p_-\rangle_{2\pi}\in \R,\\
\kappa_s & = - \langle W^{(1)}p_+,p_+\rangle_{2\pi} =- \langle W^{(1)}p_-,p_-\rangle_{2\pi} \in \R, \\
\kappa & = \begin{cases} -\langle W^{(2)}_-p_-,p_+\rangle_{2\pi}=-\langle W^{(2)}_+p_+,p_-\rangle_{2\pi} \in \R \quad &\text{if } k_+=-k_-\in(0,1/2),\\
  - \langle W^{(1)}p_-,p_+\rangle_{2\pi} =- \langle W^{(1)}p_+,p_-\rangle_{2\pi}\in \R \quad &\text{if } k_+=k_-\in\{0,1/2\},
\end{cases}\\
\alpha & = -\langle\sigma p_+^2,p_+^2\rangle_{2\pi}=-\langle\sigma p_-^2,p_-^2\rangle_{2\pi} = -\langle \sigma|p_-|^2p_+,p_+\rangle_{2\pi}=-\langle \sigma|p_+|^2p_-,p_-\rangle_{2\pi} \in \R,\\
\beta & = \begin{cases}  0  \quad &\text{if } k_+=-k_-\in(0,1/2),\\
 -\langle \sigma|p_\pm|^2p_-,p_+\rangle_{2\pi} = - \overline{\langle \sigma |p_\pm|^2p_+,p_-\rangle_{2\pi}}\quad &\text{if } k_+=k_-\in\{0,1/2\},
\end{cases}\\
\gamma & =  \begin{cases} 0  \quad &\text{if } k_+=-k_-\in(0,1/4)\cup(1/4,1/2),\\
-\langle \sigma p_-^2\overline{p_+}e^{-\ri \cdot},p_+\rangle_{2\pi} = -\overline{\langle \sigma p_+^2\overline{p_-}e^{\ri \cdot},p_-\rangle_{2\pi}} \quad &\text{if } k_+=-k_-=1/4,\\
 -\langle \sigma p_-^2\overline{p_+},p_+\rangle_{2\pi} = -\overline{\langle \sigma p_+^2\overline{p_-},p_-\rangle_{2\pi}} \quad &\text{if } k_+=k_-\in\{0,1/2\}.\\
\end{cases}
\end{aligned}
\end{equation}
The realness of $\kappa_s$ follows from the realness of $W^{(1)}$. Note that without any loss of generality we can set $\kappa_s=0$ because solutions $(A_+,A_-)$ of \eqref{E:CME} with $\kappa_s\neq 0$ can then be constructed from solutions $(A^{(0)}_+,A^{(0)}_-)$ with $\kappa_s=0$ via the multiplication by a simple phase factor, namely $(A_+,A_-)=(A^{(0)}_+,A^{(0)}_-)e^{\ri \kappa_s T}$.
The identities in $\kappa$ and its realness follow from \eqref{E:kap_real} and from $W^{(2)}_-=\overline{W^{(2)}_+}$ for $k_+=-k_-\in(0,1/2)$ and from the realness of $W^{(1)}$ for $k_+=k_-\in \{0,1/2\}$ .

\brem\label{R:gap}
As mentioned in the introduction, the linear part of system \eqref{E:CME} has a spectral gap if $\kappa\neq 0$. Indeed, the spectrum of the self-adjoint operator $\bspm \ri c_g \pa_x+\kappa_s & \kappa \\ \kappa & -\ri c_g \pa_x+\kappa_s\espm$ has the gap $(\kappa_s-|\kappa|,\kappa_s+|\kappa|)$. Hence, exponentially localized solutions $(A_+,A_-)$ can be expected. This is based on the heuristic argument that in spectral gaps the linear solution modes are exponentials and in the tails of the nonlinear solution, where the cubic nonlinearity is negligible, the linear dynamics govern. To the best of our knowledge, a rigorous proof of the existence of exponentially localized solitary waves of \eqref{E:CME} is not in the literature. However, for $\beta=\gamma=0$ explicit families of  exponentially localized solitary waves parametrized by velocity exist, see \eqref{E:1D_gap_sol}. The definition of $\kappa$ in \eqref{E:CME_coeffs} produces a necessary and sufficient condition for a spectral gap. In case (a), i.e. if $k_+=-k_-\in(0,1/2)$, this condition is
$$\langle W^{(2)}_+p_+,p_-\rangle_{2\pi} \neq 0$$
and in case  (b), i.e. if $k_+=k_-\in \{0,1/2\}$, it is
$$\langle W^{(1)}p_+,p_-\rangle_{2\pi}\neq 0.$$
A necessary condition is $W^{(2)}_\pm \neq 0$ in case (a) and $W^{(1)}\neq 0$ in case (b). Based on the definition of $W^{(2)}_\pm$ and  $W^{(1)}$ in \eqref{E:W_split}, we obtain in case (a) the necessary condition
\beq\label{E:kW_cond_a}
k_W \in 2k_+ + \frac{1}{n} \Z \quad \text{for some }n \in \{m\in \Z\setminus\{0\}: a_m\neq 0\}
\eeq
and in case (b) the necessary condition
\beq\label{E:kW_cond_b}
k_W\in \frac{1}{n}\Z \quad \text{for some }n \in \{m\in \Z\setminus\{0\}: a_m\neq 0\}.
\eeq
For a further discussion we note that these are clearly possible only if all $k_\pm,k_W$ are rational or all are irrational.

The simplest choice which satisfies conditions \eqref{E:kW_cond_a} and \eqref{E:kW_cond_b} is $k_W =2k_+$ with $a_1=a_{-1}=\tfrac{a}{2}\in \R\setminus\{0\}, a_n =0$ for all $n\in  \Z\setminus\{1,-1\}$, i.e. $W(x)=a\cos(2k_+x)$.
\erem

\bigskip

\subsection{The case of rational $k_+,k_-$}\label{S:rational}

We consider the case of rational $k_\pm$ separately because a common period of the Bloch waves and the potential $V$ (and $\sigma$) can be chosen in this case and the points $k_\pm$ get mapped to zero in the Brillouin zone corresponding to this new period. Hence, for $k_\pm \in \Q$ the problem is effectively transformed to the above case (b) of a double Bloch eigenvalue at $k=0$. 

\begin{figure}
\begin{center}
\scalebox{0.7}{\includegraphics{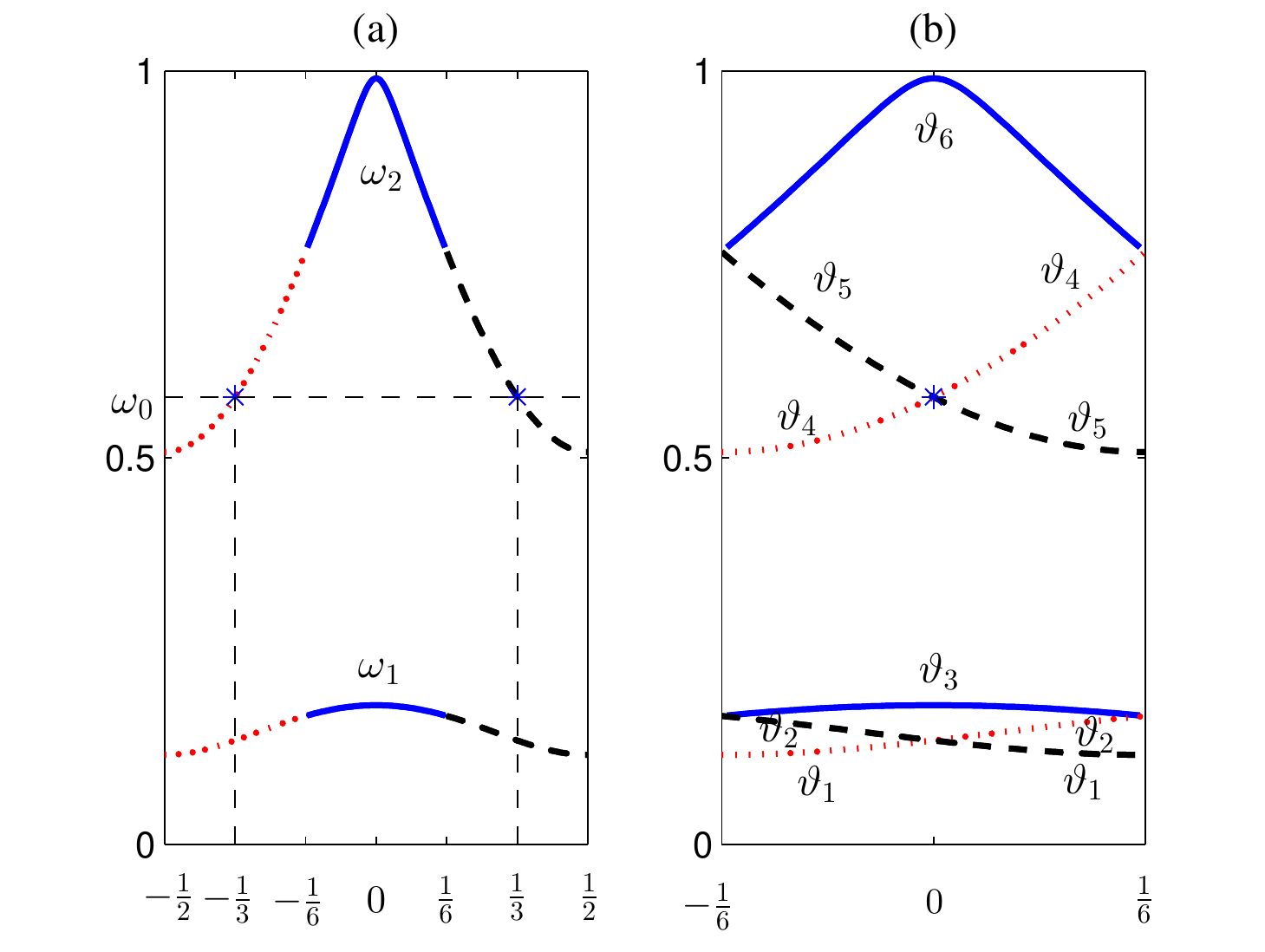}}
\end{center}
\caption{Band structure for $V(x) = \cos^5(x)$.  (a) Brillouin zone $\B_{2\pi}=(-1/2,1/2]$. The marked points are at $k_{\pm}=\pm\tfrac{1}{3}$. (b) Brillouin zone $\B_{6\pi}=(-1/6,1/6]$ corresp. to the period $6\pi$. The points $k_\pm$ are mapped to $0$ in (b). A curve segment with a given line style in (a) is mapped to that in (b). Here $n_*=4$, see \eqref{E:n_star}.} \label{F:reduce_zone}
\end{figure}

For $k_+\in \Q$ the functions $V$ and $e^{\ri k_+ \cdot}$ have a common period $[0,Q]$ with $Q=N2\pi$ for some $N\in \N$. Clearly, as either $k_-=-k_+$ or $k_-=k_+$, also $e^{\ri k_- \cdot}$  is $Q-$periodic. We use $Q$ as the working periodicity of the problem. The corresponding Brillouin zone is 
$$\B_Q:=(-\tfrac{1}{2N},\tfrac{1}{2N}].$$  
The Bloch eigenvalue problem is
\beq\label{E:Bloch_prob_Q}
\cL(x,k)q = \vartheta q, \quad x\in (0,Q)
\eeq
for $k\in\B_Q$.
The band structure on $\B_Q$ is generated from that on $\B_{2\pi}$ using the $\tfrac{1}{N}$-periodicity in the variable $k$ of the eigenvalues. The labeling of the eigenvalues changes when mapped from $\B_{2\pi}$ to $\B_Q$. We denote the band structure on $\B_Q$ by $(k, (\vartheta_n(k))_{n\in \N})$ with $k\in \B_Q$. 
Fig. \ref{F:reduce_zone} shows an example of a band structure on $\B_Q$ for a given band structure on $\B_{2\pi}$ and for $N=3,k_+=1/3$. 

Since $N k_\pm \in \Z$, we get that 
$$k_\pm = 0 \mod \tfrac{1}{N}.$$
The eigenvalue at $(k,\vartheta)=(0,\omega_0)$ is thus double. 
Because of assumption (H1) and \eqref{E:cg} the eigenvalues $\vartheta_n(k)$ can be relabeled to produce the transversal crossing of two $C^2$ eigenvalue curves at $k=0$. Our labeling of the eigenvalues is thus determined as follows. Let us denote the eigenvalues of \eqref{E:Bloch_prob_Q} at $k=0$ by  $(\lambda_n)_{n\in \N}$ ordered by size and suppose $\lambda_{n_*}=\lambda_{n_*+1}=\omega_0$. We label the eigenvalue curves $(\vartheta(k))_n$ according to size for all $n<n_*$ and all $n>n_*+1$ and for $n\in\{n_*,n_*+1\}$ we label the curves such that they are smooth at $k=0$, i.e.
\beq\label{E:n_star}
\begin{aligned}
&\vartheta_n(k) \leq \vartheta_{n+1}(k) \text{ for all } k\in \B_Q, n \in \{1,2,\dots,n_*-1,n_*+1,n_*+2,\dots\}\\
&\vartheta_{n_*}(k)< \vartheta_{n_*+1}(k) \text{ for } k<0, \ \vartheta_{n_*}(k)> \vartheta_{n_*+1}(k) \text{ for } k>0.
\end{aligned}
\eeq
By assumption (H1) is $\vartheta_{n_*}, \vartheta_{n_*+1} \in C^2(\text{int}(\B_Q))$. Note that
$$c_g= \vartheta'_{n_*}(0)=-\vartheta'_{n_*+1}(0).$$

The eigenfunctions corresponding to $\vartheta_n(k)$ are denoted by $q_n(x,k)$ and normalized via 
$$\langle q_n(\cdot,k), q_m(\cdot,k)\rangle_{Q}=\delta_{n,m}.$$ 
The $Q-$periodic eigenfunctions at $(k,\vartheta)=(0,\omega_0)$ with group velocities $\pm c_g$ are 
$$q_+(x):=q_{n_*}(\cdot,0) \quad \text{and} \quad  q_-(x):=q_{n_*+1}(\cdot,0).$$
These are related to $p_\pm$ via
$$q_\pm(x)=\frac{1}{\sqrt{N}}p_\pm(x) e^{\ri k_\pm x}, \ \text{and satisfy} \ q_-=\overline{q_+}.$$
Note that the orthogonality $\langle q_+,q_-\rangle_Q=0$ can be checked directly. It is obvious in case (b) where $k_+=k_-\in \{0,1/2\}$ and $\langle q_+,q_-\rangle_Q=\tfrac{1}{N}\langle p_+,p_-\rangle_Q=\langle p_+,p_-\rangle_{2\pi}=0$ follows from the orthogonality of eigenfunctions at each $k$. For case (a), where $k_+=-k_-\in (0,1/2)$ we argue using a Fourier series expansion. Namely, writing $p_+^2(x)=\sum_{n\in \N}b_ne^{\ri n x}$, we have
$$\langle q_+,q_-\rangle_Q=\frac{1}{N}\int_0^Q p_+^2(x)e^{2\ri k_+ x} dx=\sum_{n\in \N}b_n\int_0^{2N\pi}e^{\ri(n+2k_+) x}dx =0$$
because $n+2k_+\neq 0$ and $n+2k_+\in \tfrac{1}{N}\Z$ for all $n\in \Z$.

Because in the rational case our Bloch eigenfunctions $q_+,q_-$ are both at $k=0$, the splitting of $W$ using the parts $W^{(2)}_\pm$ as in in \eqref{E:W_split} is not suitable. It is more convenient to rewrite
Because in the rational case our Bloch eigenfunctions $q_+,q_-$ are both at $k=0$, the splitting of $W$ using the parts $W^{(2)}_\pm$ as in in \eqref{E:W_split} is not suitable. It is more convenient to rewrite
$$
W(x)=W^{(1)}_Q(x)+W_Q^{(R)}(x), \quad W^{(1)}_Q(x):=\sum_{n\in \Z_1}a_n e^{\ri nk_Wx}, \ W^{(R)}_Q(x):=\sum_{n\in \Z_R}a_n e^{\ri nk_Wx},
$$
where
$$\Z_1:= \{n\in \Z\setminus\{0\}: a_n\neq 0, nk_W\in \tfrac{1}{N}\Z\}, \quad \Z_R:= \{n\in \Z\setminus\{0\}: a_n\neq 0, nk_W\notin \tfrac{1}{N}\Z\}.$$

The formal ansatz we will use in the rational case is
\beq\label{E:uapp_rat}
\uapp(x,t) = \eps^{1/2}e^{-\ri \omega_0 t}\left(A_+(\eps x,\eps t)q_+(x) + A_-(\eps x,\eps t)q_-(x)\right).
\eeq
The effective model is again the coupled mode system \eqref{E:CME} but the nonlinear coefficients have to be modified due to the different normalization of $p_\pm$ and $q_\pm$, i.e. we have
\beq\label{E:CME_rat}
 \begin{split}
  \ri \left(\pa_T +c_g \pa_X\right)A_+ + \kappa A_- + \kappa_s A_+ +\alpha_Q (|A_+|^2+2|A_-|^2)A_+ & \\
	   + \beta_Q (|A_-|^2+2|A_+|^2)A_- + \overline{\beta_Q}A_+^2\overline{A_-} +\gamma_Q A_-^2 \overline{A_+} &= 0,\\
  \ri \left(\pa_T -c_g \pa_X\right)A_- + \kappa A_+ + \kappa_s A_- +\alpha_Q (|A_-|^2+2|A_+|^2)A_- & \\
	    + \overline{\beta_Q} (|A_+|^2+2|A_-|^2)A_+ + \beta_Q A_-^2\overline{A_+} +\overline{\gamma_Q} A_+^2 \overline{A_-} &= 0,\\
 \end{split}
\eeq
where
\begin{align*}
c_g & =2\ri \langle \pa_xq_+, q_+\rangle_Q= -2\ri \langle \pa_xq_-, q_-\rangle_Q\in \R,\\
\kappa_s & = - \langle W_Q^{(1)}q_+,q_+\rangle_Q =- \langle W_Q^{(1)}q_-,q_-\rangle_Q \in \R, \\
\kappa & = -\langle W_Q^{(1)}q_-,q_+\rangle_Q=-\langle W_Q^{(1)}q_+,q_-\rangle_Q \in \R,\\
\alpha_Q & = -\langle\sigma q_+^2,q_+^2\rangle_{Q}=-\langle\sigma q_-^2,q_-^2\rangle_{Q} = -\langle \sigma|q_-|^2q_+,q_+\rangle_Q=-\langle \sigma|q_+|^2q_-,q_-\rangle_Q \in \R,\\
\beta_Q & = -\langle \sigma|q_\pm|^2q_-,q_+\rangle_Q = -\overline{\langle \sigma |q_\pm|^2q_+,q_-\rangle_Q},\\
\gamma_Q & = -\langle \sigma q_-^2\overline{q_+},q_+\rangle_Q = -\overline{\langle \sigma q_+^2\overline{q_-},q_-\rangle_Q}.
\end{align*}
Note that the linear coefficients $c_g,\kappa_s$ and $\kappa$ are indeed identical with those defined in \eqref{E:CME_coeffs} using $p_\pm$ but for the nonlinear coefficients we have 
$$(\alpha_Q,\beta_Q,\gamma_Q)=\tfrac{1}{N}(\alpha,\beta,\gamma).$$
 These identities between the linear and nonlinear coefficients can be checked by expanding the periodic parts of the integrands in a Fourier series. E.g. to show $\langle W^{(1)}_Q q_+,q_+ \rangle_Q=\langle W^{(1)} p_+,p_+ \rangle_{2\pi}$, we expand $|q_+(x)|^2=\tfrac{1}{N}|p_+(x)|^2=\sum_{m\in \Z}b_me^{\ri m x}$ and get
\begin{align*}
\langle W^{(1)}_Q q_+,q_+ \rangle_Q&=\sum_{m\in \Z, n\in\Z_1}a_nb_m\int_0^{2N\pi} e^{\ri (m+nk_W)x}dx=\sum_{\stackrel{m\in \Z}{n\in \{j\in\Z_1:jk_W\in \Z\}}}a_nb_m\int_0^{2N\pi} e^{\ri (m+nk_W)x}dx\\
&=\frac{1}{N}\langle W^{(1)}p_+,p_+\rangle_{Q}=\langle W^{(1)} p_+,p_+ \rangle_{2\pi}
\end{align*}
because if $m\in \Z$, then $\int_0^{2N\pi} e^{\ri (m+nk_W)x}dx=0$ for any $nk_W\in \tfrac{1}{N}\Z \setminus \Z$. Similar arguments yield the other identities.


\section{Rigorous Justification of CMEs as an Effective Model}\label{S:justif}

It is clear that the discussion and the derivation of the coupled mode equations (CME) \eqref{E:CME} in Section \ref{S:formal} is only formal as terms which are not of the form of a $2\pi$-periodic function times $e^{\ri k_\pm x}$ in the residual as well as higher derivative terms with respect to $X$ have been ignored. 

To make the discussion rigorous, we will modify (and extend) the asymptotic ansatz in order to make the residual of $O(\eps^{5/2})$ in a suitable norm. After estimating the residual, we use the Gronwall lemma to show the smallness of the asymptotic error on large time scales.

We make the following basic assumptions
\begin{itemize}
\item[(H1)] In the case of a double eigenvalue at $(k,\omega)=(k_+,\omega_0)=(k_-,\omega_0)$ the functions $k\mapsto \tilde{\omega}_{\pm}(k)$ in \eqref{E:om_til} have two continuous derivatives at $k=k_+=k_-$ and the mappings $\phi_\pm:\overline{\B}\to H^2(0,2\pi), k\mapsto \tilde{p}_\pm(\cdot,k)$ with $\tilde{p}_\pm$ defined in \eqref{E:p_til} are Lipschitz continuous on $\overline{\B}$,
\item[(H2)] $W\in C^2(\R,\R), W(x+\tfrac{2\pi}{k_W})=W(x)$ for all $x\in \R$ is given by \eqref{E:W} and if $k_W\notin \Q$, then there is $M\in \N$ such that $a_n=0$ for all $|n|>M$, 
\item[(H3)] $V\in C(\R,\R), V(x+2\pi)=V(x)$ for all $x\in \R$.
\item[(H4)] $\sigma\in C(\R,\R), \sigma(x+2\pi)=\sigma(x)$ for all $x\in \R$.
\end{itemize}
Note that in (H1) the Lipschitz continuity over the whole $\overline{\B}$ requires that if $k_+=k_-=0$ and the eigenvalue at $(k,\omega)=(1/2,\tilde{\omega}_+(1/2))$ or at $(1/2,\tilde{\omega}_-(1/2))$ is double, then the eigenvalue curves need to be smoothly extendable also across $k=1/2$. In the case of simple eigenvalues at $k_+, k_-$ the $C^2$ smoothness of $\tilde{\omega}_\pm$ and the Lipschitz continuity of the Bloch eigenfunctions in (H1) always hold, see \cite{Kato-1966}.

\bthm\label{T:main}
Assume (H1)-(H4) and let $(A_+,A_-)$ be a solution of \eqref{E:CME}, \eqref{E:CME_coeffs} with $\Ahat_\pm \in C^1([0,T_0],$$L^1_{s_A}(\R)\cap L^2(\R))$ for some $T_0>0$ and some $s_A \geq 2$. There exist $c>0$ and $\eps_0>0$ such that if $u(x,0)=u_\text{app}(x,0)$ given by \eqref{E:uapp}, then the solution $u$ of \eqref{E:PNLS} satisfies $u(x,t)\to 0$ as $|x|\to \infty$ and 
$$\|u(\cdot,t)-u_\text{app}(\cdot,t)\|_{C^0_b(\R)}\leq c\eps^{3/2} \quad \text{for all } \eps\in(0,\eps_0), t\in [0,\eps^{-1}T_0].$$
\ethm
The functions $\Ahat_\pm$ are the Fourier transformations of $A_\pm$, where we define
$$ \widehat{f}(k)=\frac{1}{2\pi}\int_{-\infty}^\infty f(x)e^{-\ri kx}dx$$
with the inverse transformation $f(x)=\int_{-\infty}^\infty \widehat{f}(k)e^{\ri kx}dk$. The space $L^1_{r}, r>0$ in the theorem is defined as
$$L^1_r(\R):=\left\{f\in L^1(\R): \|f\|_{L^1_r(\Omega)}:=\int_\Omega (1+|x|)^r|f(x)|dx <\infty\right\}.$$
To the best of our knowledge, existence of solutions $(A_+,A_-)$ of \eqref{E:CME}, \eqref{E:CME_coeffs} with $\Ahat_\pm \in C^1([0,T_0],$$L^1_{s_A}(\R)\cap L^2(\R))$ is not covered in the existing literature. However, in the case $\beta=\gamma=0$ there are explicit smooth solutions (satisfying $\Ahat_\pm \in C^1([0,T_0],L^1_{s_A}(\R)\cap L^2(\R))$), see \eqref{E:1D_gap_sol}. Note that the case $\beta=\gamma=0$ is generic as it corresponds to $k_+=-k_-\in (0,1/4)\cup(1/4,1/2)$, see \eqref{E:CME_coeffs}. Note also that in general (for all values of the coefficients) local existence is guaranteed for \eqref{E:CME} in $H^1(\R)$, i.e.~for any $(A_+,A_-)(\cdot,0)\in H^1(\R)$ there is a $T_0>0$, such that \eqref{E:CME} has a $C^1((0,T_0),H^1(\R))$-solution. This follows from Theorem 1 in \cite{Reed}. Moreover, by Theorem 2 in \cite{Reed} the existence is either global (i.e.~$T_0=\infty$ can be chosen) or $\|\Ahat_+(T)\|_{H^1(\R)}+\|\Ahat_-(T)\|_{H^1(\R)}\to \infty$ as $T\to T_0$.

\bigskip
We are particularly interested in localized traveling waves $u$ of nearly constant shape. Such solutions are guaranteed if we find exact solitary waves of \eqref{E:CME}. For $\beta=\gamma=0$ system \eqref{E:CME} is exactly the classical coupled mode system for optical Bragg fibers with a small contrast, see e.g. \cite{AW89,GWH01}. If $\kappa=\langle W^{(2)}_+p_+,p_-\rangle_{2\pi}\neq 0$, this system has an explicit family of traveling solitary waves, see \cite{AW89,GWH01}. As explained above, in the generic case $k_+=-k_-\in (0,1/4)\cup(1/4,1/2)$ we always have $\beta=\gamma=0$, and hence approximate localized traveling waves $u$ of \eqref{E:PNLS} are certainly guaranteed if $k_+=-k_-\in (0,1/4)\cup(1/4,1/2)$ and $\langle W^{(2)}_+p_+,p_-\rangle_{2\pi}\neq 0$.

For $\beta,\gamma\neq 0$ traveling wave solutions may be constructed by a homotopy continuation, see \cite{D14} for a numerical implementation. As explained in Remark \ref{R:gap}, a necessary condition for exponentially localized waves of the CMEs is that either all $k_\pm, k_W$ be rational or all be irrational.
Theorem \ref{T:main}, however, holds for any $k_W\in \R$ and any $k_+=-k_-\in (0,1/2)$ or 
$k_+=k_-\in \{0,1/2\}$.  

The proof of the case of rational $k_\pm$ is technically somewhat simpler as there is a common period of $V(x)$ and $e^{\ri k_\pm x}$ and on the Brillouin zone corresponding to this common periodicity cell the wavenumbers $k_+,k_-$ both correspond (are periodic images of) the point $k=0$ as explained in Sec. \ref{S:rational}. We treat the case of rational $k_\pm$ in Sec. \ref{S:rational_pf} and irrational $k_\pm$ in Sec. \ref{S:irrational_pf}.

\subsection{Bloch Transformation, Function Spaces}\label{S:fn_spaces}
For a given $P>0$ the Bloch transformation is the operator
$$
\begin{aligned}
&\cT:H^s(\R,\C) \to L^2(\B_P,H^s((0,P),\C)), u \mapsto \util :=\cT u,\\
& \util(x,k)=(\cT u)(x,k)=\sum_{j\in \tfrac{2\pi}{P}\Z}e^{\ri j x}\uhat(k+j).
\end{aligned}
$$
By construction of $\util$ we have
\beq\label{E:util_per}
\util(x+P,k)=\util(x,k) \quad \text{and} \ \util(x,k+2\pi/P)=e^{-\ri (2\pi/P) x}\util(x,k) \quad \text{for all } x,k \in \R.
\eeq
In the following we write simply $H^s(\Omega)$ for $H^s(\Omega,\C)$. The Bloch transform $\cT:H^s(\R) \to L^2(\B_P,H^s(0,P))$ is an  isomorphism for $s\geq 0$, see \cite{RS4}, and the inverse is given by 
$$
u(x)=(\cT^{-1}\util)(x)=\int_{\B_P}e^{\ri kx}\util(x,k)dk.
$$
It is easy to see that the product of two general $H^s(\R)$ functions is mapped to a convolution by $\cT$ and that $\cT$ commutes with the multiplication by a $P-$periodic function, i.e. for all $u,v \in H^s(\R)$
$$
\begin{aligned}
&\cT(uv)(x,k)=(\util *_{\B_P} \vtil)(x,k):=\int_{\B_P}\util(x,k-l)\vtil(x,l)dl,\\
&\cT(Vu)(x,k)=V(x)\util(x,k) \quad \text{for all } V\in C(\R) \text{ such that } V(x+P)=V(x) \ \text{for all }x\in \R,
\end{aligned}
$$
where in the convolution the $k-$periodicity in \eqref{E:util_per} needs to be used when $k-l\notin \B_P$.

Our analysis does not make use of the isomorphism as we work in $L^1(\B_P,H^s(0,P))$ rather than $L^2(\B_P,H^s(0,P))$. The reason is that in $L^2$ too many powers of $\eps$ are lost such that the resulting asymptotic error is not $o(1)$ on the desired time interval $[0,O(\eps^{-1})]$, see also \cite{SU01} for the same issue. The norm in $L^1(\B_P,H^s(0,P))$ is
$$\|\util\|_{L^1(\B_P,H^s(0,P))} = \int_{\B_P} \|\util(\cdot,k)\|_{H^s(0,P)}dk.$$
Unfortunately, the isomorphism property of $\cT$ is lost when $L^2$ is replaced by $L^1$. On the other hand, for $s>1/2$ the supremum norm of $u:=\cT^{-1}\util$ can be controlled by $\|\util\|_{L^1(\B_P,H^s(0,P))}$. This means that the supremum norm of the error will be controlled if we estimate the Bloch transform of the error in $L^1(\B_P,H^s(0,P))$. Moreover, for $\util \in L^1(\B_P,H^s(0,P))$ with $s>1/2$ the function $u(x)$ decays as $|x|\to \infty$ as the next lemma shows.
\blem\label{L:sup_control}
Let $s>1/2$.  There is $c>0$ such that for all $\util \in L^1(\B_P,H^s(0,P))$ which satisfy \eqref{E:util_per}, we have for the function $u:=\cT^{-1} \util$ 
$$|u(x)|\leq c\|\util\|_{L^1(\B_P,H^s(0,P))} \quad \text{and} \ u(x)\to 0 \text{ as } |x| \to \infty.$$
\elem
\bpf
$$
\begin{aligned}
\|u\|_{C^0_b(\R)}\leq \int_{\B_P}\|\util(\cdot,k)\|_{C^0_b(0,P)}dk \leq c \int_{\B_P}\|\util(\cdot,k)\|_{H^s(0,P)}dk=c\|\util\|_{L^1(\B_P,H^s(0,P))}
\end{aligned}
$$
due to Sobolev's embedding. The proof of the decay follows the same lines as the proof of Riemann-Lebesgue's lemma. One approximates $\util$ by $v\in C^\infty(\B_P,H^s(0,P))$ with $v(x,k+2\pi/P)=e^{-\ri x 2\pi/P }v(x,k)$ and $v(x+P,k)=v(x,k)$ for all $x$ and $k$, and uses integration by parts. \epf

It is also easy to establish the following algebra property of $L^1(\B_P,H^s(0,P))$, which is needed for the treatment of the nonlinearity.
\blem\label{L:algeb_L1Hs}
Let $\util,\vtil\in L^1(\B_P,H^s(0,P))$ with $s>1/2$. Then 
$$\|\util \ast_{\B_P}\vtil\|_{L^1(\B_P,H^s(0,P))} \leq c \|\util\|_{L^1(\B_P,H^s(0,P))}\|\vtil\|_{L^1(\B_P,H^s(0,P))}.$$
\elem
\bpf
$$
\begin{aligned}
\|\util \ast_{\B_P}\vtil\|_{L^1(\B_P,H^s(0,P))} & \leq c\int_{\B_P}\int_{\B_P}\|\util(\cdot,k-l)\|_{H^s(0,P)}\|\vtil(\cdot,l)\|_{H^s(0,P)} dldk \\
&\leq c \|\util\|_{L^1(\B_P,H^s(0,P))}\|\vtil\|_{L^1(\B_P,H^s(0,P))},
\end{aligned}
$$
where the first inequality follows by the algebra property of $H^s$ in one dimension, i.e. $\|fg\|_{H^s(\Omega)}\leq c\|f\|_{H^s(\Omega)}\|g\|_{H^s(\Omega)}$ for $\Omega\subset \R$ and $s>1/2$. The second step follows by Young's inequality for convolutions.
\epf

For each $k\in \B_P$ the eigenfunctions $(p_n(\cdot,k))_{n\in \N}$ of $L(k)$ in  \eqref{E:Bloch_ev_prob} are complete in $L^2(0,P)$. Hence, for each $k\in \B_P$ fixed and $s\geq 0$ a function $\util(\cdot,k)\in H^s(0,P)$ can be expanded in $(p_n(\cdot,k))_{n\in \N}$. We denote the expansion operator by $D(k)$, i.e.
$$D(k):\util(\cdot,k)\mapsto \Uvec(k):=(\langle\util(\cdot,k),p_n(\cdot,k)\rangle_P)_{n\in \N}.$$
As shown in Lemma 3.3 of \cite{BSTU06}, $D(k)$ is an isomorphism between $H^s(0,P)$ and 
$$l_s^2=\{\vec{v}\in l^2(\C):\|\vec{v}\|^2_{l_s^2}=\sum_{n\in \N}n^{2s}|v_n|^2 <\infty\}$$ 
for all $s\geq 0$ with 
$$\|D(k)\|,\|D^{-1}(k)\|\leq c<\infty \quad \text{for all } k\in \B_P.$$ 

For the expansion coefficients $\Uvec$ we define the space
$$\cX(s):=L^1(\B_P,l^2_s) \quad \text{with the norm } \|\Uvec\|_{\cX(s)}=\int_{\B_P} \|\Uvec(k)\|_{l^2_s} dk.$$
Our working space for  $\Uvec$ will be $\cX(s)$ with $s>1/2$.

As a simple consequence of the properties of $D(k)$ we have
\blem\label{L:D_isom}
The operator $\cD:L^1(\B_P,H^s(0,P))\to \cX(s), \util \mapsto \Uvec$ is an isomorphism for any $s\geq 0$.
\elem
\bpf
There are $c_1,c_2>0$ such that for any $\Uvec(k)=D(k) \util(\cdot,k)$ we have $\|\util(\cdot,k)\|_{H^s(0,P)}\leq c_1 \|\Uvec(k)\|_{l_s^2}$ and $\|\Uvec(k)\|_{l_s^2} \leq c_2\|\util(\cdot,k)\|_{H^s(0,P)}$ for all $k\in \B_P$. Hence
$$\|\util\|_{L^1(\B_P,H^s(0,P))} \leq c_1\int_{\B_P}\|\Uvec(k)\|_{l_s^2}dk = c_1\|\Uvec\|_{\cX(s)}$$
and
$$\|\Uvec\|_{\cX(s)}\leq c_2\int_{\B_P}\|\util(\cdot,k)\|_{H^s(0,P)}dk = c_2\|\util\|_{L^1(\B_P,H^s(0,P))}.$$
\epf

Note that the $k-$periodicity in \eqref{E:util_per} implies that $\Uvec=\cD \util$ satisfies
\beq\label{E:U_per}
\Uvec(k+\tfrac{2\pi}{P})=\Uvec(k) \qquad \text{for all } k\in \R.
\eeq

The above function spaces are used in our analysis in the following way. We define an extended (compared to $\uapp$) ansatz for the approximate solution in the $\Uvec$-variables and show that it lies in $\cX(s)$ for $s<3/2$. The residual of this ansatz in equation \eqref{E:PNLS} is then estimated in $\|\cdot\|_{\cX(s)}$, where some terms are transformed by $\cD^{-1}$ to $L^1(\B_P,H^s(0,P))$ and estimated in $\|\cdot\|_{L^1(\B_P,H^s(0,P))}$. The approximation error is then estimated in the $\Uvec$-variables (i.e. in $\|\cdot\|_{\cX(s)}$), in which the equation becomes an infinite ODE system, by Gronwall's inequality. The supremum of the error in the physical $u$ variables then satisfies the same estimate due to Lemma \ref{L:sup_control} because $s>1/2$.


\subsection{Proof of Theorem \ref{T:main} for Rational $k_\pm$}\label{S:rational_pf}

Let us recall that in the rational case we work with the period $P=Q=2N\pi$, see Sec. \ref{S:rational}. Also $k_\pm \in \tfrac{1}{N}\Z$, such that in $\B_Q$ the points $k_\pm$ are identified with $0$. 

We will prove Theorem \ref{T:main} for $u_\text{app}$ in \eqref{E:uapp_rat} and $(A_+,A_-)$ a solution of \eqref{E:CME_rat}. This is, however, equivalent to the theorem since $(A_+,A_-)_{\eqref{E:CME}} = \tfrac{1}{\sqrt{N}}(A_+,A_-)_{\eqref{E:CME_rat}}$ and $p_\pm(x)e^{\ri k_\pm x}=\sqrt{N}q_\pm(x)$ such that $u_\text{app}$ in \eqref{E:uapp} and in \eqref{E:uapp_rat} are identical.

The rigorous justification of the asymptotic model \eqref{E:CME_rat} will be carried out in the Bloch variables $\Uvec=\Uvec(k,t)$. Instead of the approximate ansatz $\Uvec_{\text{app}}:=\cD\cT(\uapp)$ we use its modification $\Uvecext$ which is supported in $k$ only near the wavenumbers of the corresponding carrier waves, i.e. near $k=0$, and includes correction terms supported near the new $k-$points generated by $W^{(R)}_Q\uapp$. Note that the  nonlinearity does not generate new $k-$neighborhoods when $k_\pm =0 \mod \tfrac{1}{N}$. This is a standard approach, where the concentration points of the residual for $u_\text{app}$ are identified and the modified ansatz is chosen to be supported only in neighborhoods of these points. The problem can then be easily decomposed according to the disjoint intervals in the support on $\B_Q$. To motivate the choice of $\Uvecext$, we study the $\cT$-transform of the residual for $\uapp$ with a compact support of $\hat{A}_\pm(\cdot,T)$. Assuming that\footnote{The support $[-\eps^{-1/2},\eps^{-1/2}]$ has been chosen a-posteriori. If one assumes $\text{supp}(\hat{A}_\pm(\cdot,T))\subset [-\eps^{r-1},\eps^{r-1}]$ for some $r\in (0,1)$, it turns out in the proof of Lemma \ref{L:res_est} that $r\leq 1/2$ is optimal.} $\text{supp}(\hat{A}_\pm(\cdot,T))\subset [-\eps^{-1/2},\eps^{-1/2}]$, we get
$$
\begin{aligned}
\cT(\uapp)(x,k,t)=&\eps^{-1/2}\sum_{\pm}q_\pm(x)\sum_{\eta\in \tfrac{1}{N}\Z}\hat{A}_\pm\left(\tfrac{k+\eta}{\eps},\eps t\right)e^{\ri(\eta x- \omega_0 t)}\\
=&\eps^{-1/2}\sum_{\pm}q_\pm(x)\hat{A}_\pm\left(\tfrac{k}{\eps},\eps t\right)e^{-\ri \omega_0 t}, \quad k \in \B_Q
\end{aligned}
$$
because $\eps^{-1}(k+\eta) \in \text{supp}(\hat{A}_\pm(\cdot,T))$ for some $k\in \B_Q$ and $\eta\in \tfrac{1}{N}\Z$ implies $\eta=0$. Hence, if $\text{supp}(\hat{A}_\pm(\cdot,T))\subset [-\eps^{-1/2},\eps^{-1/2}]$, then for $k\in\B_Q$
$$\begin{aligned}
&\cT(\text{PNLS}(\uapp))(x,k,t)= \eps^{1/2}e^{-\ri \omega_0 t}\sum_{\pm}\left[\ri \pa_T\hat{A}_\pm(K,T)q_\pm(x)+2\ri K\hat{A}_\pm(K,T)q_\pm'(x)\phantom{\sum_{m\in \Z_R,\eta \in \cS_m}}\right.\\
&\left.-W_Q^{(1)}(x)q_\pm(x)\hat{A}_\pm(K,T)-\sum_{m\in \Z_R,\eta \in \cS_m} a_m\hat{A}_\pm\left(\frac{k-mk_W+\eta}{\eps},T\right)e^{\ri \eta x}q_\pm(x)\right]\\
& -\eps^{1/2}  e^{-\ri \omega_0 t}\sigma(x) \sum_{s_1,s_2,s_3\in \{+,-\}}q_{s_1}(x)\overline{q_{s_2}}(x)q_{s_3}(x)\left(\hat{A}_{s_1}\ast \hat{\overline{A}}_{s_2}\ast \hat{A}_{s_3}\right)(K,T)\\
&  - \eps^{3/2}e^{-\ri \omega_0 t}\sum_{\pm}p_\pm(x)K^2\hat{A}_\pm(K,T),
\end{aligned}
$$
where 
$$K:=\eps^{-1}k, \qquad T = \eps t, \qquad \cS_m:=\{\eta\in \tfrac{1}{N} \Z: mk_W-\eta\in \overline{\B_Q}=[-\tfrac{1}{2N},\tfrac{1}{2N}]\},$$ 
and the convolution $\ast$ is in the $K$-variable, $(\hat{f}\ast\hat{g})(K)=\int_\R\hat{f}(K-s)\hat{g}(s) ds$. Note that the set $\cS_m$ has at most two elements, namely $\cS_m\subset \{0,\tfrac{1}{N}\}$ or $\cS_m\subset \{0,-\tfrac{1}{N}\}$. 

Because of the support of $\hat{A}_\pm$ the above residual is supported in $k$ in  intervals of radius at most $3\eps^{1/2}$ centered at $k=0$, $k\in k_W \Z_R$ and at all $\tfrac{1}{N}\Z$-shifts of these points, cf. the periodicity \eqref{E:util_per} of the Bloch transform in $k$. Due to assumption (H2) in the case of irrational $k_W$ the number of distinct support centers from $k_W \Z_R + \tfrac{1}{N}\Z$ which lie in $\overline{\B_Q}$ is finite. For $k_W\in \Q$ is $\{k_W \Z_R + \tfrac{1}{N}\Z\}\cap \B_Q$ finite even if infinitely many coefficients $a_n$ are nonzero. Note that we need $\overline{\B_Q}$ rather than only $\B_Q$ because a support interval centered at $k=-\tfrac{1}{2N}$ intersects $\B_Q$ for any $\eps>0$.

We denote this \textit{finite set} of support centers by
\beq\label{E:KR}
\cK_R:=(k_W \Z_R + \tfrac{1}{N}\Z)\cap \overline{\B_Q}=\{k\in \overline{\B_Q}: k=mk_W+\eta \text{ for some } m\in \Z_R, \eta\in \cS_m\}
\eeq
and label its elements by $\kappa_j$:
$$\cK_R=\{\kappa_1,\dots, \kappa_J\} \text{ for some } J\in \N.$$

Our modified (extended) ansatz in the $\Uvec$-variables is thus for $k\in \B_Q$
\beq\label{E:Uext}
\begin{aligned}
U_{n_*}^\text{ext}(k,t) &:= \eps^{-1/2}\Atil_{n_*}(K,T)e^{-\ri \omega_0 t}+\eps^{1/2}\hspace{-0.1cm}\sum_{j=1}^J \Atil_{n_*,j}\left(\frac{k-\kappa_j}{\eps},T\right)e^{-\ri \omega_0 t}\\
U_{n_*+1}^\text{ext}(k,t) &:= \eps^{-1/2}\Atil_{n_*+1}(K,T)e^{-\ri \omega_0 t}+\eps^{1/2}\hspace{-0.1cm}\sum_{j=1}^J \Atil_{n_*+1,j}\left(\frac{k-\kappa_j}{\eps},T\right)e^{-\ri \omega_0 t}\\
U_{n}^\text{ext}(k,t) &:= \eps^{1/2}\Atil_{n}(K,T)e^{-\ri \omega_0 t}+\eps^{1/2}\hspace{-0.1cm}\sum_{j=1}^J \Atil_{n,j}\left(\frac{k-\kappa_j}{\eps},T\right)e^{-\ri \omega_0 t}, \quad n \notin \{n_*,n_*+1\}\\
\end{aligned}
\eeq
with 
$$
\begin{aligned}
&\supp(\Atil_q(\cdot,T)) \cap \eps^{-1}\B_Q \subset [-\eps^{-1/2},\eps^{-1/2}], \quad q\in \{n_*,n_*+1\},\\
&\supp(\Atil_n(\cdot,T)) \cap \eps^{-1}\B_Q \subset [-3\eps^{-1/2},3\eps^{-1/2}], \quad \supp(\Atil_{n,m}(\cdot,T))\cap \eps^{-1}\B_Q \subset [-\eps^{-1/2},\eps^{-1/2}]
\end{aligned}
$$
for all $n\in \N,m\in \Z$. As expected from the formal asymptotics and as shown in detail below, the residual for $\Uvecext$ is small if $\Atil_{n_*}$ and $\Atil_{n_*+1}$ are selected as cut-offs of $\widehat{A}_+$ and $\widehat{A}_-$ respectively. 
To comply with \eqref{E:U_per}, we define $\Uvecext$ with $\tfrac{1}{N}$-periodicity in $k$, i.e. 
$$\Uvecext\left(k+\tfrac{1}{N},t\right)=\Uvecext(k,t) \ \text{for all} \ k\in \R, t\in \R.$$
Below it will be useful to write
\beq\label{E:Uext_decomp}
\Uvecext = \vec{U}^{0,\text{ext}}+\vec{U}^{1,\text{ext}},
\eeq
where 
\beq\label{E:Uext_decomp2}
\vec{U}^{0,\text{ext}}:= \eps^{-1/2}(\Atil_{n_*}(K,T)e_{n_*}+\Atil_{n_*+1}(K,T)e_{n_*+1})e^{-\ri \omega_0 t}\quad \text{and} \ \vec{U}^{1,\text{ext}}:=\Uvecext-\vec{U}^{0,\text{ext}}.
\eeq
Here $e_n$ is the standard $n$-th Euclidean unit vector in $\R^\N$.

In order to determine the residual for $\Uvecext$, let us first formulate equation \eqref{E:PNLS} in the $\Uvec$-variables. Applying first $\cT$ to \eqref{E:PNLS}, we get
$$\left(\ri \pa_t-\cL(x,k)-\eps W^{(1)}_Q(x)\right)\util(x,k,t)-\eps \sum_{m\in \Z_R}a_m \util(x,k-mk_W,t)-\sigma(x) (\util\ast_{\B_Q}\tilde{\overline{u}}\ast_{\B_Q}\util)(x,k,t)=0.$$
Expanding $\util(x,k,t) = \sum_{n\in \N}U_n(k,t)q_n(x,k)$ leads to the (infinitely dimensional) ODE-system
\beq\label{E:PNLS_U}
\left(\ri \pa_t -\Theta(k)-\eps M^{(1)}(k)\right)\Uvec(k,t) -\eps \sum_{m\in \Z_R}M^{(R,m)}(k) \Uvec(k-mk_W,t)+\vec{F}(\Uvec,\Uvec,\Uvec)(k,t)=0
\eeq
parametrized by $k$, where
\begin{align}
&\Theta_{j,j}(k):=\vartheta_j(k) \ \forall j\in \N, \ \Theta_{i,j}:=0 \ \forall i,j\in \N, i\neq j,\notag \\
&M^{(1)}_{i,j}(k):=\langle W^{(1)}_Qq_j(\cdot,k), q_i(\cdot,k)\rangle_{Q}, \quad M^{(R,m)}_{i,j}(k):=a_m\langle q_j(\cdot,k-mk_W), q_i(\cdot,k)\rangle_{Q},\notag\\
&F_j(\Uvec,\Uvec,\Uvec)(k,t) :=-\langle \sigma(\cdot) (\util\ast_{\B_Q}\tilde{\overline{u}}\ast_{\B_Q}\util)(\cdot,k,t),q_j(\cdot,k)\rangle_Q, \ \util(x,k,t)=\sum_{n\in \N}U_n(k,t)q_n(x,k).\label{E:Fj}
\end{align}
Note that $\tilde{\overline{u}}(x,k,t)=\sum_{n\in \N}\overline{U_n}(-k,t)q_n(x,k)$.
Substituting $\Uvec^{\text{ext}}$ in \eqref{E:PNLS_U}, we get the residual
$$
\begin{aligned}
\vec{\Res}(k,t):=&\left(\ri \pa_t -\Theta(k)-\eps M^{(1)}(k)\right)\Uvecext(k,t) -\eps \sum_{m\in \Z_R}M^{(R,m)}(k) \Uvecext(k-mk_W,t)\\
&+\vec{F}(\Uvecext,\Uvecext,\Uvecext)(k,t)
\end{aligned}
$$ 
where in \eqref{E:Fj} $\util^\text{ext}(x,k,t) := \sum_{n\in \N}U_n^\text{ext}(k,t)q_n(x,k)$ in $\vec{F}$. For $k\in \B_Q$ we have
$$
\begin{aligned}
&\Res_{n_*}(k,t)= \eps^{1/2}\left[\left(\ri \pa_T -\eps^{-1}(\vartheta_{n_*}(k)-\omega_0)-M^{(1)}_{n_*,n_*}(k)\right)\Atil_{n_*}(K,T) \phantom{\sum_{m\in \Z_R,\eta\in \cS_m}}\right.\\
& + \eps^{-1/2}F_{n_*}(\UvecOext,\UvecOext,\UvecOext)(k,t)e^{\ri \omega_0 t} - M^{(1)}_{n_*,n_*+1}(k)\Atil_{n_*+1}(K,T)  \\
&- \left(\vartheta_{n_*}(k)-\omega_0\right)\sum_{j=1}^J \Atil_{n_*,j}\left(\frac{k-\kappa_j}{\eps},T\right)\\
&-\sum_{m\in \Z_R,\eta\in \cS_m}M^{(R,m)}_{n_*,n_*}(k)\Atil_{n_*}\left(\frac{k-mk_W+\eta}{\eps},T\right)\\
&\left.-\sum_{m\in \Z_R,\eta\in \cS_m}M^{(R,m)}_{n_*,n_*+1}(k)\Atil_{n_*+1}\left(\frac{k-mk_W+\eta}{\eps},T\right)\right]e^{-\ri \omega_0 t} + \text{h.o.t.},
\end{aligned}
$$
$$
\begin{aligned}
&\Res_{n_*+1}(k,t)= \eps^{1/2}\left[\left(\ri \pa_T -\eps^{-1}(\vartheta_{n_*+1}(k)-\omega_0)-M^{(1)}_{n_*+1,n_*+1}(k)\right)\Atil_{n_*+1}(K,T) \phantom{\sum_{m\in \Z_R,\eta\in \cS_m}}\right.\\
& + \eps^{-1/2}F_{n_*+1}(\UvecOext,\UvecOext,\UvecOext)(k,t)e^{\ri \omega_0 t} - M^{(1)}_{n_*+1,n_*}(k)\Atil_{n_*}(K,T)  \\
&-  \left(\vartheta_{n_*+1}(k)-\omega_0\right)\sum_{j=1}^J\Atil_{n_*+1,j}\left(\frac{k-\kappa_j}{\eps},T\right)\\
&-\sum_{m\in \Z_R,\eta\in \cS_m}M^{(R,m)}_{n_*+1,n_*}(k)\Atil_{n_*}\left(\frac{k-mk_W+\eta}{\eps},T\right)\\
&\left.-\sum_{m\in \Z_R,\eta\in \cS_m}M^{(R,m)}_{n_*+1,n_*+1}(k)\Atil_{n_*+1}\left(\frac{k-mk_W+\eta}{\eps},T\right)\right]e^{-\ri \omega_0 t} + \text{h.o.t.},
\end{aligned}
$$
and for $n\in \N\setminus\{n_*,n_*+1\}$
$$
\begin{aligned}
&\Res_{n}(k,t)= \eps^{1/2}\left[(\omega_0-\vartheta_{n}(k))\Atil_{n}(K,T) -(M_{n,n_*}^{(1)}(k)\Atil_{n_*}(K,T) +M_{n,n_*+1}^{(1)}(k)\Atil_{n_*+1}(K,T))\phantom{\sum_{m\in \Z_R,\eta\in \cS_m}}\right.\\
&\left.+\eps^{-1/2}F_{n}(\UvecOext,\UvecOext,\UvecOext)(k,t)e^{\ri \omega_0 t} - \left(\vartheta_{n}(k)-\omega_0\right)\sum_{j=1}^J\Atil_{n,j}\left(\frac{k-\kappa_j}{\eps},T\right)\right.\\
&-\sum_{m\in \Z_R,\eta\in \cS_m}M^{(R,m)}_{n,n_*}(k)\Atil_{n_*}\left(\frac{k-mk_W+\eta}{\eps},T\right)\\
&\left.-\sum_{m\in \Z_R,\eta\in \cS_m}M^{(R,m)}_{n,n_*+1}(k)\Atil_{n_*+1}\left(\frac{k-mk_W+\eta}{\eps},T\right)\right]e^{-\ri \omega_0 t} + \text{h.o.t.}.
\end{aligned}
$$
Note that
$$F_{n}(\UvecOext,\UvecOext,\UvecOext)=\eps^{-3/2}\sum_{\alpha,\beta,\gamma\in\{n_*,n_*+1\}}F_{n}(\Atil_{\alpha}e_\alpha,\Atil_{\beta}e_\beta,\Atil_{\gamma}e_\gamma)e^{-\ri \omega_0 t}$$
and
$$
\begin{aligned}
F_{n}&(\Atil_{\alpha}e_\alpha,\Atil_{\beta}e_\beta,\Atil_{\gamma}e_\gamma)(k,t) \\
& = \int_{-2\eps^{1/2}}^{2\eps^{1/2}}\int_{-\eps^{1/2}}^{\eps^{1/2}} b_{\alpha\beta\gamma}^{(n)}(k,k-s,s-l,l)\Atil_\alpha\left(\tfrac{k-s}{\eps},T\right)\tilde{\overline{A}}_\beta\left(\tfrac{s-l}{\eps},T\right)\Atil_\gamma\left(\tfrac{l}{\eps},T\right) dl ds\\
&=\eps^2\int_{-2\eps^{-1/2}}^{2\eps^{-1/2}}\int_{-\eps^{-1/2}}^{\eps^{-1/2}} b_{\alpha\beta\gamma}^{(n)}(\eps K,\eps(K-s),\eps(s-l),\eps l)\Atil_\alpha(K-s,T)\tilde{\overline{A}}_\beta(s-l,T)\Atil_\gamma(l,T) dl ds
\end{aligned}
$$
with 
$$
b_{\alpha\beta\gamma}^{(n)}(k,r,s,l):=\langle \sigma(\cdot) q_\alpha(\cdot,r)\overline{q_\beta}(\cdot,-s)q_\gamma(\cdot,l),q_n(\cdot,k)\rangle_Q.
$$
Hence $\eps^{-1/2}F_{n}(\UvecOext,\UvecOext,\UvecOext)(k,t)$ is $O(1)$. Also the terms $\eps^{-1}(\vartheta_{j}(k)-\omega_0)\Atil_{j}$ for $j=n_*,n_*+1$ are $O(1)$ as can be seen by the Taylor expansion $\vartheta_{j}(k)\sim\omega_0 \pm kc_g=\omega_0\pm\eps K c_g  \ (\eps \to 0, k\in \supp \Atil_j(\cdot,T)\cap \eps^{-1}\B_Q)$ for $j=n_*,n_*+1$ respectively.

The ``h.o.t.'' in the residual stands for terms of higher order in $\eps$ and consists of the following terms
\beq\label{E:hot}
\begin{aligned}
&\eps \pa_T  \vec{U}^{1,\text{ext}}, \quad \eps M^{(1)}\vec{U}^{1,\text{ext}}, \quad \eps \sum_{m\in \Z_R}M^{(R,m)}\vec{U}^{1,\text{ext}}(\cdot-mk_W,t),\\
& \vec{F}(\UvecOext,\vec{U}^{1,\text{ext}},\UvecOext)+2\vec{F}(\UvecOext,\UvecOext,\vec{U}^{1,\text{ext}}),
\end{aligned}
\eeq
and nonlinear terms quadratic or cubic in $\vec{U}^{1,\text{ext}}$.

Note that if we approximate in the first two lines of $\Res_{n_*}$ the function $\vartheta_{n_*}(k)$ by $\vartheta_{n_*}(0)+k c_g$, the function $M^{(1)}_{n_*,j}(k)$ by $M^{(1)}_{n_*,j}(0)$ for $j=n_*,n_*+1$ and in $F_{n_*}$ approximate $q_j(\cdot,k)$ by $q_j(\cdot,0)$, such that $b_{\alpha\beta\gamma}^{(n_*)}(\eps K,\eps(K-s),\eps(s-l),\eps l)$ is replaced by $b_{\alpha\beta\gamma}^{(n_*)}(0,0,0,0)$, we recover the left hand side of the first equation in the CMEs \eqref{E:CME_rat}. Similarly for the first two lines in $\Res_{n_*+1}$ and the second equation in \eqref{E:CME_rat}. This folows from the identities
\beq\label{E:coef_ident}
\begin{aligned}
&M^{(1)}_{n_*,n_*+1}(0)=M^{(1)}_{n_*+1,n_*}(0)=-\kappa, \quad  M^{(1)}_{n_*,n_*}(0)=M^{(1)}_{n_*+1,n_*+1}(0)=-\kappa_s,\\
&b^{(n_*)}_{n_*,n_*,n_*}(0,0,0)=b^{(n_*+1)}_{n_*+1,n_*+1,n_*+1}(0,0,0)=b^{(n_*)}_{n_*,n_*+1,n_*+1}(0,0,0)=b^{(n_*+1)}_{n_*,n_*,n_*+1}(0,0,0)=-\alpha_Q\\
&b^{(n_*)}_{n_*,n_*,n_*+1}(0,0,0)=b^{(n_*)}_{n_*+1,n_*+1,n_*+1}(0,0,0)=b^{(n_*+1)}_{n_*,n_*,n_*}(0,0,0)=b^{(n_*+1)}_{n_*+1,n_*+1,n_*}(0,0,0)=-\beta_Q\\
&b^{(n_*)}_{n_*+1,n_*,n_*+1}(0,0,0)=b^{(n_*+1)}_{n_*,n_*+1,n_*}(0,0,0)=-\gamma_Q.
\end{aligned}
\eeq
Therefore, the first two lines in $\Res_{n_*}$, $\Res_{n_*+1}$ will be close to zero (precisely $O(\eps^{5/2})$ in $\|\cdot\|_{L^1(\B_Q)}$ as shown below) after setting
\beq\label{E:Atil_def}
\Atil_{n_*}(K,T):= \chi_{[-\eps^{-1/2},\eps^{-1/2}]}(K)\hat{A}_+(K,T), \quad \Atil_{n_*+1}(K,T):= \chi_{[-\eps^{-1/2},\eps^{-1/2}]}(K)\hat{A}_-(K,T).
\eeq
The remaining three lines in the formal $O(\eps^{1/2})$ part of $\Res_{n_*}$ and $\Res_{n_*+1}$ can then be made exactly zero by setting
\beq\label{E:Atilcor_def}
\begin{aligned}
\Atil_{q,j}(K,T):=(&\omega_0-\vartheta_q(\kappa_j+\eps K))^{-1}\sum_{\stackrel{m\in \Z_R}{mk_W\in \kappa_j+\tfrac{1}{N}\Z}}\left(M^{(R,m)}_{q,n_*}(\kappa_j+\eps K)\Atil_{n_*}\left(K,T\right)\right.\\
&\left.+M^{(R,m)}_{q,n_*+1}(\kappa_j+\eps K)\Atil_{n_*+1}\left(K,T\right)\right), \quad q\in\{n_*,n_*+1\}, j\in\{1,\dots,J\},
\end{aligned}
\eeq  
where the $1/N$-periodicity of $\vartheta_q$ and $M_{q,n}^{(R,m)}$ was used. Note that for $\eps>0$ small enough we have
$$\min_{q\in\{n_*,n_*+1\}, j\in \{1,\dots,J\}, |K|\leq \eps^{-1/2}}\left|\vartheta_q(\kappa_j+\eps K)-\omega_0\right|>c>0$$
because $\omega_0=\vartheta_{n_*}(0)=\vartheta_{n_*+1}(0)$, the functions $\vartheta_{q}(k)$ are monotonous on $(-\tfrac{1}{2N},0)$ and $(0,\tfrac{1}{2N}]$ and $\tfrac{1}{N}$-periodic and because ${\rm dist}(\cK_R,\tfrac{1}{N}\Z)> \delta >0$. The last property follows from the finiteness of $\cK_R$.

Finally, the formal $O(\eps^{1/2})$ part of $\Res_{n}, n\notin\{n_*,n_*+1\}$ vanishes if we set for $n\in \N\setminus\{n_*,n_*+1\}$
\beq\label{E:Atiln_def}
\begin{aligned}
\Atil_{n}(K,T):=&(\omega_0-\vartheta_n(\eps K))^{-1}\left(M^{(1)}_{n,n_*}(\eps K)\Atil_{n_*}(K,T)+M^{(1)}_{n,n_*+1}(\eps K)\Atil_{n_*+1}(K,T) \right.\\
&-\int_{-2\eps^{-1/2}}^{2\eps^{-1/2}}\int_{-\eps^{-1/2}}^{\eps^{-1/2}} \sum_{\alpha,\beta,\gamma\in \{n_*,n_*+1\}}b_{\alpha\beta\gamma}^{(n)}(\eps K,\eps(K-s),\eps(s-l),\eps l)\times\\
 &\left. \qquad \Atil_\alpha(K-s,T)\tilde{\overline{A}}_\beta(s-l,T)\Atil_\gamma(l,T) dl ds \right)
\end{aligned}
\eeq  
and
\beq\label{E:Atilncor_def}
\begin{aligned}
\Atil_{n,j}(K,T):=(&\omega_0-\vartheta_n(\kappa_j+\eps K))^{-1}\sum_{\stackrel{m\in \Z_R}{nk_W\in \kappa_j+\tfrac{1}{N}\Z}}\left(M^{(R,m)}_{n,n_*}(\kappa_j+\eps K)\Atil_{n_*}(K,T)\right.\\
&\left.+M^{(R,m)}_{n,n_*+1}(\kappa_j+\eps K)\Atil_{n_*+1}(K,T) \right), \quad j\in \{1,\dots,J\}.
\end{aligned}
\eeq  
The factors $(\omega_0-\vartheta_n(\eps K))^{-1}$ and $(\omega_0-\vartheta_n(\kappa_j+\eps K))^{-1}$ are again bounded because
$$\min_{n\in\N\setminus\{n_*,n_*+1\}, k\in \R}\left|\omega_0-\vartheta_n(k)\right|>c>0.$$
This follows from the transversal crossing of $\vartheta_{n_*}$ and $\vartheta_{n_*+1}$ at $(0,\omega_0)$ and because eigenvalue functions do not overlap in one dimension.

In the remainder of the argument, first in Lemma \ref{L:res_est} we estimate $\|\vec{\Res}(\cdot,t)\|_{\cX(s)}$ under the conditions \eqref{E:Uext}, \eqref{E:Atil_def}, \eqref{E:Atilcor_def}, \eqref{E:Atiln_def}, and \eqref{E:Atilncor_def}, and then use a Gronwall argument and an estimate of $\|\Uvecext(\cdot,t)-\Uvec^\text{app}(\cdot,t)\|_{\cX(s)}$ in Lemma \ref{L:diff_app_ext} to conclude the proof of Theorem \ref{T:main}.

Below we will need some asymptotics of $\left|\vartheta_n(k)-\omega_0\right|$ for $n\to \infty$. In particular, see p. 55 in \cite{Hoerm_85}, there are constants $c_1,c_2>0$ such that
\beq\label{E:band_as}
c_1n^2 \leq \vartheta_n(k) \leq c_2n^2 \quad \text{for all} \ k\in \B_Q, n\in \N.
\eeq

The ansatz components $\UvecOext$ and $\vec{U}^{1,\text{ext}}$ can be estimated as follows.
\blem\label{L:ext_est}
Assume (H2)-(H4) and let $(A_+,A_-)(X,T)$ be a solution of \eqref{E:CME_rat} such that $\widehat{A}_\pm\in C([0,T_0],L^1(\R))$ and assume \eqref{E:Atil_def}, \eqref{E:Atilcor_def}, \eqref{E:Atiln_def} and \eqref{E:Atilncor_def}. Then there is $\eps_0>0$ and 
$$c=c\left(\max_{T\in[0,T_0]}\|\hat{A}_+(\cdot,T)\|_{L^1(\R)},\max_{T\in[0,T_0]}\|\hat{A}_-(\cdot,T)\|_{L^1(\R)}\right)>0$$ 
such that for all $s\in (0,3/2),  \eps\in(0,\eps_0)$ and all $t\in [0,\eps^{-1}T_0]$
$$\|\UvecOext\|_{\cX(s)}\leq c\eps^{1/2}, \qquad \|\vec{U}^{1,\text{ext}}\|_{\cX(s)}\leq c\eps^{3/2}.$$
\elem
\bpf
From \eqref{E:Atil_def} we get
$$\|\UvecOext(\cdot,t)\|_{\cX(s)} \leq c \eps^{1/2}\sum_{\pm}\|\hat{A}_\pm(\cdot, \eps t)\|_{L^1(\R)}.$$

Next, note that due to the $C^2$ smoothness of $W$ there is a constant $c>0$ such that
\beq\label{E:MRm_est}
|M^{(R,m)}_{n,j}(k)|\leq \frac{c}{m^2} \ \text{for all }  j\in\{n_*,n_*+1\}, n\in \N \text{ and } k\in \B_Q.
\eeq
Also, one obviously has
\beq\label{E:Mb_est}
|M^{(1)}_{n,j}(k)|\leq c, |b_{\alpha,\beta,\gamma}^{(n)}(k,l,s,r)|\leq c
\eeq
for all $k,l,s,r\in \B_Q$, all $j,\alpha,\beta,\gamma\in \{n_*,n_*+1\}$ and all $n\in \N$

Finally, we use \eqref{E:MRm_est}, \eqref{E:Mb_est}, \eqref{E:Atilcor_def},\eqref{E:Atiln_def},\eqref{E:Atilncor_def}, and \eqref{E:band_as} to estimate 
$$
\begin{aligned}
&\|\vec{U}^{1,\text{ext}}\|_{\cX(s)}= \left\|\left(\sum_{n\in \N}n^{2s}|U_n^{1,\text{ext}}(\cdot,t)|^2\right)^{1/2}\right\|_{L^1(\B_Q)}\\
&\leq c\eps^{1/2} \left(\sum_{n\in \N}n^{2s-4}\right)^{1/2}\left(\left(\left\|\Atil_{n_*}(\eps^{-1}\cdot,T)\right\|_{L^1(\B_Q)}+\left\|\Atil_{n_*+1}(\eps^{-1}\cdot,T)\right\|_{L^1(\B_Q)}\right)\left(1+\sum_{m\in \Z_R}\frac{c}{m^2}\right)\right.\\
&\left.+\sum_{\alpha,\beta,\gamma\in \{n_*,n_*+1\}}\left\|\left(|\Atil_\alpha|*_{\B_Q}|\Atil_\beta|*_{\B_Q}|\Atil_\gamma|\right)(\eps^{-1}\cdot,T)\right\|_{L^1(\B_Q)}\right)\\
&\leq c\eps^{3/2} \left(\sum_{\pm}\|\Ahat_\pm(\cdot,T)\|_{L^1(\R)}+\sum_{\xi,\zeta,\theta\in \{+,-\}}\|\Ahat_\xi(\cdot,T)\|_{L^1(\R)}\|\Ahat_\zeta(\cdot,T)\|_{L^1(\R)}\|\Ahat_\theta(\cdot,T)\|_{L^1(\R)}\right),
\end{aligned}
$$
where the factor $cn^{-4}$ in the sum over $n\in \N$ comes from $(\omega_0 -\vartheta_n(\eps K))^{-2}$ using \eqref{E:band_as}. The assumption $s<3/2$ implies the summability of the series in $n$. The bound on the convolution follows from Young's inequality for convolutions. 

The estimates $\|\UvecOext(\cdot,t)\|_{\cX(s)}\leq c\eps^{1/2}, \|\vec{U}^{1,\text{ext}}(\cdot,t)\|_{\cX(s)}\leq c\eps^{3/2}$ for some $c>0$ and all $t\in [0,\eps^{-1}T_0]$ and $s\in (0,3/2)$ now follow from $\widehat{A}_\pm\in C([0,T_0],L^1(\R))$.
\epf

\blem\label{L:res_est}
Assume (H1)-(H4). Let $s\in(1/2,3/2)$ and let $(A_+,A_-)(X,T)$ be a solution of \eqref{E:CME_rat} such that $\hat{A}_\pm\in C^1([0,T_0],L^1_{s_A}(\R))$ for some $T_0>0$ and some $s_A\geq 2$. Choosing $\Uvecext$ according to \eqref{E:Uext}, \eqref{E:Atil_def}, \eqref{E:Atilcor_def}, \eqref{E:Atiln_def}, and \eqref{E:Atilncor_def}, there exists $\eps_0>0$ and 
$$
\begin{aligned}
C_\text{Res}=&C_\text{Res}\left(\max_{\pm,T\in[0,T_0]}\|\hat{A}_\pm(\cdot,T)\|_{L^1_{s_A}(\R)},\max_{\pm,T\in[0,T_0]}\|\pa_T\hat{A}_\pm(\cdot,T)\|_{L^1_{s_A}(\R)}\right)>0
\end{aligned}
$$ 
such that for all $\eps \in (0,\eps_0)$ and all $t\in [0,\eps^{-1}T_0]$ the residual $\vec{\Res}$ of $\Uvecext$ in \eqref{E:PNLS_U} satisfies
$$
\|\vec{\Res}(\cdot,t)\|_{\cX(s)} \leq C_\text{Res} \eps^{5/2}.
$$
\elem
\bpf
First we show that for $n=n_*,n_*+1$ the formal $O(\eps^{1/2})$ part in $\Res_{n}$ satisfies $\|\Res_{n}(\cdot,t)\|_{L^1(\B_Q)}\leq c \eps^{5/2}$ for all $t\in [0,\eps^{-1}T_0]$. Choosing $\Atil_{n}$ and $\Atil_{n,m}$ as in \eqref{E:Atil_def}, \eqref{E:Atilcor_def}, \eqref{E:Atiln_def}, and \eqref{E:Atilncor_def} with $(A_+,A_-)$ being a solution of \eqref{E:CME_rat}, we get
$$\|\Res_{n_*}(\cdot,t)\|_{L^1(\B_Q)}\leq \eps^{-1/2}I_1(T)+\eps^{1/2}\sum_{j=2}^4I_j(T) + \text{h.o.t.},$$
where 
$$
\begin{aligned}
I_1(T)&:=\int_{-\eps^{1/2}}^{\eps^{1/2}}|(\vartheta_{n_*}(k)-\omega_0-c_g\tfrac{k}{\eps})\Atil_{n_*}\left(\tfrac{k}{\eps},T\right)|dk, \\
I_2(T)&:=\int_{-\eps^{1/2}}^{\eps^{1/2}}|(M^{(1)}_{n_*,n_*}(k)-M^{(1)}_{n_*,n_*}(0))\Atil_{n_*}(\tfrac{k}{\eps},T)|dk, \\
I_3(T)&:=\int_{-\eps^{1/2}}^{\eps^{1/2}}|(M^{(1)}_{n_*,n_*+1}(k)-M^{(1)}_{n_*,n_*+1}(0))\Atil_{n_*+1}\left(\tfrac{k}{\eps},T\right)|dk, \\
I_4(T)&:=\hspace{-0.4cm}\sum_{\alpha,\beta,\gamma\in\{n_*,n_*+1\}}\hspace{-0.1cm}J_{\alpha\beta\gamma}(T), \quad J_{\alpha\beta\gamma}(T):=\int_{-3\eps^{1/2}}^{3\eps^{1/2}}\left| \int_{-2\eps^{-1/2}}^{2\eps^{-1/2}}\int_{-\eps^{-1/2}}^{\eps^{-1/2}} b_{\alpha\beta\gamma}^{(n_*)}(k,k-\eps s,\eps(s-l),\eps l) \times \right.\\
&\left. \Atil_\alpha\left(\tfrac{k}{\eps}-s,T\right)\tilde{\overline{A}}_\beta(s-l,T)\Atil_\gamma(l,T) dl ds  - b_{\alpha\beta\gamma}^{(n_*)}(0,0,0,0) (\Ahat_{\xi_\alpha}\ast
\widehat{\overline{A}}_{\xi_\beta}\ast \Ahat_{\xi_\gamma})\left(\tfrac{k}{\eps},T\right) \right|dk,
\end{aligned}
$$
where 
$$\xi_\delta =\text{ ``+'' if  }\delta=n_*\text{  and  }\xi_\delta =\text{ ``-'' if  }\delta=n_*+1.$$
For $I_1$ we use the $C^2$-smoothness assumption on  $\tilde{\omega}_\pm(k)$ in (H1) and have
\begin{align}
\eps^{-1/2}I_1(T)&=\eps^{1/2}\int_{-\eps^{-1/2}}^{\eps^{-1/2}}|(\vartheta_{n_*}(\eps K)-\omega_0-c_g K)\Atil_{n_*}(K,T)|dK \notag\\
&\leq c \eps^{5/2}\int_{-\eps^{-1/2}}^{\eps^{-1/2}} K^2 |\Ahat_+(K,T)|dK\leq c\eps^{5/2}\|\Ahat_+(\cdot,T)\|_{L^1_2(\R)}.\label{E:I1_est}
\end{align}
For $I_2$ we use the Lipschitz property in (H1). Hence there is $L>0$ such that
$$
\begin{aligned}
&|M^{(1)}_{n_*,n_*}(k)-M^{(1)}_{n_*,n_*}(0)| \\
&\leq |\langle W_Q^{(1)}(q_{n_*}(\cdot, k)-q_{n_*}(\cdot, 0)), q_{n_*}(\cdot, k)\rangle_Q | + |\langle W_Q^{(1)}q_{n_*}(\cdot, 0), (q_{n_*}(\cdot, k)-q_{n_*}(\cdot, 0))\rangle_Q|\\
& \leq L|k| \|W^{(1)}_Q\|_{L^\infty}\left(\|q_{n_*}(\cdot,k)\|_{L^2((0,Q))}+\|q_{n_*}(\cdot,0)\|_{L^2((0,Q))}\right)= 2L|k| \|W^{(1)}_Q\|_{L^\infty}
\end{aligned}
$$
for all $|k|\leq \eps^{1/2}$. Hence
$$\eps^{1/2}I_2(T)\leq 2L\|W^{(1)}_Q\|_{L^\infty}\eps^{5/2} \int_{-\eps^{-1/2}}^{\eps^{-1/2}}|K||\Ahat_+(K,T)|dK \leq c \eps^{5/2} \|\Ahat_+(\cdot,T)\|_{L^1_1(\R)}.$$
Similarly
$$\eps^{1/2}I_3(T)\leq c \eps^{5/2} \|\Ahat_-(\cdot,T)\|_{L^1_1(\R)}.$$
For $I_4$ we first write
$$
\begin{aligned}
&J_{\alpha\beta\gamma}\leq \eps\int_{-3\eps^{-1/2}}^{3\eps^{-1/2}}\int_{-2\eps^{-1/2}}^{2\eps^{-1/2}}\int_{-\eps^{-1/2}}^{\eps^{-1/2}} \left| b_{\alpha\beta\gamma}^{(n_*)}(\eps K,\eps(K- s),\eps(s-l),\eps l) - b_{\alpha\beta\gamma}^{(n_*)}(0,0,0,0)\right| \times \\
&\left| \Atil_\alpha(K-s)\tilde{\overline{A}}_\beta(s-l)\Atil_\gamma(l)\right| dl ds dK \\
&+ \eps |b_{\alpha\beta\gamma}^{(n_*)}(0,0,0,0)| \|\Atil_\alpha \ast \tilde{\overline{A}}_\beta\ast \Atil_\gamma-\Ahat_{\xi_\alpha}\ast
\widehat{\overline{A}}_{\xi_\beta}\ast \Ahat_{\xi_\gamma}\|_{L^1(-3\eps^{-1/2},3\eps^{-1/2})},
\end{aligned}
$$
where we have left out the $T$-dependence of $\Atil_\alpha$ and $\Ahat_\alpha$ for brevity. By the triangle inequality
$$
\begin{aligned}
&\left| b_{\alpha\beta\gamma}^{(n_*)}(k,\lambda,\mu,\nu) - b_{\alpha\beta\gamma}^{(n_*)}(0,0,0,0)\right|\leq \left| b_{\alpha\beta\gamma}^{(n_*)}(k,\lambda,\mu,\nu) -b_{\alpha\beta\gamma}^{(n_*)}(0,\lambda,\mu,\nu)\right| \\
&+\left|b_{\alpha\beta\gamma}^{(n_*)}(0,\lambda,\mu,\nu) - b_{\alpha\beta\gamma}^{(n_*)}(0,0,\mu,\nu)\right|+\dots+\left|b_{\alpha\beta\gamma}^{(n_*)}(0,0,0,\nu)-b_{\alpha\beta\gamma}^{(n_*)}(0,0,0,0)\right|.
\end{aligned}
$$
The differences on the right hand side can be estimated using the Lipschitz continuity in $k$ of Bloch waves, the Cauchy-Schwarz inequality and the algebra property of $H^1(\R)$. For instance, for the first difference we have
$$
\begin{aligned}
&\left| b_{\alpha\beta\gamma}^{(n_*)}(\eps K,\lambda,\mu,\nu) -b_{\alpha\beta\gamma}^{(n_*)}(0,\lambda,\mu,\nu)\right| = \left|\langle q_\alpha(\cdot,\lambda) \overline{q}_\beta(\cdot,\mu)q_\gamma(\cdot,\nu), q_{n_*}(\cdot,\eps K)-q_{n_*}(\cdot,0)\rangle_Q\right|\\
&\leq \eps L|K|\|q_{\alpha}(\cdot,\lambda)\overline{q}_{\beta}(\cdot,\mu)q_{\gamma}(\cdot,\nu)\|_{L^2(0,Q)}\\
&\leq \eps L|K| \|q_{\alpha}(\cdot,\lambda)\|_{H^1(0,Q)} \|q_{\beta}(\cdot,\mu)\|_{H^1(0,Q)}\|q_{\gamma}(\cdot,\nu)\|_{H^1(0,Q)}.
\end{aligned}
$$
Further, to estimate $\|\Atil_\alpha \ast \tilde{\overline{A}}_\beta\ast \Atil_\gamma-\Ahat_{\xi_\alpha}\ast
\widehat{\overline{A}}_{\xi_\beta}\ast \Ahat_{\xi_\gamma}\|_{L^1(-3\eps^{-1/2},3\eps^{-1/2})}$, we set
$$\Atil_\alpha(K):= \Ahat_{\xi_\alpha}(K)+e_{\xi_\alpha}(K), \text{  where } e_\pm(K):=(\chi_{[-\eps^{-1/2},\eps^{-1/2}]}(K)-1)\Ahat_\pm(K),$$
and similarly for $\Atil_\beta$ and $\Atil_\gamma$.
Since $e_\pm$ satisfy
\beq\label{E:epm_est}
\|e_\pm\|_{L^1(\R)}\leq \sup_{|K|>\eps^{-1/2}}|(\chi_{[-\eps^{-1/2},\eps^{-1/2}]}(K)-1)(1+|K|)^{-s_A}|\|\Ahat_\pm\|_{L^1_{s_A}(\R)}\leq c \eps^{s_A/2}\|\Ahat_\pm\|_{L^1_{s_A}(\R)},
\eeq
and because $\|\Atil_\alpha \ast \tilde{\overline{A}}_\beta\ast \Atil_\gamma-\Ahat_{\xi_\alpha}\ast
\widehat{\overline{A}}_{\xi_\beta}\ast \Ahat_{\xi_\gamma}\|_{L^1}\leq \|e_{\xi_\alpha} \ast \widehat{\overline{A}}_{\xi_\beta}\ast \Ahat_{\xi_\gamma}\|_{L^1} +\|\Ahat_{\xi_\alpha}\ast
\overline{e_{\xi_\beta}}(-\cdot)\ast \Ahat_{\xi_\gamma}\|_{L^1}+\|\Ahat_{\xi_\alpha}\ast
\widehat{\overline{A}}_{\xi_\beta}\ast e_{\xi_\gamma}\|_{L^1}  + $ terms quadratic or cubic in $e_\pm$, we get, using Young's inequality for convolutions,
$$\|\Atil_\alpha \ast \tilde{\overline{A}}_\beta\ast \Atil_\gamma-\Ahat_{\xi_\alpha}\ast
\widehat{\overline{A}}_{\xi_\beta}\ast \Ahat_{\xi_\gamma}\|_{L^1(-3\eps^{-1/2},3\eps^{-1/2})}\leq c\eps^{s_A/2}\|\Ahat_{\xi_\alpha}\|_{L^1_{s_A}(\R)}\|\Ahat_{\xi_\beta}\|_{L^1_{s_A}(\R)}\|\Ahat_{\xi_\gamma}\|_{L^1_{s_A}(\R)}.$$
Hence, we arrive at
$$\eps^{1/2}I_4(T)\leq c\eps^{5/2}\sum_{\xi,\zeta,\theta\in \{+,-\}}\|\Ahat_{\xi}(\cdot,T)\|_{L^1_{s_A}(\R)}\|\Ahat_{\zeta}(\cdot,T)\|_{L^1_{s_A}(\R)}\|\Ahat_{\theta}(\cdot,T)\|_{L^1_{s_A}(\R)}$$ 
because $s_A \geq 2$.

Similarly, one can show that all the formal $O(\eps^{1/2})$ terms in $\Res_{n_*+1}$ are $O(\eps^{5/2})$ in the $L^1(\B_Q)$-norm.

Because in $\Res_n, n \in \N\setminus\{n_*,n_{*+1}\}$ the whole formal $O(\eps^{1/2})$-part vanishes by the choice of $\Atil_j$ and $\Atil_{j,m}$, it remains to discuss the h.o.t. terms in \eqref{E:hot}.

Firstly, for $s\in (0,3/2)$
$$\|\eps \pa_T \vec{U}^{1,\text{ext}}\|_{\cX(s)}  \leq c\eps^{5/2},$$
where $c=c(\max_{T\in [0,T_0]}\|\pa_T \Ahat_\pm(\cdot,T)\|_{L^1(\R)},\max_{T\in [0,T_0]}\|\Ahat_\pm(\cdot,T)\|_{L^1(\R)})$ similarly to the proof of Lemma \ref{L:ext_est}. Next, because $M^{(1)}\vec{U}^{1,\text{ext}}=\cD(W_Q^{(1)}\util^{1,\text{ext}})$ with $\util^{1,\text{ext}}:=\cD^{-1}\vec{U}^{1,\text{ext}}$, the isomorphic property of $\cD$ (Lemma \ref{L:D_isom}) yields
\beq\label{E:M1_U1ext}
\|\eps M^{(1)} \vec{U}^{1,\text{ext}}\|_{\cX(s)}\leq c\eps \|W_Q^{(1)} \util^{1,\text{ext}}\|_{L^1(\B_Q,H^s(0,Q))}\leq c\eps \|\util^{1,\text{ext}}\|_{L^1(\B_Q,H^s(0,Q))}\leq c\eps \|\vec{U}^{1,\text{ext}}\|_{\cX(s)},
\eeq
where $c=c(\|\tfrac{d^{\lceil s \rceil}}{dx^{\lceil s \rceil}}W_Q^{(1)}\|_{L^2(\B_Q)})>0$. Hence, by Lemma \ref{L:ext_est},
$$\|\eps M^{(1)}(\cdot) \vec{U}^{1,\text{ext}}(\cdot,t)\|_{\cX(s)}\leq c\eps^{5/2} \quad \text{for all  }t\in [0,T_0\eps^{-1}].$$

Similarly
\beq\label{E:MR_U1ext}
\begin{aligned}
\left\|\eps \sum_{m\in \Z_R}M^{(R,m)}(\cdot) \vec{U}^{1,\text{ext}}(\cdot-mk_W,t)\right\|_{\cX(s)}&\leq c \eps \sum_{m\in \Z_R}|a_m|\|\util^{1,\text{ext}}(\cdot,\cdot-mk_W,t)\|_{L^1(\B_Q,H^s(0,Q))}  \\
&\leq c \eps \sum_{m\in \Z_R} \|\vec{U}^{1,\text{ext}}(\cdot-mk_W,t)\|_{\cX(s)}m^{-2}\\
&\leq c\eps^{5/2} \quad \text{for all  }t\in [0,T_0\eps^{-1}]
\end{aligned}
\eeq
where we have used the decay $|a_m|\leq cm^{-2}$ for $W\in C^2$, the identity $\|\vec{U}^{1,\text{ext}}(\cdot-mk_W,t)\|_{\cX(s)}=\|\vec{U}^{1,\text{ext}}(\cdot,t)\|_{\cX(s)}$ due to the $1/N$-periodicity of $\Uvec^{1,\text{ext}}$ in $k$, and Lemmas \ref{L:D_isom} and \ref{L:ext_est}.

Finally,
$$
\begin{aligned}
\|\vec{F}(\UvecOext,\vec{U}^{1,\text{ext}},\UvecOext)\|_{\cX(s)}&\leq c|\sigma| \|\util^{0,\text{ext}}*_{\B_Q}\overline{\util}^{1,\text{ext}}*_{\B_Q}\util^{0,\text{ext}}\|_{L^1(\B_Q,H^s(0,Q))} \\
& \leq c \|\util^{1,\text{ext}}\|_{L^1(\B_Q,H^s(0,Q))}\|\util^{0,\text{ext}}\|^2_{L^1(\B_Q,H^s(0,Q))}\\
&\leq c \|\vec{U}^{1,\text{ext}}\|_{\cX(s)} \|\UvecOext\|_{\cX(s)}^2 \leq c \eps^{5/2} \quad \text{for all  }t\in [0,T_0\eps^{-1}]
\end{aligned}
$$
using the algebra property in Lemma \ref{L:algeb_L1Hs} for $s>1/2$ and, again, Lemma \ref{L:ext_est}. Other nonlinear terms are treated analogously.

We conclude that there is $C_{\text{Res}}>0$ and $\eps_0>0$ such that $\|\vec{\Res}(\cdot,t)\|_{\cX(s)}\leq C_{\text{Res}}\eps^{5/2}$ for all $\eps\in(0,\eps_0)$ and all $t\in [0,T_0\eps^{-1}]$. 
\epf


Finally, we complete the proof of Theorem \ref{T:main}. Writing the exact solution of \eqref{E:PNLS_U} as a sum of the extended ansatz $\Uvecext $ and the error $\Evec$, we have $\Uvec = \Uvecext + \Evec$. In \eqref{E:PNLS_U} this produces
\beq\label{E:err_eq}
 \partial_t \Evec = -\ri \Theta(k) \Evec +\ri \vec{G}(\Uvecext,\Evec)
\eeq
with 
$$
\begin{aligned}
\vec{G}(\Uvecext,\Evec)=&\vec{\Res}+\vec{F}(\Uvec,\Uvec,\Uvec) -\vec{F}(\Uvecext,\Uvecext,\Uvecext)\\ 
&-\varepsilon \left(M^{(1)}(k) \Evec +  \sum \limits_{m \in \mathbb{Z}_R} M^{(R,m)}(k) \Evec(k-mk_W,t)\right).
\end{aligned}
$$
Due to the cubic structure of $\vec{F}$, Lemma \ref{L:ext_est}, the algebra property of $L^1(\B_Q,H^s(0,Q))$ and the isomorphism $\cD$ we have for $s>1/2$ the existence of $c>0$ such that 
\beq\label{E:F_est}
\|\vec{F}(\Uvec,\Uvec,\Uvec) -\vec{F}(\Uvecext,\Uvecext,\Uvecext)\|_{\cX(s)} \leq c(\|\Uvecext\|^2_{\cX(s)}\|\Evec\|_{\cX(s)}+\|\Uvecext\|_{\cX(s)}\|\Evec\|^2_{\cX(s)}+\|\Evec\|^3_{\cX(s)}).
\eeq

Similarly to \eqref{E:M1_U1ext} and \eqref{E:MR_U1ext} we get
\beq\label{E:M_est}
\| \varepsilon ( M^{(1)}(\cdot) \Evec(\cdot,t) +  \sum \limits_{m \in \mathbb{Z}_R} M^{(R,m)}(\cdot)  \Evec(\cdot-mk_W,t) )\|_{\cX(s)} \leq c \varepsilon \| \Evec(\cdot, t)\|_{\cX(s)}.
\eeq
Next, because $\Uvecext=\UvecOext+\vec{U}^{1,\text{ext}}$, Lemma \ref{L:ext_est} implies 
\beq\label{E:Uext_est}
\|\Uvecext(\cdot,t)\|_{\cX(s)}\leq c\eps^{1/2} \qquad \text{for all } t\in [0,\eps^{-1}T_0]. 
\eeq

Combining \eqref{E:F_est}, \eqref{E:M_est}, \eqref{E:Uext_est} and Lemma \ref{L:res_est} provides the estimate 
$$\|\vec{G}(\Uvecext,\Evec)\|_{\cX(s)} \leq c_1\eps \|\Evec\|_{\cX(s)}+c_2\eps^{1/2} \|\Evec\|_{\cX(s)}^2+c_3\|\Evec\|_{\cX(s)}^3+C_{\text{Res}}\eps^{5/2}$$
on the time interval $[0,\eps^{-1}T_0]$ with some $c_1,c_2,c_3>0$ independent of $\eps$ and $t$.

The operator $-\ri\Theta(k)$ generates a strongly continuous unitary group $S(t)=e^{-\ri \Theta t}: \cX(s) \rightarrow \cX(s)$ and equation \eqref{E:err_eq} reads
$$
\Evec(t) = \Evec(0)+\int \limits_0^t S(t-\tau)\vec{G}(\Uvecext,\Evec)(\tau) ~d\tau.
$$
With the above estimates we get
$$
\begin{aligned}
\| \Evec(t) \|_{\cX(s)} \leq \| \Evec(0) \|_{\cX(s)}+\int \limits^t_0 c_1 \varepsilon \| \Evec(\tau) \|_{\cX(s)} + c_2 \varepsilon^{1/2} \| \Evec(\tau) \|^2_{\cX(s)}
+c_3 \| \Evec(\tau) \|^3_{\cX(s)} + C_{\text{Res}} \varepsilon^{5/2} ~d\tau.
\end{aligned}
$$
Because $\Uvec(0)=\Uvec^\text{app}(0)$, we get $\Evec(0)=\Uvec^\text{app}(0)-\Uvecext(0)$ and Lemma \ref{L:diff_app_ext} provides $\|\Evec(0)\|_{\cX(s)}\leq C_0\eps^{3/2}$ for all $s\in (0,3/2)$ and $\eps\in (0,\eps_0)$ with some $\eps_0>0$.

Given an $M>C_0$ there exists $T>0$ such that $\|\Evec(t)\|_{\cX(s)}\leq M \eps^{3/2}$ for all $t\in [0,T]$. Next, we use Gronwall's lemma and a bootstrapping argument to choose $\eps_0>0$ and $M>0$ such that if $\eps\in (0,\eps_0)$, then $\|\Evec(t)\|_{\cX(s)}\leq M \eps^{3/2}$ for all $t\in [0,\eps^{-1}T_0]$.

If $\|\Evec(t)\|_{\cX(s)} \leq M\eps^{3/2}$, then 
$$ 
\begin{aligned}
\| \Evec(t) \|_{\cX(s)} & \leq ~C_0 \eps^{3/2}+\int \limits^t_0 c_1 \varepsilon \| \Evec(\tau) \|_{\cX(s)}  ~d\tau + 
t\left( c_2 \varepsilon^{7/2} M^2 +c_3 \varepsilon^{9/2} M^3  +  C_{\text{Res}} \varepsilon^{5/2} \right)\\
& \leq ~ \varepsilon^{3/2} \left[ C_0 + t \eps\left( c_2 \varepsilon M^2 + c_3 \varepsilon^{2}  M^3 +  C_{\text{Res}} \right) \right] e^{c_1 \varepsilon t},
\end{aligned}
$$
where the second inequality follows from Gronwall's lemma. In order to achieve the desired estimate on $t\in [0,\eps^{-1}T_0]$, we redefine
$$M:=C_0+T_0(C_\text{Res}+1)e^{c_1T_0}$$
and choose $\eps_0$ so small that $ c_2 \varepsilon_0 M^2 + c_3 \varepsilon_0^{2}  M^3 \leq 1$. Then, clearly,
$$
\sup \limits_{t \in [0,\eps^{-1}T_0]} \| \Evec(t) \|_{\cX(s)} \leq M \varepsilon^{3/2}. 
$$
Using Lemma \ref{L:diff_app_ext} and the triangle inequality produces $\sup_{t\in[0,\eps^{-1}T_0]}\|\Uvec(\cdot,t)-\Uvecapp(\cdot,t)\|_{\cX(s)} \leq c \eps^{3/2}$. The estimate in the $C^0_b$-norm and the decay for $|x|\to \infty$ follow from Lemmas \ref{L:sup_control} and \ref{L:D_isom} since $s>1/2$. This completes the proof of Theorem \ref{T:main} for the case of $k_\pm \in \Q$. \hfill $\square$

\blem\label{L:diff_app_ext}
Assume (H1) and let $\Ahat_\pm(\cdot,T) \in L^1_{s_A}(\R)\cap L^2(\R)$ with $s_A\geq 2$ for all $T\in [0,T_0]$. Then there exist $c=c(\max_{T\in [0,T_0]}\|\Ahat_+(\cdot,T)\|_{L^1_{s_A}(\R)},\max_{T\in [0,T_0]}\|\Ahat_-(\cdot,T)\|_{L^1_{s_A}(\R)})$ and $\eps_0>0$  such that for $u_\text{app}$ and $\Uvecext$ given by  \eqref{E:uapp_rat} and \eqref{E:Uext} respectively we have
$$\|\Uvecext(\cdot,t)-\Uvec^\text{app}(\cdot,t)\|_{\cX(s)}\leq c \eps^{3/2}$$
for all $s\in (0,3/2), \eps \in (0,\eps_0)$ and all $t\in [0,\eps^{-1}T_0]$.
\elem
\bpf
Because $\Ahat_\pm(\cdot,T) \in L^2(\R)$, it is also $A_\pm(\cdot,T)\in L^2(\R)$ and hence $u^\text{app}(\cdot,t)\in L^2(\R)$ and one can apply $\cD\cT$ to $u^\text{app}$ producing $\Uvec^\text{app}$ with
$$U^\text{app}_n(k,t)=\eps^{-1/2}e^{-\ri \omega_0 t} \sum_\pm \Ahat_\pm\left(\tfrac{k}{\eps},\eps t\right)\pi_n^\pm(k), \quad \text{where } \pi_n^\pm(k):=\langle q_\pm(\cdot), q_n(\cdot,k)\rangle_Q.$$
Similarly to the decomposition of $\Uvecext$ in \eqref{E:Uext_decomp} and \eqref{E:Uext_decomp2} we write
$$\Uvecapp=\UvecOapp+\vec{U}^{1,\text{app}},$$ 
where
$$\UvecOapp(k,t):=\eps^{-1/2}e^{-\ri \omega_0 t}\left(\Ahat_{+}(K,T)\pi_{n_*}^+(k)e_{n_*}+\Ahat_{-}(K,T)\pi_{n_*+1}^-(k)e_{n_*+1}\right), \ K=\eps^{-1}k.$$
Since $\|\vec{U}^{1,\text{ext}}(\cdot,t)\|_{\cX(s)} \leq c\eps^{3/2}$ for all $t\in [0,\eps^{-1}T_0]$ (see Lemma \ref{L:ext_est}), it remains to show that $\|\vec{U}^{1,\text{app}}(\cdot,t)\|_{\cX(s)} \leq c\eps^{3/2}$ and $\|(\UvecOext-\UvecOapp)(\cdot,t)\|_{\cX(s)}\leq c\eps^{3/2}$ for all $t\in [0,\eps^{-1}T_0]$.

For $\|\vec{U}^{1,\text{app}}(\cdot,t)\|_{\cX(s)}$ note first that
$$
\begin{aligned}
&U^{1,\text{app}}_{n_*}(k,t)=\eps^{-1/2}e^{-\ri \omega_0 t} \Ahat_-(K,T)\pi^-_{n_*}(k), \ U^{1,\text{app}}_{n_*+1}(k,t)=\eps^{-1/2}e^{-\ri \omega_0 t} \Ahat_+(K,T)\pi^+_{n_*+1}(k),\\
&U^{1,\text{app}}_{n}(k,t)=\eps^{-1/2}e^{-\ri \omega_0 t} \sum_\pm \Ahat_\pm(K,T)\pi^\pm_{n}(k) \quad \text{for } n \notin \{n_*,n_*+1\}.
\end{aligned}
$$
Using the Lipschitz continuity in (H1), there is some $L>0$, such that, for instance, $|\pi^-_{n_*}(k)|\leq |\langle q_-(\cdot),q_{n_*}(\cdot,0)\rangle_Q|+L |k|=L|k|$ because $q_-(x)=q_{n_*+1}(x,0)$. Hence
$$\|U^{1,\text{app}}_{n_*}(\cdot,t)\|_{L^1(\B_Q)}\leq L\eps^{3/2} \int_\R |\Ahat_-(K,T)||K| dK \leq c \eps^{3/2} \|\Ahat_-(\cdot,T)\|_{L^1_1(\R)}.$$
Similarly, $\|U^{1,\text{app}}_{n_*+1}(\cdot,t)\|_{L^1(\B_Q)}\leq c\eps^{3/2} \|\Ahat_+(\cdot,T)\|_{L^1_1(\R)}$. Because 
$$
\begin{aligned}
\|\vec{U}^{1,\text{app}}(\cdot,t)\|_{\cX(s)} \leq & c\left(\sum_{n\in \{n_*,n_*+1\}}\|U^{1,\text{app}}_{n}(\cdot,t)\|_{L^1(\B_Q)} \phantom{\left\|\left(\sum_{n\in \N\setminus\{n_*,n_*+1\}}n^{2s}|U^{1,\text{app}}_n(\cdot,t)|^2\right)^{1/2}\right\|_{L^1(\B_Q)}}\right.\\
& \left. + \left\|\left(\sum_{n\in \N\setminus\{n_*,n_*+1\}}n^{2s}|U^{1,\text{app}}_n(\cdot,t)|^2\right)^{1/2}\right\|_{L^1(\B_Q)}\right),
\end{aligned}
$$
it remains to consider $|U^{1,\text{app}}_n|$ for $n\in \N\setminus\{n_*,n_*+1\}$. For $\pi^+_n(k)$ we have
\begin{align}
\pi^+_n(k)&=\langle q_{n_*}(\cdot,k),q_n(\cdot,k)\rangle_Q+\langle q_{n_*}(\cdot,0)-q_{n_*}(\cdot,k),q_n(\cdot,k)\rangle_Q =\langle q_{n_*}(\cdot,0)-q_{n_*}(\cdot,k),q_n(\cdot,k)\rangle_Q \notag\\
&=\frac{1}{\vartheta_n(k)}\langle \cL(\cdot,k)(q_{n_*}(\cdot,0)-q_{n_*}(\cdot,k)),q_n(\cdot,k)\rangle_Q.\label{E:pi+_L}
\end{align}
The Lipschitz continuity in (H1) now provides 
$$\|\cL(\cdot,k)(q_{n_*}(\cdot,0)-q_{n_*}(\cdot,k))\|_{L^2(0,Q)}\leq \|q_{n_*}(\cdot,0)-q_{n_*}(\cdot,k)\|_{H^2(0,Q)}\leq L|k|$$
such that (using \eqref{E:band_as})
$$|\pi^+_n(k)|\leq \frac{c}{n^2}|k| \quad \text{for all } n\in \N \ \text{and } k\in \B_Q$$
and similarly for $|\pi^-_n|$. As a result $\left\|\left(\sum_{n\in \N\setminus\{n_*,n_*+1\}}n^{2s}|U^{1,\text{app}}_n(\cdot,t)|^2\right)^{1/2}\right\|_{L^1(\B_Q)} $ can be estimated by
$$c\eps^{1/2}\left(\sum_{n\in \N}n^{2s-4}\right)^{1/2}\sum_\pm\int_{\B_Q}|\Ahat_\pm\left(\tfrac{k}{\eps},T\right)||\tfrac{k}{\eps}|dk \leq c\eps^{3/2}\sum_\pm \|\Ahat_\pm(\cdot,T)\|_{L^1_1(\R)},$$
where the last inequality uses $s<3/2$. In summary 
$$\|\vec{U}^{1,\text{app}}(\cdot,t)\|_{\cX(s)} \leq c \eps^{3/2} \sum_\pm \|\Ahat_\pm(\cdot,T)\|_{L^1_1(\R)}.$$

For $\|(\UvecOext-\UvecOapp)(\cdot,t)\|_{\cX(s)}$  we have
$$
\begin{aligned}
\|(U^{0,\text{ext}}_{n_*}-& U^{0,\text{app}}_{n_*})(\cdot,t)\|_{L^1(\B_Q)} \\
&\leq \eps^{-1/2}\left(\int_{-\eps^{1/2}}^{\eps^{1/2}}|\Ahat_+\left(\tfrac{k}{\eps},\eps t\right)||1-\pi_{n_*}^+(k)|dk+\int_{\B_Q\setminus (-\eps^{1/2},\eps^{1/2})}|\Ahat_+\left(\tfrac{k}{\eps},\eps t\right)||\pi_{n_*}^+(k)|dk\right).
\end{aligned}
$$
Once again, by the Lipschitz continuity it is $|1-\pi_{n_*}^+(k)|=|\langle q_+(\cdot),q_{n_*}(\cdot,0)-q_{n_*}(\cdot,k)\rangle_Q|\leq L|k|$. Clearly, also $|\pi_{n_*}^+(k)|\leq 1$ for all $k\in \B_Q$. Hence
$$
\begin{aligned}
\|(U^{0,\text{ext}}_{n_*}-U^{0,\text{app}}_{n_*})(\cdot,t)\|_{L^1(\B_Q)}&\leq c\eps^{3/2}\|\Ahat_+(\cdot,\eps t)\|_{L^1_1(\R)} + c\eps^{1/2}\|\Ahat_+(\cdot,\eps t)\|_{L^1(\R\setminus (-\eps^{-1/2}, \eps^{-1/2}))}\\
&\leq c\eps^{3/2}\|\Ahat_+(\cdot,\eps t)\|_{L^1_1(\R)} + c\eps^{1/2+s_A/2}\|\Ahat_+(\cdot,\eps t)\|_{L^1_{s_A}(\R)},
\end{aligned}
$$
where in the second step we have used $\|\Ahat_+(\cdot,\eps t)\|_{L^1(\R\setminus (-\eps^{-1/2}, \eps^{-1/2}))}\leq c\eps^{s_A/2}\|\Ahat_+(\cdot,\eps t)\|_{L^1_{s_A}(\R)}$, see \eqref{E:epm_est}. Similarly, one gets 
$$
\|(U^{0,\text{ext}}_{n_*+1}-U^{0,\text{app}}_{n_*+1})(\cdot,t)\|_{L^1(\B_Q)}\leq c\eps^{3/2}\|\Ahat_-(\cdot,\eps t)\|_{L^1_1(\R)} + c\eps^{1/2+s_A/2}\|\Ahat_-(\cdot,\eps t)\|_{L^1_{s_A}(\R)}.
$$
For $s_A\geq 2$ is $1/2+s_A/2\geq 3/2$ and the lemma is proved.
\epf


\subsection{Proof of Theorem \ref{T:main} for Irrational $k_\pm$}\label{S:irrational_pf}

The method of proof in the case of irrational $k_\pm$ is the same as in the rational case. The main difference is in the choice of the extended ansatz $\Uvec^\text{ext}(k,t)$. Therefore, we concentrate on explaining the choice of $\Uvec^\text{ext}(k,t)$ and describe where the proof differs from that in Section \ref{S:rational_pf}.

When $k_\pm \notin \Q$, then clearly case (a) in Section \ref{S:Bloch} applies, i.e. we have simple Bloch eigenvalues at $k=k_+=k_0$ and $k=k_-=-k_0$ for some $k_0\in (0,1/2)\setminus \Q$. The Bloch eigenfunctions are 
$$p_\pm(x)=p_{n_0}(x,\pm k_0).$$
Because $k_0\notin \Q$, there is no common period of the Bloch waves $p_\pm(x)e^{\pm \ri k_0 x}$ and of $V$. Hence, we use the period $P=2\pi$ corresponding to $V$ as the working period with the corresponding Brillouin zone $\B_{2\pi}=(-1/2,1/2]$. Note that the effective coupled mode equations are now \eqref{E:CME} with $\beta=\gamma=0$.

Similarly to the beginning of Section \ref{S:rational_pf}, in order to motivate the choice of the extended ansatz for the approximate solution $\Uvec^\text{ext}$, we study first the residual of the formal approximate ansatz $u_\text{app}$ with $\text{supp}(\hat{A}_\pm(\cdot,T))\subset [-\eps^{-1/2},\eps^{-1/2}]$. We have
$$
\begin{aligned}
\cT(\uapp)(x,k,t)=&\eps^{-1/2}\sum_{\pm}p_\pm(x)\sum_{\eta\in \Z}\hat{A}_\pm\left(\tfrac{k\mp k_0+\eta}{\eps},\eps t\right)e^{\ri(\eta x- \omega_0 t)}\\
=&\eps^{-1/2}\sum_{\pm}p_\pm(x)\hat{A}_\pm\left(\tfrac{k\mp k_0}{\eps},\eps t\right)e^{-\ri \omega_0 t}, \quad k \in \B_{2\pi}
\end{aligned}
$$
because $\eps^{-1}(k\mp k_0+\eta) \in \text{supp}(\hat{A}_\pm(\cdot,T))$ for some $k\in \B_Q$ and $\eta\in \Z$ and with $k_0 \in (0,1/2)$ is possible only if $\eta=0$. Hence, if $\text{supp}(\hat{A}_\pm(\cdot,T))\subset [-\eps^{-1/2},\eps^{-1/2}]$, then for $k\in\B_{2\pi}$
$$\begin{aligned}
&\cT(\text{PNLS}(\uapp))(x,k,t)= \eps^{1/2}\sum_{\pm}\left[\ri \pa_T\hat{A}_\pm\left(\tfrac{k\mp k_0}{\eps},T\right)p_\pm(x)+2(\ri k_0p_\pm +p_\pm')\ri \tfrac{k\mp k_0}{\eps}\hat{A}_\pm\left(\tfrac{k\mp k_0}{\eps},T\right)\phantom{\sum_{m\in \Z_3^\pm,\eta \in \cS_m^\pm}}\right.\\
&\left.-W^{(1)}(x)p_\pm(x)\hat{A}_\pm\left(\tfrac{k\mp k_0}{\eps},T\right)-W^{(2)}_\pm p_\pm(x)\hat{A}_\pm\left(\tfrac{k\pm k_0}{\eps},T\right)\right.\\
&\left.- p_\pm(x)\sum_{m\in \Z_3^\pm,\eta \in \cS_m^\pm} a_m\hat{A}_\pm\left(\frac{k\mp k_0-mk_W+\eta}{\eps},T\right)e^{\ri \eta x}\right]e^{-\ri \omega_0 t}\\
& -\eps^{1/2}  e^{-\ri \omega_0 t}\sigma(x) \sum_{\stackrel{\xi,\zeta,\theta\in \{+,-\}}{\eta\in \cS_{\xi,\zeta,\theta}}}p_{\xi}(x)\overline{p_{\zeta}}(x)p_{\theta}(x)\left(\hat{A}_{\xi}\ast \hat{\overline{A}}_{\zeta}\ast \hat{A}_{\theta}\right)\left(\tfrac{k-(\xi k_0-\zeta k_0+\theta k_0) +\eta}{\eps},T\right)e^{\ri \eta x}\\
&  - \eps^{3/2}e^{-\ri \omega_0 t}\sum_{\pm}p_\pm(x)\left(\tfrac{k\mp k_0}{\eps}\right)^2\hat{A}_\pm\left(\tfrac{k\mp k_0}{\eps},T\right),
\end{aligned}
$$
where 
$$\begin{aligned}
&T = \eps t, \quad \Z_3^\pm := \{m\in \Z\setminus\{0\}: a_m\neq 0, mk_W\notin \Z, mk_W\pm 2k_0 \notin \Z\},\\ 
& \cS_m^\pm:=\{\eta\in\Z: \pm k_0 +mk_W-\eta\in \overline{\B_{2\pi}}=[-\tfrac{1}{2},\tfrac{1}{2}]\}, \\
&\cS_{\xi,\zeta,\theta}:=\{\eta\in \Z: \xi k_0-\zeta k_0+\theta k_0-\eta\in \overline{\B_{2\pi}}\} \ \text{for} \  \xi,\zeta,\theta\in\{+,-\}.
\end{aligned}
$$ 
It is $\cS_{\xi,\zeta,\theta}\subset \{0,1\}$ or $\cS_{\xi,\zeta,\theta}\subset \{0,-1\}$. 

In the $k-$variable the support of $\cT(\text{PNLS}(\uapp))$ within $\B_{2\pi}$ consists of intervals (with radius at most $3\eps^{1/2}$) centered at $k_0,-k_0, k_0+mk_W$ with $m\in \Z_3^{+}$, $-k_0+mk_W$ with $m\in \Z_3^{-}$, and at $3k_0$ and $-3k_0$ as well as at  integer shifts of these points. Note that within $\B_{2\pi}$ there are only finitely many support intervals because $\pm k_0+mk_W+\Z, m\in \Z_3^\pm$ generates only finitely many distinct points in $\B_{2\pi}$ due to assumption (H2). Similarly to \eqref{E:KR} we define these sets of points by
$$
\cK_R^\pm:=(\pm k_0+k_W+\Z_3^\pm + \Z)\cap \overline{\B_{2\pi}}=\{k\in \overline{\B_{2\pi}}: k=\pm k_0+mk_W+\eta \text{ for some } m\in \Z_3^\pm, \eta\in \cS_m^\pm\}
$$
and their elements by
$$\cK_R^\pm=\{\kappa^\pm_1,\dots, \kappa^\pm_{J^\pm}\} \text{ with some } J^\pm\in \N.$$

The choice of the splitting in $W$ using $W^{(2)}_+$ and $W^{(3)}_+$ or $W^{(2)}_-$ and $W^{(3)}_-$ is motivated at the beginning of Sec. \ref{S:formal}. The choice is made in order to isolate the parts of $W$ responsible for the coupling of the two modes.

Analogously to \eqref{E:Uext} we are lead to the following extended ansatz for $k\in \B_{2\pi}$
\beq\label{E:Uext_rat}
\begin{aligned}
&U_{n_0}^\text{ext}(k,t) := e^{-\ri \omega_0 t}\sum_\pm\left(\eps^{-1/2}\Atil_{n_0}^\pm\left(\tfrac{k\mp k_0}{\eps},T\right)
+\eps^{1/2}\sum_{j=1}^{J^\pm} \Atil^\pm_{n_0,j}\left(\frac{k-\kappa_j^\pm}{\eps},T\right) \right)\\
&+\eps^{1/2}e^{-\ri \omega_0 t} \left(\sum_{\eta\in \cS_{+,-,+}}\Atil^+_{n_0,NL}\left(\frac{k- 3k_0+\eta}{\eps},T\right)+\sum_{\eta\in \cS_{-,+,-}}\Atil^-_{n_0,NL}\left(\frac{k+3k_0+\eta}{\eps},T\right)\right),\\
&U_{n}^\text{ext}(k,t) := \eps^{1/2}e^{-\ri \omega_0 t}\sum_\pm\left(\Atil_{n}^\pm\left(\tfrac{k\mp k_0}{\eps},T\right)
+\sum_{j=1}^{J^\pm} \Atil^\pm_{n,j}\left(\frac{k-\kappa_j^\pm}{\eps},T\right) \right)\\
&+\eps^{1/2}e^{-\ri \omega_0 t}\left(\sum_{\eta\in \cS_{+,-,+}}\Atil^+_{n,NL}\left(\frac{k- 3k_0+\eta}{\eps},T\right)+\sum_{\eta\in \cS_{-,+,-}}\Atil^-_{n,NL}\left(\frac{k+ 3k_0+\eta}{\eps},T\right)\right)
\end{aligned}
\eeq
for $n\in \N\setminus\{n_0\}$, where
$$
\begin{aligned}
&\supp (\Atil_{n_0}^\pm(\cdot, T)) \cap \eps^{-1}\B_{2\pi}, \quad  \supp (\Atil_{n,j}^\pm(\cdot, T)) \cap \eps^{-1}\B_{2\pi} \subset [-\eps^{-1/2},\eps^{-1/2}],\\
&\supp (\Atil_{m}^\pm(\cdot, T)) \cap \eps^{-1}\B_{2\pi}, \quad  \supp (\Atil_{n,NL}^\pm(\cdot,T)) \cap \eps^{-1}\B_{2\pi} \subset [-3\eps^{-1/2},3\eps^{-1/2}]
\end{aligned}
$$
for all $n\in \N, m\in \N\setminus\{n_0\}$ and $j\in \{1,\dots,J^\pm\}$ and where 
$$\Uvecext\left(k+1,t\right)=\Uvecext(k,t) \ \text{for all} \ k\in \R, t \in \R.$$
Similarly to \eqref{E:Uext_decomp}, \eqref{E:Uext_decomp2} we decompose $\Uvecext = \vec{U}^{0,\text{ext}}+\vec{U}^{1,\text{ext}}$, where 
\beq\label{E:Uext_decomp2_rat}
\begin{aligned}
&\vec{U}^{0,\text{ext}}=\vec{U}_+^{0,\text{ext}}+\vec{U}_-^{0,\text{ext}}, \\ 
&\vec{U}^{0,\text{ext}}_+:= \eps^{-1/2}e_{n_0}\Atil^+_{n_0}\left(\tfrac{k-k_0}{\eps},T\right)e^{-\ri \omega_0 t}, \vec{U}^{0,\text{ext}}_-:=\eps^{-1/2}e_{n_0}\Atil^-_{n_0}\left(\tfrac{k+k_0}{\eps},T\right)e^{-\ri \omega_0 t},\quad \text{and} \\ 
&\vec{U}^{1,\text{ext}}:=\Uvecext-\vec{U}^{0,\text{ext}}.
\end{aligned}
\eeq

In analogy to \eqref{E:PNLS_U} we get
\beq\label{E:PNLS_U_rat}
\left(\ri \pa_t -\Omega(k)-\eps M'^{(1)}(k)\right)\Uvec(k,t) -\eps \sum_{mk_W\notin \Z}M'^{(R,m)}(k) \Uvec(k-mk_W,t)+\vec{F'}(\Uvec,\Uvec,\Uvec)(k,t)=0,
\eeq
where
$$
\begin{aligned}
&\Omega_{j,j}(k):=\omega_j(k), \ \Omega_{i,j}:=0 \text{ if } i\neq j,\\
&M'^{(1)}_{i,j}(k):=\langle W^{(1)}(\cdot)p_j(\cdot,k), p_i(\cdot,k)\rangle_{2\pi}, \quad M'^{(R,m)}_{i,j}(k):=a_m\langle p_j(\cdot,k-mk_W), p_i(\cdot,k)\rangle_{2\pi},\\
&F'_j(\Uvec,\Uvec,\Uvec)(k,t) :=-\langle\sigma(\cdot) (\util\ast_{\B_{2\pi}}\tilde{\overline{u}}\ast_{\B_{2\pi}}\util)(\cdot,k,t),p_j(\cdot,k)\rangle_{2\pi}, \ \util(x,k,t)=\sum_{n\in \N}U_n(k,t)p_n(x,k).
\end{aligned}
$$
When studying the residual near $k=k_0$ or $k=-k_0$, we exploit both ways of splitting $W(x)$ in \eqref{E:W_split}. Namely, we have the following two equivalent reformulations of \eqref{E:PNLS_U_rat}
\beq\label{E:PNLS_U_split1}
\begin{aligned}
&\left(\ri \pa_t -\Omega(k)-\eps M'^{(1)}(k)\right)\Uvec(k,t) -\eps M^{(2_+)}(k)\Uvec(k+2k_0,t)\\
&-\eps \sum_{m \in \Z_3^+}M'^{(R,m)}(k) \Uvec(k-mk_W,t)+\vec{F'}(\Uvec,\Uvec,\Uvec)(k,t)=0
\end{aligned}
\eeq
and 
\beq\label{E:PNLS_U_split2}
\begin{aligned}
&\left(\ri \pa_t -\Omega(k)-\eps M'^{(1)}(k)\right)\Uvec(k,t) -\eps M^{(2_-)}(k)\Uvec(k-2k_0,t)\\
&-\eps \sum_{m \in \Z_3^-}M'^{(R,m)}(k) \Uvec(k-mk_W,t)+\vec{F'}(\Uvec,\Uvec,\Uvec)(k,t)=0,
\end{aligned}
\eeq
where
$$
\begin{aligned}
&M^{(2_\pm)}_{i,j}(k):=\langle W^{(2)}_\pm p_j(\cdot,k\pm 2k_0), p_i(\cdot,k)\rangle_{2\pi}.
\end{aligned}
$$
As explained in Sec. \ref{S:formal}, the splitting of $W$ with $W^{2_-}$ will be used near $k=k_0$ because it extracts the part of $W$ which shifts $\Uvec^\text{ext}(k,t)$ in $k$ by $2k_0$ to the right and thus produces the linear $\Atil_{n_0}^-$-term in the residual near $k=k_0$. Similarly, the splitting with $W^{2_+}$ will be used near $k=-k_0$. 

The residual of $\Uvec^\text{ext}$ on $k\in (\pm k_0-3\eps^{1/2},\pm k_0+3\eps^{1/2})$ is given by 
$$
\begin{aligned}
&\Res_{n_0}(k,t)=\\ 
&\eps^{1/2}\left[\left(\ri \pa_T -\eps^{-1}(\omega_{n_0}(k)-\omega_0)-M'^{(1)}_{n_0,n_0}(k)\right)\Atil^\pm_{n_0}\left(\tfrac{k\mp k_0}{\eps},T\right) - M^{(2_\mp)}_{n_0,n_0}(k)\Atil^\mp_{n_0}\left(\tfrac{k\mp k_0}{\eps},T\right) \right]e^{-\ri \omega_0 t}\\
&+ F'_{n_0}(\UvecOext_\pm,\UvecOext_\pm,\UvecOext_\pm)(k,t)+2F'_{n_0}(\UvecOext_\mp,\UvecOext_\mp,\UvecOext_\pm)(k,t) + \text{h.o.t.}
\end{aligned}
$$
respectively and for $n\in \N\setminus \{n_0\}$ by
$$
\begin{aligned}
&\Res_{n}(k,t)= \\
&\eps^{1/2}\left[(\omega_0-\omega_{n}(k))\Atil^\pm_{n}\left(\tfrac{k\mp k_0}{\eps},T\right) -M'^{(1)}_{n,n_0}(k)\Atil^\pm_{n_0}\left(\tfrac{k\mp k_0}{\eps},T\right) - M^{(2_\mp)}_{n,n_0}(k)\Atil^\mp_{n_0}\left(\tfrac{k\mp k_0}{\eps},T\right) \right]e^{-\ri \omega_0 t}\\
&+F'_{n}(\UvecOext_\pm,\UvecOext_\pm,\UvecOext_\pm)(k,t)+2F'_{n}(\UvecOext_\mp,\UvecOext_\mp,\UvecOext_\pm)(k,t)  + \text{h.o.t.}.
\end{aligned}
$$
Note that no $M^{R',m}(k)$-terms (with $m\in \Z_3^\pm$) appear because these are not supported near $\pm k_0$.

For $k\in (\kappa_j^\pm-\eps^{1/2},\kappa_j^\pm+\eps^{1/2})\cap\B_{2\pi}$ with $j\in \{1,\dots,J^\pm\}$ respectively we use the residual as given by the left hand side of \eqref{E:PNLS_U_rat} and get
$$
\begin{aligned}
&\Res_{n}(k,t)= \eps^{1/2}\left[(\omega_0-\omega_{n}(k))\Atil^\pm_{n,j}\left(\tfrac{k-\kappa_j^\pm}{\eps},T\right) - \sum_{\stackrel{m\in \Z_3^\pm}{\pm k_0+mk_W \in \kappa_j^\pm+\Z}}M'^{(R,m)}_{n,n_0}(k)\Atil^\pm_{n_0}\left(\tfrac{k-\kappa_j^\pm}{\eps},T\right) \right.\\
&\left. - \sum_{\stackrel{m\in \Z_3^\mp}{\mp k_0+mk_W \in \kappa_j^\pm+\Z}}M'^{(R,m)}_{n,n_0}(k)\Atil^\mp_{n_0}\left(\tfrac{k-\kappa_j^\pm}{\eps},T\right) \right]e^{-\ri \omega_0 t} + \text{h.o.t.}, \quad n \in \N.
\end{aligned}
$$
When $\kappa_j^+\in \pm 3k_0+\Z$ for some $\kappa_j^+\in \cK_R^+$ or $\kappa_j^-\in \pm 3k_0+\Z$ for some $\kappa_j^-\in \cK_R^-$, then also nonlinear terms $F'_n$ appear in this part of the residual. We treat, however, the neighborhoods of $\pm 3k_0$ separately below. Hence, all terms in the residual  are accounted for. 

Finally, we consider the residual for $k\in (\pm 3k_0-\eta-3\eps^{1/2},\pm 3k_0-\eta+3\eps^{1/2})\cap\B_{2\pi}$ with $\eta\in \cS_{\pm,\mp,\pm}$ respectively. Here
$$
\begin{aligned}
\Res_{n}(k,t)= &\eps^{1/2}e^{-\ri \omega_0 t}(\omega_0-\omega_{n}(k))\Atil^\pm_{n,NL}\left(\tfrac{k\mp 3k_0+\eta}{\eps},T\right) +F'_n(\UvecOext_\pm,\UvecOext_\mp,\UvecOext_\pm)(k,t)\\
& + \text{h.o.t.}, \quad n \in \N.
\end{aligned}
$$

In all other neighborhoods of its support the residual is of higher order in $\eps$, i.e. falls into the ``h.o.t.'' part. Similarly to \eqref{E:hot} the ``h.o.t.'' part consists of the following terms
\beq\label{E:hot_rat}
\begin{aligned}
&\eps \pa_T  \vec{U}^{1,\text{ext}}, \quad \eps M'^{(1)}\vec{U}^{1,\text{ext}}, \quad \eps \sum_{mk_W\notin \Z}M'^{(R,m)}\vec{U}^{1,\text{ext}}(\cdot-mk_W,t),\\
& \vec{F}'(\UvecOext,\vec{U}^{1,\text{ext}},\UvecOext)+2\vec{F}'(\UvecOext,\UvecOext,\vec{U}^{1,\text{ext}}),
\end{aligned}
\eeq
and nonlinear terms quadratic or cubic in $\vec{U}^{1,\text{ext}}$.

In analogy to \eqref{E:Atil_def},\eqref{E:Atilcor_def}, \eqref{E:Atiln_def}, and \eqref{E:Atilncor_def} we make the residual small by choosing
\beq\label{E:Atil_rat_def}
\Atil^\pm_{n_0}(K,T):=\chi_{[-\eps^{-1/2},\eps^{-1/2}]}(K)\Ahat_\pm(K,T),
\eeq
where $(A_+,A_-)(X,T)$ is a solution of \eqref{E:CME},
\beq\label{E:Atil_rat_def2}
\begin{aligned}
&\Atil^\pm_{n}(K,T)\\
&:=(\omega_0-\omega_n(\pm k_0+\eps K))^{-1}\left[M'^{(1)}_{n,n_0}(\pm k_0+\eps K)\Atil^\pm_{n_0}(K,T) + M^{(2_\mp)}_{n,n_0}(\pm k_0+\eps K)\Atil^\mp_{n_0}(K,T) \right.\\
& \left. - \eps^{-1/2}e^{\ri \omega_0 t}\left(F'_{n}(\UvecOext_\pm,\UvecOext_\pm,\UvecOext_\pm)(\pm k_0+\eps K,t)+2F'_{n}(\UvecOext_\mp,\UvecOext_\mp,\UvecOext_\pm)(\pm k_0+\eps K,t) \right)\right]
\end{aligned}
\eeq
for all $n \in \N \setminus \{n_0\}$,
\beq\label{E:Atilcor_rat_def}
\begin{aligned}
\Atil_{n,j}^\pm(K,T):=&(\omega_0-\omega_n(\kappa_j^\pm+\eps K))^{-1}\left[\Atil^\pm_{n_0}(K,T)\sum_{\stackrel{m\in \Z_3^\pm}{\pm k_0+mk_W \in \kappa_j^\pm+\Z}}M'^{(R,m)}_{n,n_0}(\pm \kappa_j^\pm+\eps K) \right.\\
&\left.+\Atil^\mp_{n_0}(K,T)\sum_{\stackrel{m\in \Z_3^\mp}{\mp k_0+mk_W \in \kappa_j^\pm+\Z}}M'^{(R,m)}_{n,n_0}(\pm \kappa_j^\pm+\eps K)
\right],
\end{aligned}
\eeq
and
\beq\label{E:Atilcor_ratNL_def}
\Atil_{n,NL}^\pm(K,T):=(\omega_0-\omega_n(\pm 3k_0+\eps K))^{-1}\eps^{-1/2}e^{\ri \omega_0 t}F'_n(\UvecOext_\mp,\UvecOext_\mp,\UvecOext_\pm)(\pm 3k_0+\eps K, t), \ n\in \N.
\eeq

The estimate of the residual and the Gronwall argument are completely analogous to the rational case and the proofs are omitted. Lemmas \ref{L:ext_est} and \ref{L:res_est} hold again with \eqref{E:Uext}, \eqref{E:Atil_def}, \eqref{E:Atilcor_def}, \eqref{E:Atiln_def}, and \eqref{E:Atilncor_def} replaced by \eqref{E:Uext_rat},
\eqref{E:Atil_rat_def}, \eqref{E:Atil_rat_def2}, \eqref{E:Atilcor_rat_def}, and \eqref{E:Atilcor_ratNL_def}. 

Also Lemma \ref{L:diff_app_ext} holds in the irrational case - with \eqref{E:uapp_rat} and \eqref{E:Uext} replaced by \eqref{E:uapp} and \eqref{E:Uext_rat} respectively. But some notational changes are needed in the proof. We list them next. For $u^\text{app}$ with $k_\pm=\pm k_0, k_0\in(0,1/2)$ and $p_\pm(x)=p_{n_0}(x,\pm k_0)$ we have
$$U^\text{app}_n(k,t)=\eps^{-1/2}e^{-\ri \omega_0 t} \sum_\pm \Ahat_\pm\left(\tfrac{k\mp k_0}{\eps},\eps t\right)\pi_n^\pm(k), \quad \text{where } \pi_n^\pm(k):=\langle p_\pm(\cdot), p_n(\cdot,k)\rangle_{2\pi}.$$
We decompose
$$\Uvecapp=\UvecOapp+\vec{U}^{1,\text{app}},$$ 
where
$$\UvecOapp(k,t):=e_{n_0}U^\text{app}_{n_0}(k,t).$$
Again, we need to show that $\|\vec{U}^{1,\text{app}}(\cdot,t)\|_{\cX(s)} \leq c\eps^{3/2}$ and $\|(\UvecOext-\UvecOapp)(\cdot,t)\|_{\cX(s)}\leq c\eps^{3/2}$ for all $t\in [0,\eps^{-1}T_0]$.

We have for any $n\neq n_0$
\begin{align}
\pi^\pm_n(k)&=\langle p_{n_0}(\cdot,k),p_n(\cdot,k)\rangle_{2\pi}+\langle p_{n_0}(\cdot,\pm k_0)-p_{n_0}(\cdot,k),p_n(\cdot,k)\rangle_{2\pi} \notag\\
&=\frac{1}{\omega_n(k)}\langle \cL(\cdot,k)(p_{n_0}(\cdot,\pm k_0)-p_{n_0}(\cdot,k)),p_n(\cdot,k)\rangle_{2\pi}\label{E:pi+-}
\end{align}
and by the $H^2$-Lipschitz continuity (in $k$) of the Bloch waves and using \eqref{E:band_as} (which holds also for $Q=2\pi$ with $\vartheta_n$ replaced by $\omega_n$)
$$|\pi^\pm_n(k)|\leq \frac{c}{n^2}|k\mp k_0| \quad \text{for all } n\in \N \ \text{and } k\in \B_{2\pi}.$$
As a result (for $s<3/2$)
$$
\begin{aligned}
&\|\vec{U}^{1,\text{app}}(\cdot,t)\|_{\cX(s)} = \left\|\left(\sum_{n\in \N\setminus\{n_0\}}n^{2s}|U^{1,\text{app}}_n(\cdot,t)|^2\right)^{1/2}\right\|_{L^1(\B_{2\pi})} \\
&\leq c\eps^{1/2}\left(\sum_{n\in \N}n^{2s-4}\right)^{1/2}\sum_\pm\int_{\B_{2\pi}}|\Ahat_\pm\left(\tfrac{k\mp k_0}{\eps},\eps t\right)||\tfrac{k\mp k_0}{\eps}|dk \leq c\eps^{3/2}\sum_\pm \|\Ahat_\pm(\cdot,\eps t)\|_{L^1_1(\R)}.
\end{aligned}
$$

For $\|(\UvecOext-\UvecOapp)(\cdot,t)\|_{\cX(s)}$  we have
$$
\begin{aligned}
\|(U^{0,\text{ext}}_{n_0}-& U^{0,\text{app}}_{n_0})(\cdot,t)\|_{L^1(\B_{2\pi})} \\
&\leq \eps^{-1/2}\sum_\pm\left(\int_{\pm k_0-\eps^{1/2}}^{\pm k_0+\eps^{1/2}}|\Ahat_\pm\left(\tfrac{k\mp k_0}{\eps},\eps t\right)||1-\pi_{n_0}^\pm(k)|dk\right.\\
& \left.  +\int_{\B_{2\pi}\setminus (\pm k_0-\eps^{1/2},\pm k_0+\eps^{1/2})}|\Ahat_\pm\left(\tfrac{k\mp k_0}{\eps},\eps t\right)||\pi_{n_0}^\pm(k)|dk\right).
\end{aligned}
$$
By the Lipschitz continuity it is $|1-\pi_{n_0}^\pm(k)|\leq L|k\mp k_0|$ and 
like in \eqref{E:epm_est} 
$$\|\Ahat_\pm(\cdot,\eps t)\|_{L^1(\R\setminus (-\eps^{-1/2}, \eps^{-1/2}))}\leq c\eps^{s_A/2}\|\Ahat_\pm(\cdot,\eps t)\|_{L^1_{s_A}(\R)}.$$ 
Hence
$$
\begin{aligned}
\|(\UvecOext-\UvecOapp)(\cdot,t)\|_{\cX(s)}&\leq c\eps^{3/2}\sum_\pm\|\Ahat_\pm(\cdot,\eps t)\|_{L^1_1(\R)} + c\eps^{1/2}\sum_\pm\|\Ahat_\pm(\cdot,\eps t)\|_{L^1(\R\setminus (-\eps^{-1/2}, \eps^{-1/2}))}\\
&\leq c\eps^{3/2}\sum_\pm\|\Ahat_\pm(\cdot,\eps t)\|_{L^1_1(\R)} + c\eps^{1/2+s_A/2}\sum_\pm\|\Ahat_\pm(\cdot,\eps t)\|_{L^1_{s_A}(\R)}.
\end{aligned}
$$

\section{Discussion}\label{S:discuss}

In fact, we have proved slightly more than the supremum norm estimate of the error. Namely, our proof estimates $\|u(\cdot,t)-u_\text{app}(\cdot,t)\|_{H^s(\R)}$ for all $s\in (1/2,3/2)$. In order to provide an estimate in $H^s$ for $s\geq 3/2$, the $l^2_s$-summability of $\Uvec^{1,\text{ext}}$ and $\Uvec^{1,\text{app}}$ has to be ensured. For this one would need in the rational case a sufficient decay of $M^{(1)}_{n,j}, M^{(R,m)}_{n,j}$ and $b^{(n)}_{\alpha,\beta,\gamma}$ as $n\to \infty$ instead of the simple estimates \eqref{E:MRm_est} and \eqref{E:Mb_est}. This can be achieved by replacing $q_n(x,k)$ in the inner products by $(\vartheta_n(k))^{-r}\cL^r(x,k)q_n(x,k)$ with $r\in \N$ sufficiently large and moving the self-adjoint $\cL^r$ to the other argument of the inner product. This would require sufficient smoothness of $q_{n_*}$ and $q_{n_*+1}$ in $x$. Similarly, $\cL$ would be applied more times in \eqref{E:pi+_L} to get a faster decay of $|\pi_n^+|$ (and similarly for $|\pi_n^-|$). This would require replacing $H^2$ by some smaller space $H^s, s >2$ in the Lipschitz assumption (H1). Analogous requirements would be needed in the irrational case on the respective matrices and vectors. Alternatively, the residual and the error can be estimated in $L^1(\B_Q,H^s(0,Q))$ instead of $\cX(s)$ in order to avoid the $l^2_s$-summability issue \cite{SU_book}. Because our result is sufficient to provide the physically relevant $C_b^0$ as well as $H^1$ bounds, we do not pursue the straightforward improvement to $s\geq 3/2$ here. 

Note also that the analysis can be easily carried over to other nonlinear equations with periodic coefficients. For instance, for the wave equation with periodic coefficients and the cubic nonlinearity $u^3$ as studied in \cite{BSTU06} the approach is completely analogous except for first rewriting the equation as a system of two first order equations. To obtain a real solution $u$ of the wave equation the ansatz is extended by adding the complex conjugate to \eqref{E:uapp}. Nevertheless, this does not change the analysis as no other $k$-points are generated by $u^3$  applied to such ansatz compared to $|u|^2u$ applied to \eqref{E:uapp}.


\section{Numerical Examples}\label{S:num}
We present numerical examples for both cases (a) and (b) from Section \ref{S:Bloch}: in Section \ref{S:num_a} for case (a) with simple Bloch eigenvalues at $k=k_+$ and $k=k_-=-k_+, k_+\in(0,1/2)$ and in Section \ref{S:num_b} for case (b) with a double eigenvalue at $k=k_+=k_-\in \{0,\pi/P\}$.

We solve \eqref{E:PNLS} by the Strang splitting method of second order in time, see e.g. \cite{WH86}. The equation is split into the part  $\ri\pa_t u=-\pa_x^2u$, which is solved with a spectral accuracy in Fourier space, and into the ODE part $\ri\pa_t u =(V(x)+\eps W(x))u+\sigma(x)|u|^2u$, which is solved exactly: $u(x,t)=e^{\ri(V(x)+\eps W(x)+\sigma(x)|u_0(x)|^2)t}u_0(x)$ (for initial data $u(x,0)=u_0(x)$). We discretize with $dx=0.05$ and $dt=0.02$ and solve \eqref{E:PNLS} up to $t=2\eps^{-1} $ with the initial data $u(x,0)=u_\text{app}(x,0)$ for a range of values of $\eps$ in order to study the error convergence. The choice of a solution $(A_+,A_-)$ of the CMEs is specified in each case below. 

\subsection{Case (a): simple eigenvalues at $k=k_\pm=\pm k_0, k_0\in (0,1/2)$}\label{S:num_a}
\begin{figure}[h!]
\begin{center}
\scalebox{0.5}{\includegraphics{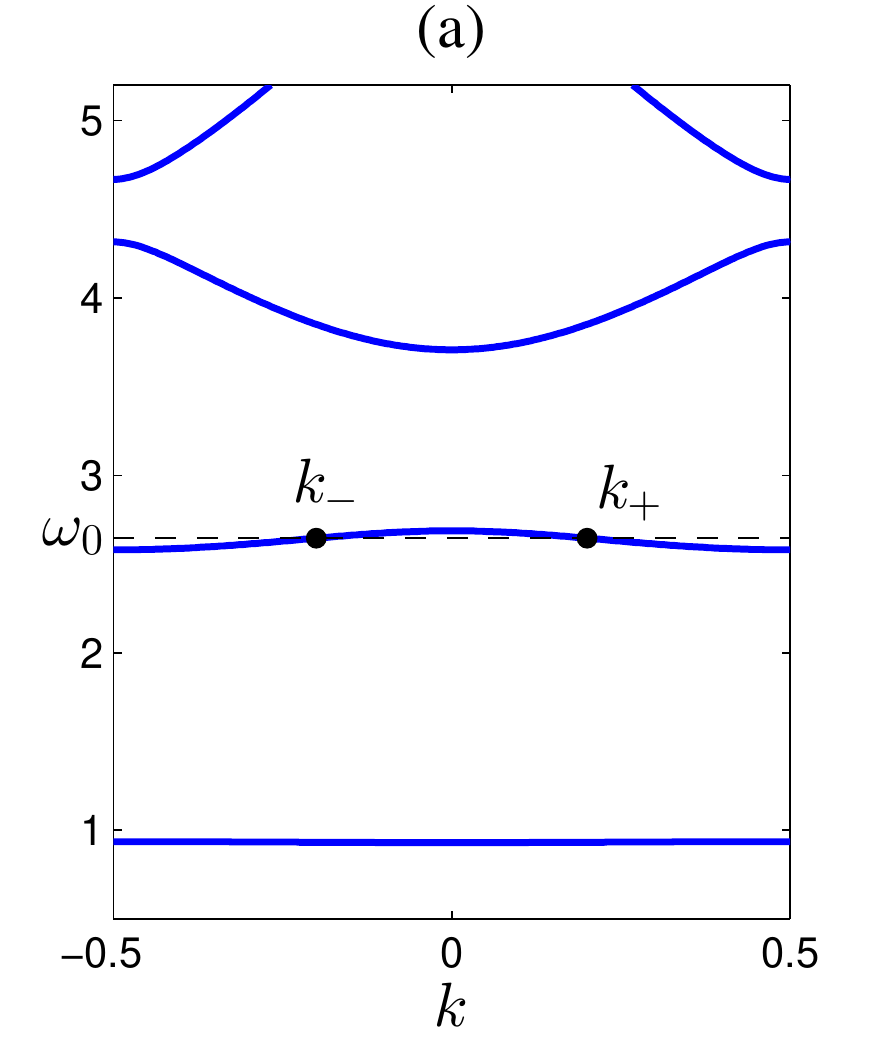}}
\scalebox{0.5}{\includegraphics{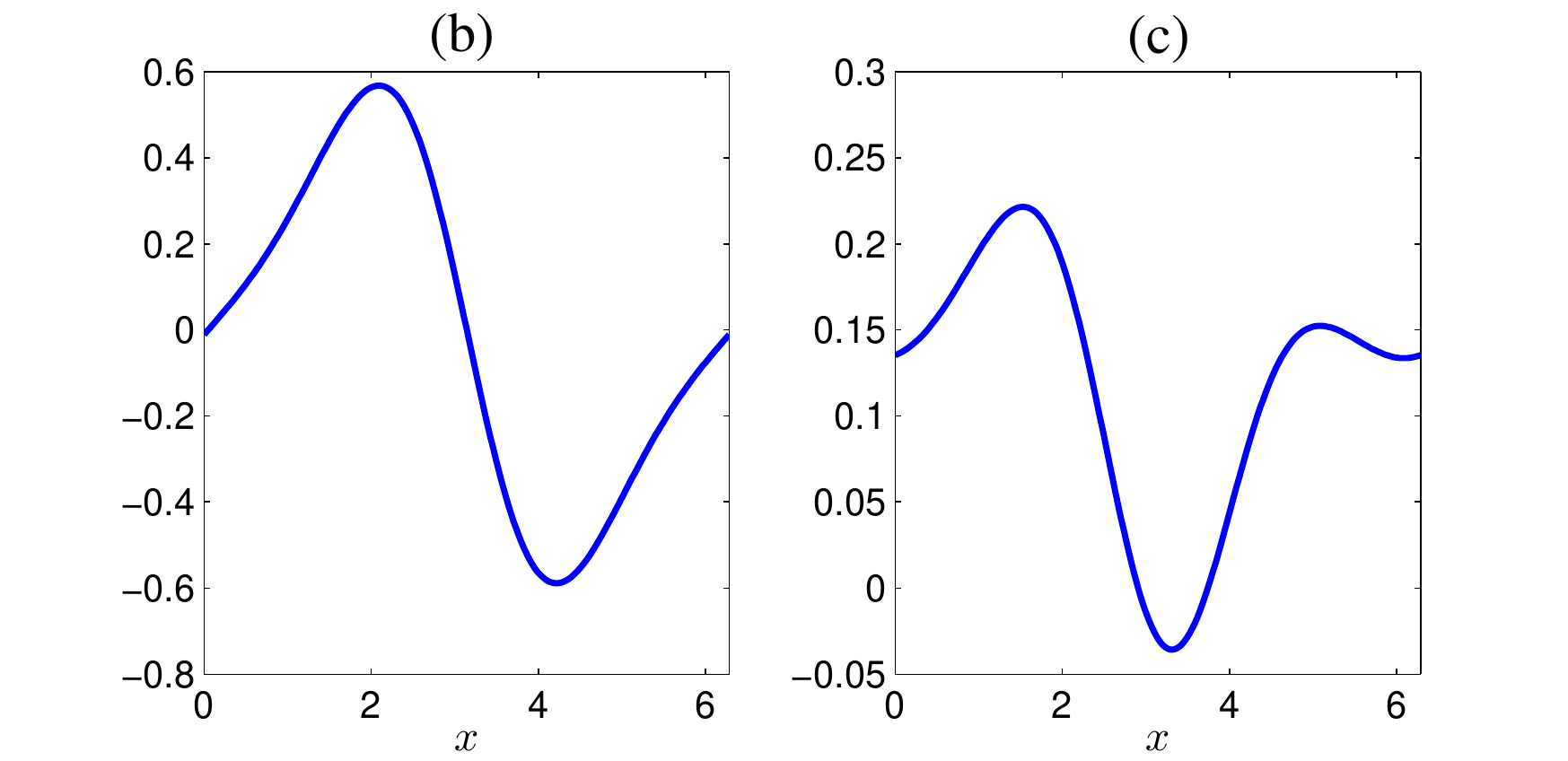}}
\end{center}
\caption{(a) band structure for $V(x) = 2(\cos(x)+1)$. The marked points are $(k_\pm,\omega_0)$ with $k_+=-k_-=0.2, \omega_0\approx 2.645$; (b) $\text{Re}(p_+)$; (c) $\text{Im}(p_+)$ where $p_+(x)=p_2(x,0.2)$.} \label{F:band_Bloch_cos}
\end{figure}

We choose here $V(x)=2(\cos(x)+1),\sigma \equiv -1$ such that $P=2\pi$. The band structure is plotted in Fig. \ref{F:band_Bloch_cos} (a). The carrier Bloch waves are given by the choice $k_0=0.2, n_0=2$ resulting in $\omega_0 \approx 2.645$. The Bloch function $p_+(x)=p_2(x,0.2)$ is plotted in Fig. \ref{F:band_Bloch_cos} (b) and (c).

Because $k_0\in(0,1/4)\cup(1/4,1/2)$, the resulting CMEs \eqref{E:CME} have $\beta=\gamma=0$ and hence reduce to the classical CMEs for envelopes of pulses in the nonlinear wave equation with an infinitesimal contrast periodicity \cite{AW89,GWH01}. There is the following two-parameter family of explicit solitary waves \cite{AW89,GWH01,D14}
\begin{equation}\label{E:1D_gap_sol}
\begin{split}
A_+(X,T)&= \nu a e^{\ri\eta}\sqrt{\frac{|\kappa|}{2|\alpha|}}\sin(\delta)\Delta^{-
1}e^{\ri\nu\zeta} \text{sech}(\theta - \ri\nu\delta/2),\\
A_-(X,T)&= -a e^{\ri\eta}\sqrt{\frac{|\kappa|}{2|\alpha|}}\sin(\delta)\Delta e^{\ri\nu\zeta} \text{sech}(\theta  +\ri\nu\delta/2),
\end{split}
\end{equation}
where
\begin{align*}
\nu &=\text{sign}(\kappa \alpha), \quad  a = \sqrt{\frac{2(1-v^2)}{3-v^2}}, \quad \Delta = \left(\frac{1-v}{1+v}\right)^{1/4}, \quad e^{\ri\eta} = \left(-\frac{e^{2\theta}+e^{-\ri\nu\delta}}{e^{2\theta}+e^{\ri\nu\delta}}
\right)^{\frac{2v}{3-v^2}},\\
\theta & = \mu \kappa
\sin(\delta)\left(\frac{X}{c_g}-vT\right), \quad \zeta = \mu
\kappa \cos(\delta)\left(\frac{v}{c_g}X-T\right), \quad \mu = (1-v^2)^{-1/2}
\end{align*}
with the velocity $v\in (-1,1)$ and ``detuning'' $\delta \in [0,\pi]$.

Next, we select two examples of the perturbations $\eps W$ of the periodic structure.
\subsubsection{\underline{$W(x)=\cos(2k_0 x)$}}\label{S:W_cos_2k0x}
\begin{figure}[h!]
\begin{center}
\scalebox{0.5}{\includegraphics{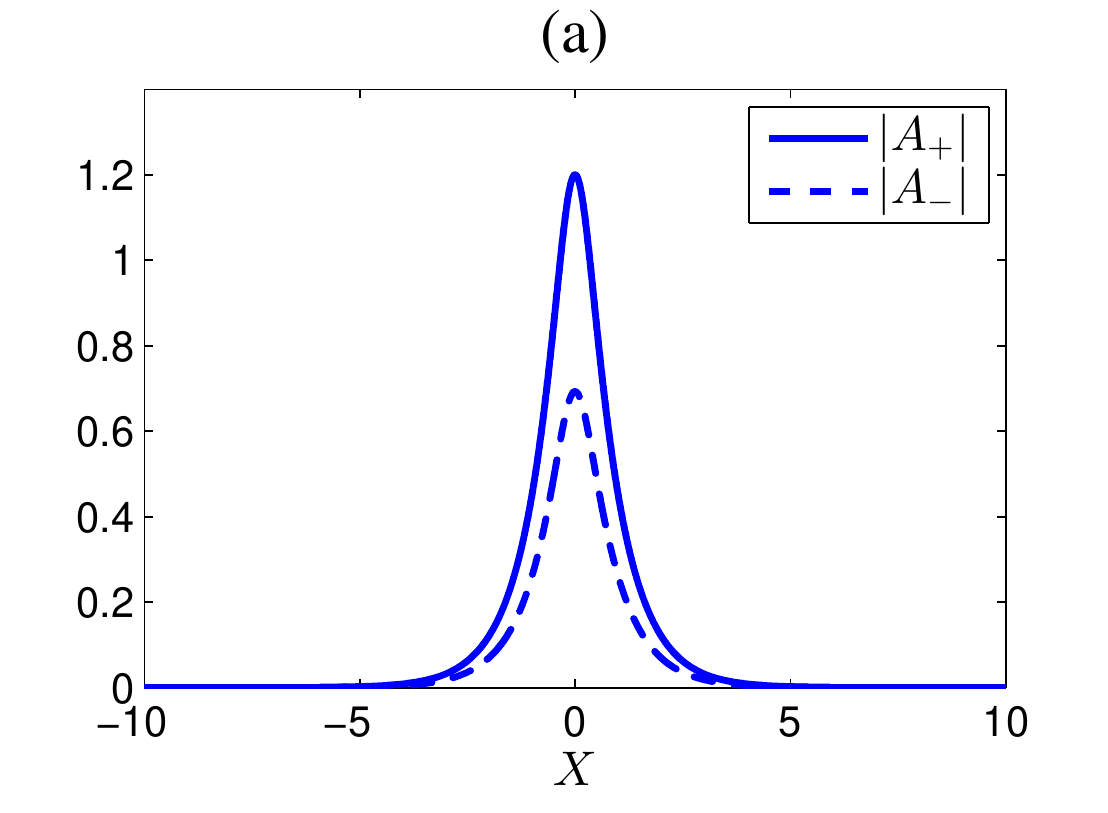}}
\scalebox{0.5}{\includegraphics{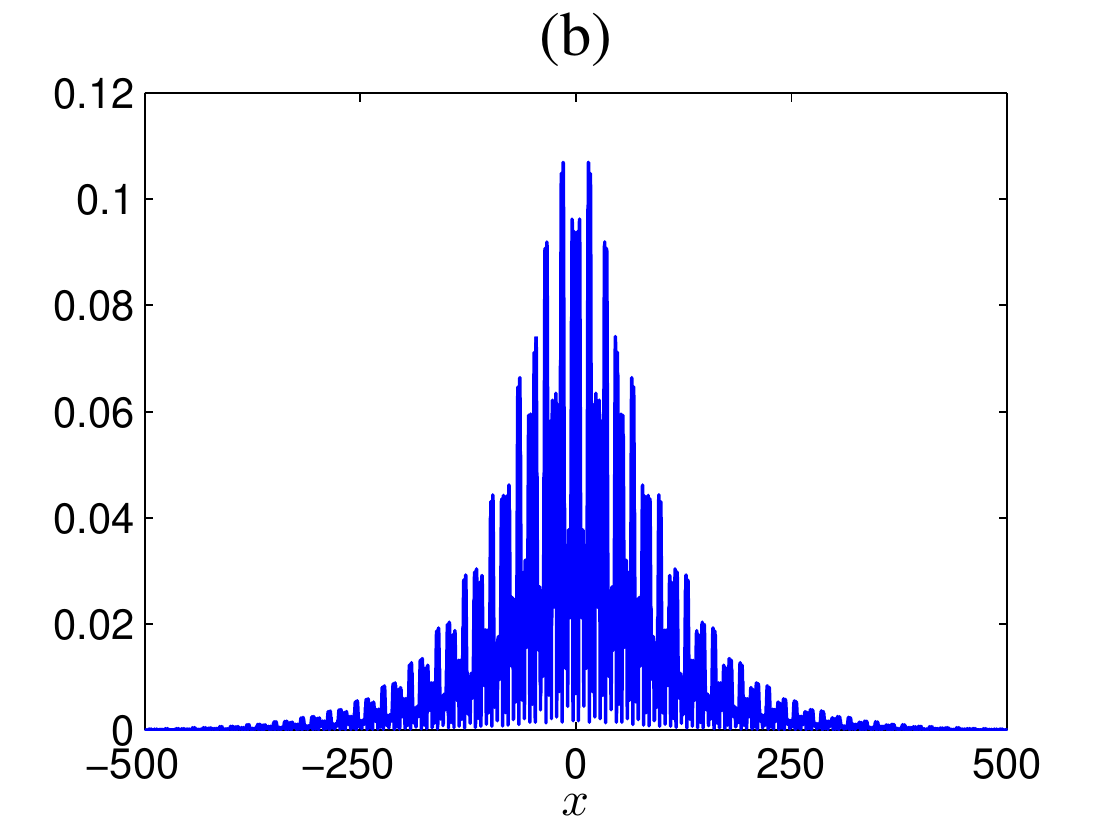}}
\end{center}
\caption{(a) $|A_\pm(X,0)|$ for the coefficients in Sec. \ref{S:W_cos_2k0x} and $\delta=\pi/2,v=0.5$; (b) $|u_\text{app}|$ with $A_\pm$ from (a) and with $\eps =0.01$.} \label{F:GS_profile}
\end{figure}

 For this $W$ the splitting in \eqref{E:W_split} is given by $W^{(1)}\equiv 0, W^{(2)}_\pm \equiv \tfrac{1}{2}, W^{(3)}_\pm(x)=\tfrac{1}{2}e^{\pm 2\ri k_0 x}$. The resulting coefficients of the CMEs are
\beq\label{E:CME_coeffs_ex}
\begin{aligned}
c_g \approx -0.3341, \quad \kappa\approx 0.3826, \quad \kappa_s=0, \quad \alpha\approx  0.2509, \quad \beta=\gamma=0.
\end{aligned}
\eeq
We choose the solution \eqref{E:1D_gap_sol} with $\delta =\pi/2$ and $v=0.5$. The modulus of $A_+$ and $A_-$ as well as of the approximation $u_\text{app}$ for $\eps =0.01$ are plotted in Fig. \ref{F:GS_profile}.

\begin{figure}[h!]
\begin{center}
\scalebox{0.5}{\includegraphics{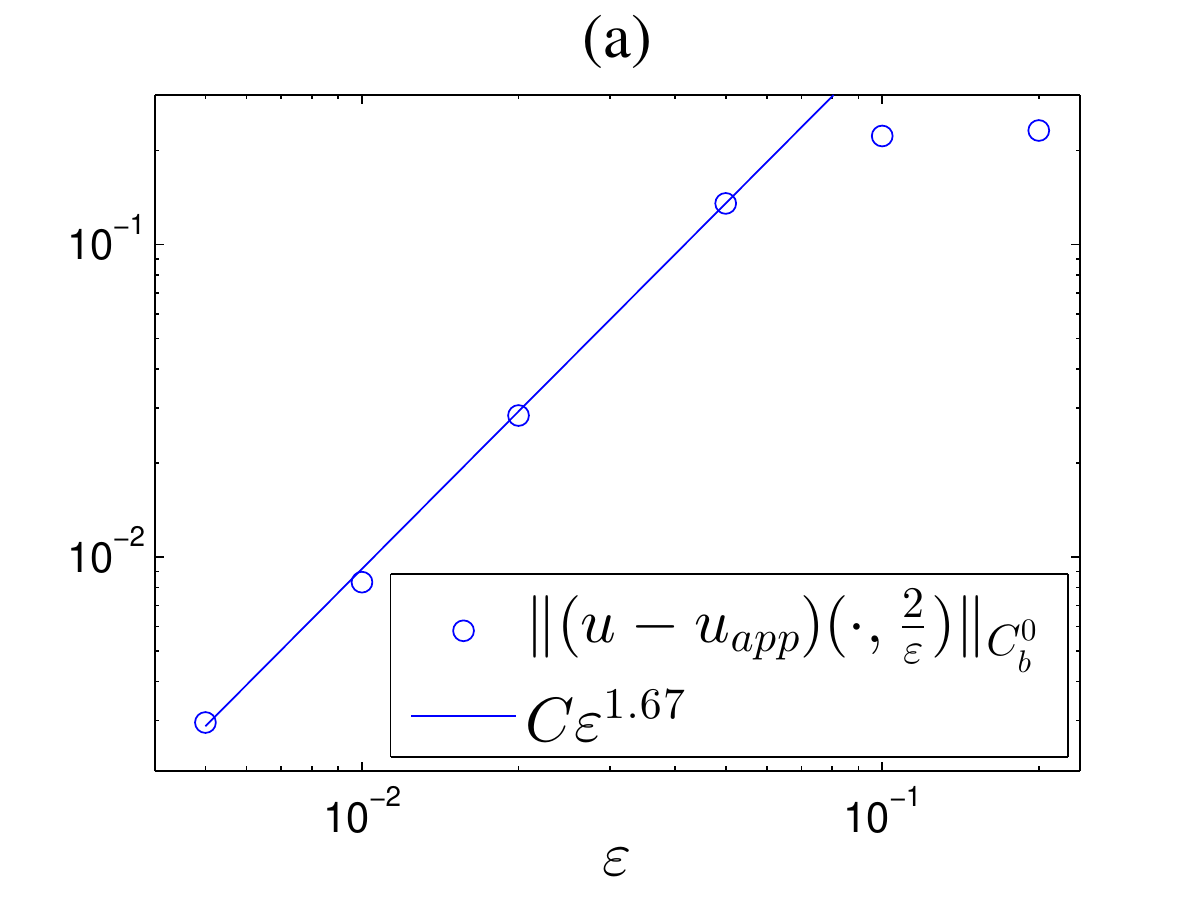}}
\scalebox{0.55}{\includegraphics{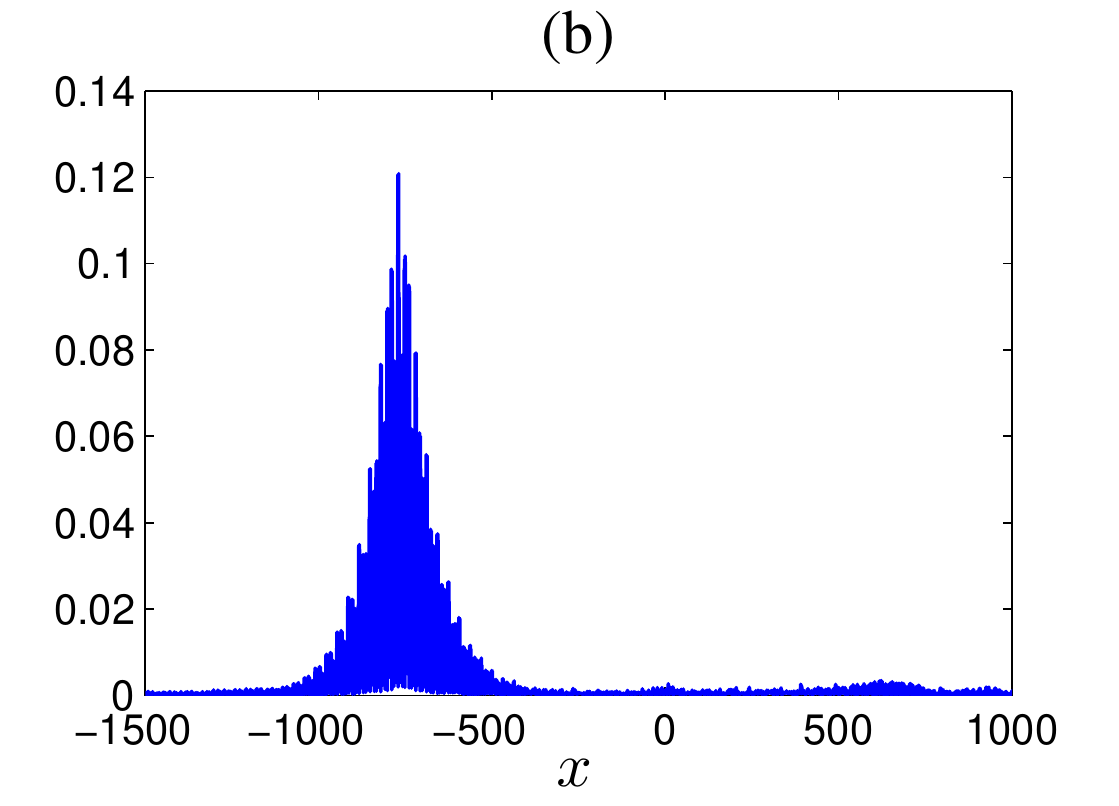}}
\end{center}
\caption{(a) Convergence of the supremum norm of the approximation error at $t=2\eps^{-1}$ with parameters from Sec. \ref{S:W_cos_2k0x}; (b) The modulus of the numerical solution $u(x,t)$ for $\eps=0.01$ at $t=5000$.} \label{F:eps_conv_Vcos}
\end{figure}

In Fig. \ref{F:eps_conv_Vcos} (a) we show that the error convergence in Theorem \ref{T:main} is confirmed by studying the error at $t=2\eps^{-1}$. The observed convergence rate is $\eps^{1.67}$. In order to demonstrate the approximation quality on a very large time interval Fig. \ref{F:eps_conv_Vcos} (b) shows the numerical solution $|u(x,t)|$ for $\eps=0.01$ and $u(x,0)=u_\text{app}(x,0)$ at $t=5000=50\eps^{-1}$. The solitary wave shape  is still well preserved.

\subsubsection{\underline{$W(x)=\cos(2k_0 x)+\tfrac{1}{2}\cos(4k_0 x)+\tfrac{1}{3}\cos(10 k_0 x)$}}\label{S:W_3terms}
\begin{figure}[h!]
\begin{center}
\scalebox{0.5}{\includegraphics{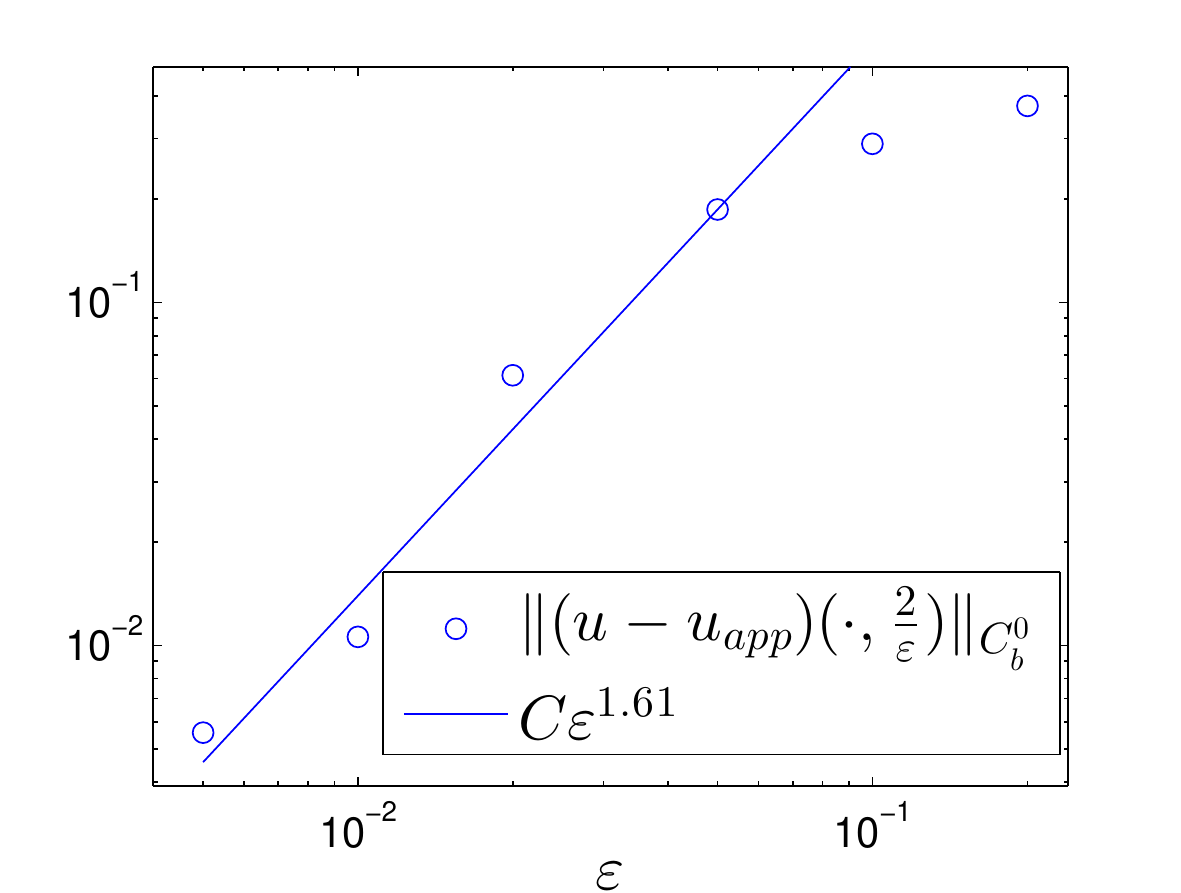}}
\end{center}
\caption{Convergence of the supremum norm of the approximation error at $t=2\eps^{-1}$ with parameters from Sec. \ref{S:W_3terms}.} \label{F:eps_conv_Vcos_W_3terms}
\end{figure}

For $k_0=0.2$ the splitting in \eqref{E:W_split} is $W^{(1)}(x)=\tfrac{1}{3}\cos(2x), W^{(2)}_\pm(x)=\tfrac{1}{2}, W^{(3)}_\pm(x)=\tfrac{1}{2}\cos(\tfrac{4}{5}x)$. The CME coefficients are given by \eqref{E:CME_coeffs_ex} except for $\kappa_s\approx 0.1324$, such that we can choose as $(A_+,A_-)$ the solutions of Sec. \ref{S:W_cos_2k0x} multiplied by $e^{\ri \kappa_s T}$. Also here the convergence rate is confirmed, see Fig. \ref{F:eps_conv_Vcos_W_3terms}, where the observed convergence rate is $\eps^{1.61}$.

\subsection{Case (b): double eigenvalue at $k=k_+=k_-\in \{0,1/2\}$}\label{S:num_b}

\begin{figure}[h!]
\begin{center}
\scalebox{0.5}{\includegraphics{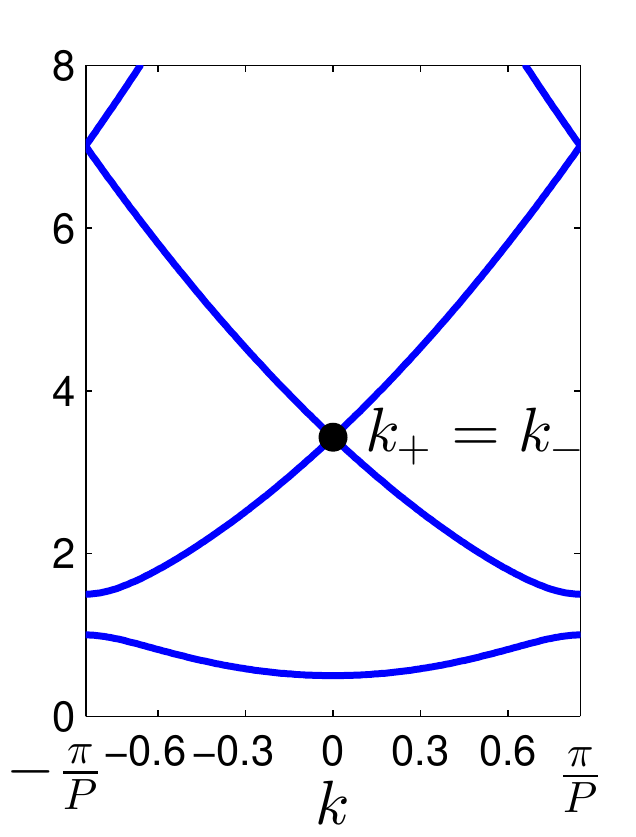}}
\end{center}
\caption{Band structure for the finite band potential $V(x) = \text{sn}^2(x;1/2)$ in Sec. \ref{S:num_b}. The marked point is $(k_\pm,\omega_0)$ with $k_+=k_-=0, \omega_0\approx 3.428$} \label{F:band_str_FB}
\end{figure}

\begin{figure}[h!]
\begin{center}
\scalebox{0.5}{\includegraphics{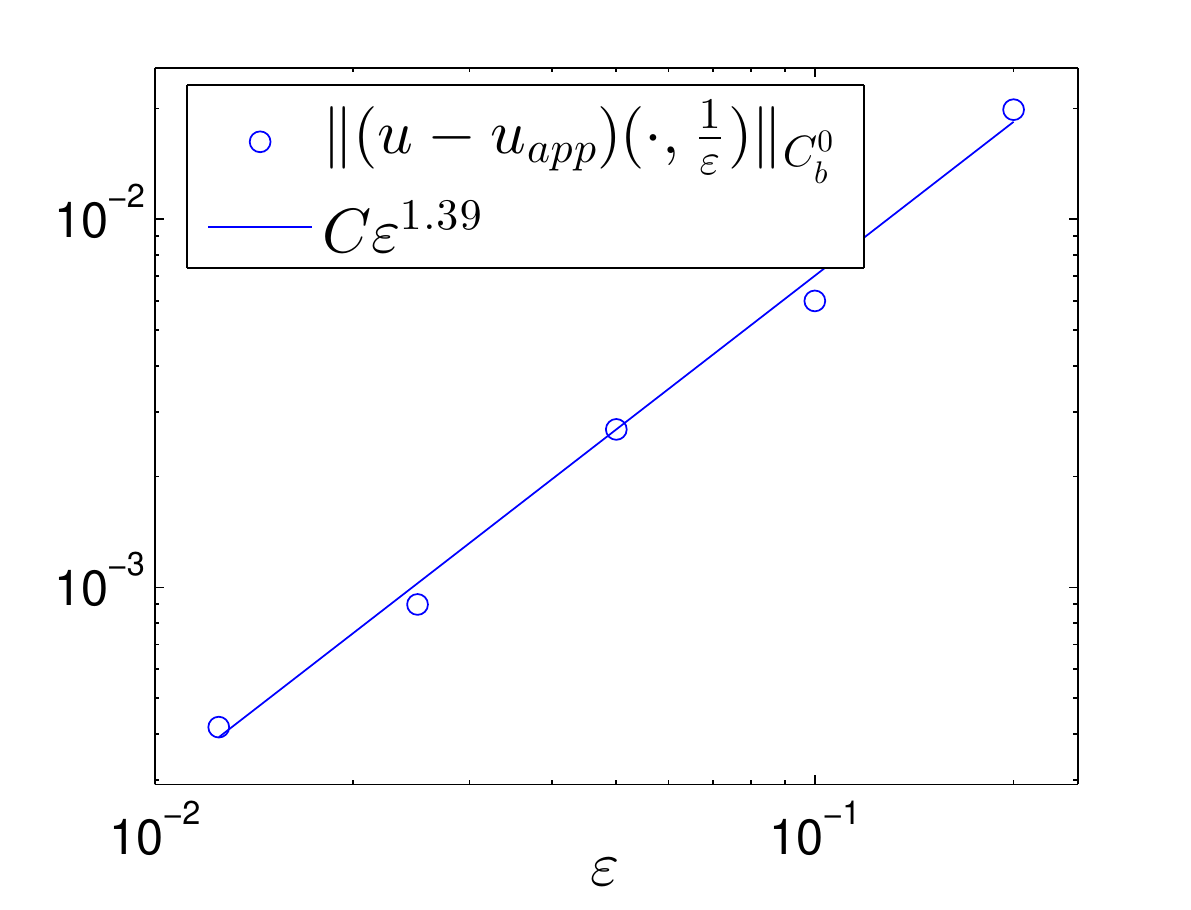}}
\end{center}
\caption{Convergence of the supremum norm of the approximation error at $t=\eps^{-1}$ with parameters from Sec. \ref{S:num_b}.} \label{F:err_conv_FB}
\end{figure}

This case was previously studied numerically in \cite{D14}. We choose here the example 6.2 of \cite{D14}, where $V$ is the finite band potential $V(x)=\text{sn}^2(x;1/2)$ and $\sigma \equiv -1$ such that $P\approx 3.7081$. The band structure from Fig. 2(a) in \cite{D14} is reproduced in Fig. \ref{F:band_str_FB}.
We choose the point $k_\pm=0, \omega_0\approx 3.428$. With $W(x)=\cos(\tfrac{4\pi}{P}x)$ we get a $P$-periodic perturbation of the potential $V$ such that setting $k_W=0$, we have in \eqref{E:W_split} $W^{(1)}=W$ after adjusting the definition of $W^{(1)},W^{(2)}_\pm,W^{(3)}_\pm$ to the $P-$periodic case (i.e. replacing $nk_W\in \Z$ by $nk_W\in\frac{2\pi}{P}\Z$). The CME coefficients are 
given in Sec. 6.2 of \cite{D14}. Although $\beta$ and $\gamma$ are very small: $\beta\approx 6.5*10^{-4}, \gamma\approx 7*10^{-6}$, it is shown that the error convergence is suboptimal when $\beta$ and $\gamma$ are set to zero, see Sec. 6.1 in \cite{D14}. Solutions of CMEs for $\beta\neq 0$ or $\gamma\neq 0$ can be found by a numerical parameter continuation technique starting from the explicit CME solutions \eqref{E:1D_gap_sol} at $\beta=\gamma=0$, see \cite{D14}. The error convergence for the velocity $v=0.5$ and detuning $\delta =\pi/2$ is plotted in Fig. \ref{F:err_conv_FB}. The observed rate $\eps^{1.39}$ is close to the predicted $\eps^{1.5}$. 

Note that in \cite{D14} the convergence  of the error was studied in the $L^2$-norm and the rate $\eps^{0.91}$, i.e. approximately $1/2$-smaller than in $\|\cdot\|_{C^0_b}$, was observed. Although our result does not guarantee the convergence in $L^2$, this is expected because heuristically, the error should be of the form $\eps^{3/2}B(\eps x,\eps t)f(x,t)$ for some bounded function $f$ and an $L^2$-function $B$. Due to the scaling of the $L^2$-norm, we then get an $O(\eps^1)$ estimate of the $L^2$-error.

\section*{Acknowledgments}
This research is supported by the \emph{German Research Foundation}, DFG grant No. DO1467/3-1. The authors thank Guido Schneider for fruitful discussions.

\bibliographystyle{plain}
\bibliography{biblio_CME1D}

\end{document}